%% file: main.tex
\documentclass{custom-arxiv-paper}

\usepackage[utf8]{inputenc}
\usepackage[T1]{fontenc}
\usepackage[british]{babel}
\usepackage{libertine}
\usepackage[libertine]{newtxmath}

\usepackage{jbmath3}
\usepackage{amsfonts}
\usepackage{siunitx}

\newcommand*{\grad}{\nabla}
\newcommand*{\Dp}{\mathcal{D}}
\newcommand*{\Vp}{\mathcal{V}}
\newcommand*{\Real}{\mathbb{R}}
\newcommand*{\Integer}{\mathbb{Z}}
\newcommand*{\Gw}{\Gamma_\text{w}}
\newcommand*{\Gp}{\Gamma_\text{p}}
\newcommand*{\Lp}{\mathcal{L}}
\newcommand*{\etol}{\varepsilon_\text{tol}}
\DeclareMathOperator*{\argmin}{arg~min}
\DeclareMathOperator{\sgn}{sgn}

\usepackage{placeins}
\usepackage{float}
\usepackage{booktabs}

\usepackage{xcolor}
\input{custom-colors}
\usepackage{graphicx}
\usepackage{grffile}
\usepackage{tikz}
\usepackage[framemethod=tikz]{mdframed}
\newmdenv[roundcorner=5pt, linecolor=myRed, linewidth=1.3pt,
          innerleftmargin=6pt, innerrightmargin=6pt,
          innertopmargin=6pt, innerbottommargin=6pt]{framedbox}

\newcommand*{\bibauthor}[1]{\textsc{#1}}
\newcommand*{\bibtitle}[1]{\textit{#1}}
\newcommand*{\bibarticle}[1]{#1}
\newcommand*{\bibvolume}[1]{\textbf{#1}}
\newcommand*{\bibissue}[1]{(#1)}
\newcommand*{\bibdoi}[1]{\href{https://doi.org/#1}{\textcolor{blue}{doi:#1}}}

\usepackage[hidelinks]{hyperref}
\title{Highly accurate special quadrature methods\\
for Stokesian particle suspensions in confined geometries}
\shorttitle{Highly accurate special quadrature methods for
Stokesian particle suspensions in confined geometries}
\author{Joar Bagge$^*$ \et Anna-Karin Tornberg$^*$}
\shortauthor{J.~Bagge and A.-K.~Tornberg}
\institution{Department of Mathematics, Linné Flow Centre/Swedish
e-Science Research Centre,\\
KTH Royal Institute of Technology, 100~44~Stockholm, Sweden}
\date{\today}

\begin{document}

\maketitle%
\renewcommand*{\thefootnote}{$*$}\footnotetext{\textit{E-mail
addresses:} \href{mailto:joarb@kth.se}{joarb@kth.se} (J.~Bagge),
\href{mailto:akto@kth.se}{akto@kth.se} (A.-K.~Tornberg).}%
\renewcommand*{\thefootnote}{\arabic{footnote}}%

\renewcommand*{\abstractname}{Abstract.\space}
\begin{abstract}
  Boundary integral methods are highly suited for problems with
  complicated geometries, but require special quadrature methods
  to accurately compute the singular and nearly singular layer
  potentials that appear in them. This paper presents a boundary
  integral method that can be used to study the motion of rigid
  particles in three-dimensional periodic Stokes flow with
  confining walls. A centrepiece of our method is the highly
  accurate special quadrature method, which is based on a
  combination of upsampled quadrature and quadrature by expansion
  (QBX), accelerated using a precomputation scheme. The method is
  demonstrated for rodlike and spheroidal particles, with the
  confining geometry given by a pipe or a pair of flat walls.
  A parameter selection strategy for the special quadrature
  method is presented and tested. Periodic interactions are
  computed using the Spectral Ewald (SE) fast summation method,
  which allows our method to run in $O(n \log n)$ time for $n$
  grid points, assuming the number of geometrical objects grows
  while the grid point concentration is kept fixed.

  \vspace{\baselineskip}\noindent
\end{abstract}

\chapter{Introduction}
\label{sec:introduction}

Microhydrodynamics is the study of fluid flow at low Reynolds
numbers, also known as Stokes flow or creeping flow. Applications
are found in biology, for example in the swimming of
micro\-organisms~\cite{guasto12} and in blood flow~\cite{lu19},
as well as in the field of microfluidics, which concerns the
design and construction of miniaturized fluid devices~\cite{squires05}.
Suspensions of rigid particles in Stokes flow are important both
in various applications and in fundamental fluid
mechanics~\cite{reddig13,gontijo15,mittal18,guasto10}. In
this paper, we describe a boundary integral method that can be
used to study the motion of rigid particles of different shapes
in Stokes flow. The particle suspension may also be confined in a
container geometry, such as a pipe or a pair of flat walls. The
flow in the fluid domain (i.e.\ within the container but outside
the particles) is governed by the \emph{Stokes equations}, which
for an incompressible Newtonian fluid take the form
\begin{align}
  \label{eq:stokes-1}
  \grad p - \mu \grad^2 \vec{u} &= \vec{f}, \\
  \label{eq:stokes-2}
  \grad \cdot \vec{u} &= 0.
\end{align}
Here, $p$ is the pressure, $\vec{u}$ is the flow velocity,
$\vec{f}$ is the body force per unit volume and
$\mu$ is the viscosity of the fluid. The Stokes equations arise
as a linearization of the Navier--Stokes equations in the case
where fluid inertia can be neglected, i.e.\ when the Reynolds
number is much less than 1.

On the surfaces of the particles and walls, no-slip boundary
conditions are prescribed. A problem of physical interest is the
\emph{resistance problem}: given the velocities of all particles,
compute the forces and torques (caused by viscous resistance)
acting on them by the fluid. The inverse problem is called the
\emph{mobility problem}: given the forces and torques acting on
all particles by the fluid, compute the particle velocities. The
mobility problem is useful in the case of noninertial particles,
since then the net force on each particle must be zero, so any
external forces (such as gravity) must be balanced by viscous
forces from the fluid; given external forces and torques, one can
then compute the motion of the particles.

Since the governing equations
\eqref{eq:stokes-1}--\eqref{eq:stokes-2} are linear, boundary
integral methods can be used to solve them. In these methods, the
flow is expressed in terms of layer potentials, which are integrals
over the boundary of the fluid domain (i.e.\ over the container
walls and particle surfaces). This reduces the dimensionality of
the problem from three to two, and leads to a smaller linear
system compared to
methods that must discretize the whole volume (such as the finite
difference or finite element methods). It it also easy to move
the particles, since no remeshing is needed. For a detailed
discussion on the properties of boundary integral methods, we
refer to the books by Pozrikidis \cite{pozrikidis92}, Atkinson
\cite{atkinson97} and Kress \cite{kress14}. Of special importance
are Fredholm integral equations of the second kind, which when
discretized, for example using the Nyström
method~\cite[ch.~4]{atkinson97}\cite[sec.~12.2]{kress14}, are
known to remain well-conditioned as the system size increases
\cite[p.~113]{atkinson97}\cite[p.~282]{kress14}.

The linear system resulting from the discretization of a
boundary integral equation is dense, and thus naive Gaussian
elimination would require $O(N^3)$ operations to solve a system
of $N$ unknowns. Using an iterative solution method such as
the generalized minimal residual method (GMRES)~\cite{saad86},
the complexity is reduced to $O(N^2)$ since the condition number
and thus the number of iterations are independent of the system
size (but may depend on the geometry). For a large system, this
complexity is still prohibitive. This can be overcome by using a
fast summation method such as the fast multipole method
(FMM)~\cite{greengard87,greengard97} or a fast Ewald summation
method~\cite{deserno98,lindbo11,klinteberg17b}, which reduce the
complexity further to $O(N)$ or $O(N \log N)$, respectively.

One of the challenges of boundary integral methods is the need
for accurate special quadrature methods for singular and nearly
singular integrands. These are necessary when evaluating the
layer potentials at a point on the boundary (where the integral
kernel is singular) or close to the boundary (where the kernel
is nearly singular, i.e.\ hard to resolve using a quadrature rule
designed for smooth integrands). Such special quadrature methods
are the main focus of this paper.

\section{Overview of related work}
\label{sec:intro-related-work}

In two dimensions, there are excellent special quadrature methods
available, such as the one introduced by Helsing and Ojala~\cite{helsing08},
which has been adapted to simulations of clean~\cite{ojala15} and
surfactant-covered~\cite{palsson19} drops in Stokes flow, as well
as vesicles~\cite{barnett15}. However, this method is based on a
complex variable formulation and not easy to generalize to three
dimensions.

In three dimensions, the development of an accurate and efficient
special quadrature method is still an active research problem,
especially in the nearly singular case. For an overview of methods
that have been used, we refer to \cite[sec.~1]{klockner13} and
\cite[sec.~1]{rahimian18}.
One of the most promising methods which is still under
development is \emph{quadrature by
expansion} (QBX), first introduced by Klöckner et
al.~\cite{klockner13} and Barnett~\cite{barnett14} and applied to
the Helmholtz equation in two dimensions. This method is based on
the observation that the layer potentials are smooth all the way
up to the boundary, and can therefore be locally expanded around a
point away from the boundary. The expansions can be evaluated at a
point closer to the boundary, or even on the boundary itself. The
convergence theory of QBX was developed in \cite{epstein13},
while \cite{klinteberg17} analyzed the error from the underlying
quadrature rule used to compute the expansion coefficients. A
strength of QBX is that it separates source points and target
points; the source points enter only in the computation of the
expansion coefficients, which can then be used to evaluate the
layer potential in all target points within a ball of
convergence. QBX has been applied to spheroidal particles in
three-dimensional Stokes flow by af Klinteberg and
Tornberg~\cite{klinteberg16b}, using a geometry-specific
precomputation scheme to accelerate the computation of the
coefficients.

A different approach that has been taken to accelerate QBX is to
couple it to a customized FMM, which has been done in two
dimensions \cite{rachh15,rachh17,wala18} and more recently in
three dimensions \cite{wala19,wala20}. This coupling is a natural
step to take since the FMM uses expansions of the same kind as
QBX, but it requires nontrivial modifications to the FMM. The
resulting method has complexity $O(N)$ and works for any smooth
geometry. The work published so far has been for the Laplace and
Helmholtz equations, but it is likely to be extended to more
kernels, including the ones needed for Stokes flow.

The QBX-FMM methods above all use \emph{global} QBX, in which
\emph{all} source
points are included in forming the local expansion. An
alternative is \emph{local} QBX, in which only source points that are
close to the expansion centre are included. Yet another variant
is found in \cite{klinteberg16b}, where all source points on a
single particle is used when forming expansions close to that
particle; we call this variant \emph{particle-global}. Local QBX is
typically combined with a patch-based discretization of the
geometry. While it reduces the cost of the method, it also poses
a challenge since the local layer potential from a single patch
may not be as smooth as the global layer potential from the
whole geometry (or a whole particle). Different versions of local
QBX have been described in two dimensions
\cite{barnett14,rachh15} and three dimensions \cite{siegel18}.
The latter paper also uses target-specific expansions, which need
only $O(p)$ terms to obtain the same accuracy as a QBX expansion
based on spherical harmonics with $O(p^2)$ terms. However, they
sacrifice the separation of source and target that is otherwise
present in QBX. This separation is in principle necessary in the
QBX-FMM methods, but also in these methods can target-specific
expansions be used to lower the computational cost of the method
\cite{wala20}.

Some of the recent work have focused on automating the parameter
selection based on a given error tolerance, resulting in the
\emph{adaptive} QBX method \cite{klinteberg18}. The results have
so far not been generalized to three dimensions. There has also been work on a
kernel-independent version of QBX, called \emph{quadrature by
kernel-independent expansion} (QBKIX) \cite{rahimian18} and meant
to be combined with the kernel-independent FMM. The published
work is in two dimensions, but a generalization to three
dimensions is expected to follow.

Other methods, which are not based on QBX, have also been used
successfully as special quadrature methods in three dimensions.
One example is the ``line interpolation method'' introduced by
Ying et al.\ \cite{ying06}. In this method, a line is constructed
through the target point, which is close to the boundary, and its
projection onto the boundary. The layer potential is evaluated at
points further away from the boundary along this line, and also
at the projection point where the line intersects the boundary if a separate
singular integration method is available. The value at the target
point is then computed using interpolation along the line (or
extrapolation if no singular integration method is available).
Like QBX, the success of this method hinges on the fact that they
layer potential is smooth in the domain, so that it can be well
interpolated (or extrapolated). It has been applied to
surfactant-covered drops~\cite{sorgentone18} and
vesicles~\cite{malhotra18} in three-dimensional Stokes flow. The
extrapolatory method used in \cite{lu19} falls into the same
category. Other types of methods are based on
regularizing the kernel and adding corrections
\cite{beale01,tlupova13,beale16,tlupova19}, density interpolation
techniques \cite{perez19}, coordinate rotations and a subtraction
method \cite{carvalho18b}, asymptotic approximations
\cite{carvalho18c}, analytical expressions available only for spheres
\cite{corona18} or floating partitions of unity
\cite{bruno01,zhao10,klinteberg14}. Many of these methods are
target-specific, and their cost grows rapidly if there are many
nearly singular target points.

\section{Scope of this paper}

In this paper, we present a boundary integral method based on the
Stokes double layer potential, which can be used to solve the
mobility and resistance problems for a system of rigid particles
in incompressible three-dimensional Stokes flow, possibly
confined within a container geometry. Our formulation leads to a
Fredholm integral equation of the second kind. We use QBX for singular integration,
and a combination of QBX and upsampled quadrature for nearly singular
integration. Our QBX implementation is based on the work by af
Klinteberg and Tornberg \cite{klinteberg16b}, which we have
extended to rodlike particles, plane walls and pipes (using
particle-global QBX for the particles and local QBX for the two
wall geometries). A precomputation scheme is used to greatly
accelerate the QBX computations for all geometries. For this
precomputation scheme to be feasible, we require that each
particle or wall is rigid and has some degree of symmetry, such
as axisymmetry or reflective symmetry. Nonetheless, we have
chosen this route since the implementation is relatively simple
compared to e.g.\ a QBX-FMM method. When container walls are
present, we restrict ourselves to periodicity in all three
spatial directions and use a fast Ewald summation method called
the Spectral Ewald method \cite{lindbo10,klinteberg14,klinteberg16a}
to accelerate computations.\footnote{
  The implementation of the Spectral Ewald method that we use is
  publicly available at \cite{ewald-package}.
} In this situation our method scales
as $O(N \log N)$ in the number of unknowns $N$, assuming fixed
grid point concentration. The container geometry is
restricted to a periodic straight pipe or a pair of periodic
plane walls. Our contributions include:
\begin{itemize}
  \item
    The combined special quadrature method based on QBX and
    upsampling, which we have implemented for spheroidal and
    rodlike particles, plane walls and pipes. (The QBX
    implementation for spheroids is reused from
    \cite{klinteberg16b}. Our initial work on QBX for plane walls
    is published in the conference proceedings \cite{bagge17}.)
  \item
    A strategy for experimentally selecting the parameters of the
    special quadrature method to meet a given error tolerance. We
    also demonstrate that the boundary integral method in full
    meets the given error tolerance and scales as $O(N \log N)$.
  \item
    Construction of fully smooth rodlike particles. We
    demonstrate the effect of smoothness on the convergence of
    the local expansions in this particular case.
  \item
    Derivation of a stresslet identity for an infinite pipe and a
    pair of infinite plane walls. This is used as an exact
    solution to test the special quadrature method.
  \item
    An outline of how streamlines can be efficiently computed for
    periodic problems using the Spectral Ewald method, by reusing
    data. (This idea was used, but not explicitly described, in
    \cite{klinteberg16b}.)
  \item
    The so-called completion sources that appear in our
    formulation are distributed along the axis of symmetry of
    rodlike particles, and we have studied how the
    number of completion sources influences the accuracy.
\end{itemize}

\section{Organization of the paper}

In section~\ref{sec:formulation}, we introduce the mathematical
formulation of the problem, including the boundary integral
formulation and the boundary integral equations for the
resistance and mobility problems. In
section~\ref{sec:disc-quadrature}, we describe the discretization
of the geometry and the quadrature method, including the combined
special quadrature. The details on our QBX method are then given
in section~\ref{sec:qbx-detail}, including the precomputation
scheme. In section~\ref{sec:periodicity}, we describe how
periodicity is treated and how the special quadrature is combined
with the Spectral Ewald method. Then, in
section~\ref{sec:parameters}, our parameter selection strategy
for the special quadrature is described and demonstrated.
Numerical results are given in section~\ref{sec:results}, to
demonstrate the accuracy and scaling of our method. Finally, in
section~\ref{sec:smoothness}, we demonstrate the effect of
nonsmooth geometries on the convergence. The appendices include a
derivation of the stresslet identity for plane walls and pipes,
details on the construction of the smooth rodlike particles, and
a note on streamline computation.

\FloatBarrier
\chapter{Mathematical formulation}
\label{sec:formulation}

\begin{figure}[h!]
  \centering\small
  \begin{minipage}[b]{0.5\textwidth}%
    \centering
    \includegraphics{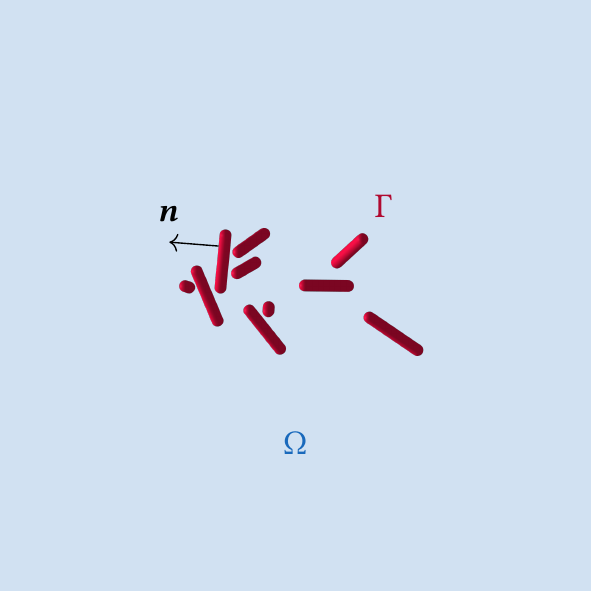}\\
    (a) Free-space problem
  \end{minipage}%
  \begin{minipage}[b]{0.5\textwidth}%
    \centering
    \includegraphics{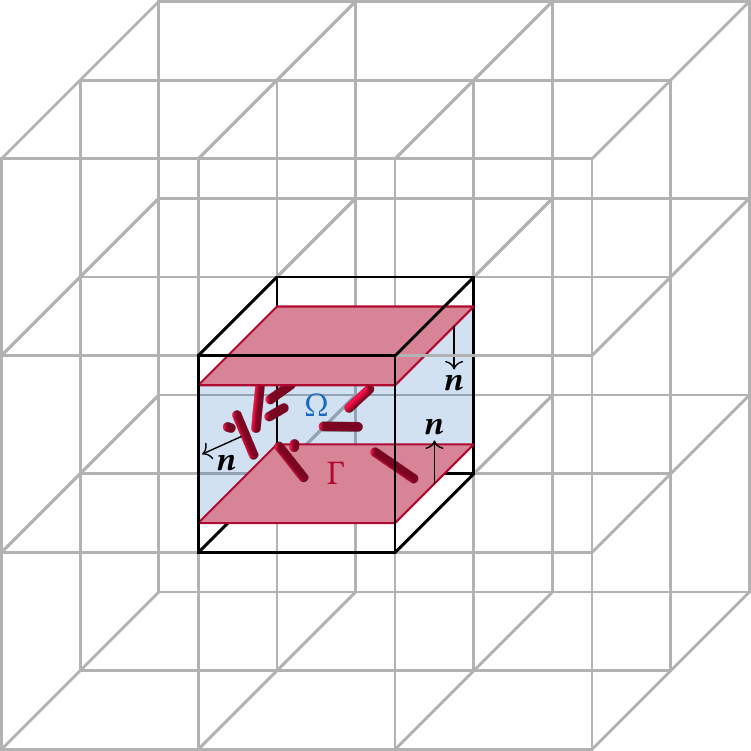}\\
    (b) Fully periodic problem
  \end{minipage}%
  \caption{The geometry for (a) a free-space problem, and (b) a fully
  periodic problem. In (b), the primary cell is marked with a
  darker outline than the other cells. The geometry is shown
  only in the primary cell. The lattice of periodic cells fills
  the whole space; only a small part is shown here.}
  \label{fig:periodicity}
\end{figure}

We consider two different kinds of problems: free-space problems
and fully periodic problems. In a \emph{free-space problem},
$M$~particles (spheroids or rods) are located in a fluid
extending to infinity. We denote the fluid domain by $\Omega$ and
its boundary, i.e.\ the union of all particle surfaces, by
$\Gamma$. The Stokes equations
\eqref{eq:stokes-1}--\eqref{eq:stokes-2} with $\vec{f}=\vec{0}$
hold in $\Omega$, while
no-slip boundary conditions are imposed on $\Gamma$. The unit
normal vector $\vec{n}$ of $\Gamma$ is defined to point into the
fluid domain $\Omega$, as shown in Figure~\ref{fig:periodicity}~(a).

A \emph{fully periodic problem}, on the other hand, is periodic
in all three spatial directions. The \emph{primary cell} is a box
with side lengths $\vec{B} = (B_1, B_2, B_3)$, which is
considered to be replicated periodically in all three spatial
directions, as illustrated by Figure~\ref{fig:periodicity}~(b).
Let the number of particles in the primary cell be $M$. In this
case we also allow a container consisting either of a pair of
plane walls or a pipe. Only the geometry inside the primary cell
is discretized, which means that $\Gamma$ consists of the union
of the $M$~particle surfaces and the parts of the wall surfaces
that lie in the primary cell.\footnote{
  For a single infinitely large plane wall, the method of images
  can be used \cite{blake71,gimbutas15,srinivasan18}, which has
  the advantage that the wall itself does not need to be
  discretized. However, that method does not work when there are
  more than one wall, or when the wall is curved, which are the
  cases we consider here. Therefore we must discretize the walls.
}
The fluid domain $\Omega$ lies within the container but outside
the particles; the flow is thus external to the particles, but
internal to the surrounding walls. The unit normal vector
$\vec{n}$ of $\Gamma$ is always defined to point into $\Omega$.

Below, we introduce our boundary integral formulation in the
free-space setting, or for the primary cell without periodicity.
Full treatment of the periodic problem is deferred to
section~\ref{sec:periodicity}.

\section{Boundary integral formulation}
\label{sec:bif}

Any flow field $\vec{u}$ that satisfies the Stokes equations
\eqref{eq:stokes-1}--\eqref{eq:stokes-2} with $\vec{f}=\vec{0}$
can be expressed in terms of
integrals over the boundary of the fluid domain $\Omega$, as
described for example by Pozrikidis \cite[ch.~4]{pozrikidis92}
and Kim and Karrila \cite[ch.~14--16]{kimkarrila91}. The boundary
integral formulation that we use is based on the
\emph{Stokes double layer potential}~$\vec{\Dp}$, which in free space is given by%
\footnote{The Einstein summation convention is used in this
paper, meaning that indices appearing twice in the same term are
to be summed over the set $\{1,2,3\}$. The remaining free indices
take values in the same set.}
\begin{equation}
  \label{eq:double-layer}
  \Dp_i[\Gamma,\vec{q}](\vec{x}) = \int_\Gamma
  T_{ijk}(\vec{x}-\vec{y}) q_j(\vec{y}) n_k(\vec{y}) \, \D
  S(\vec{y}), \qquad
  T_{ijk}(\vec{r}) = -6\frac{r_i r_j r_k}{\absi{\vec{r}}^5}.
\end{equation}
Here, $\Gamma$ and $\vec{n}$ are as in Figure~\ref{fig:periodicity},
and the \emph{double layer density} $\vec{q}$ is a continuous
vector field defined on $\Gamma$. The tensor kernel $\mat{T}$
in \eqref{eq:double-layer} is known as the \emph{stresslet}. The
potential $\vec{\Dp}$ has a jump discontinuity as $\vec{x}$
passes over $\Gamma$. More specifically, for $\vec{x} \in
\Gamma$ it holds that \cite[p.~110]{pozrikidis92}
\begin{equation}
  \label{eq:double-jump}
  \lim_{\varepsilon \to 0^{+}} \vec{\Dp}[\Gamma, \vec{q}](\vec{x}
  \pm \varepsilon \vec{n}) = \vec{\Dp}[\Gamma, \vec{q}](\vec{x})
  \mp 4\pi \vec{q}(\vec{x}).
\end{equation}
For any closed Lyapunov surface $\tilde{\Gamma}
\subseteq \Gamma$ and any constant vector $\tilde{\vec{q}}$, the
\emph{stresslet identity} \cite[p.~28]{pozrikidis92}
\begin{equation}
  \label{eq:stresslet-identity}
  \vec{\Dp}[\tilde{\Gamma}, \tilde{\vec{q}}](\vec{x}) =
  \begin{cases}
    \vec{0}, & \text{if $\vec{x}$ is outside the domain enclosed
    by $\tilde{\Gamma}$}, \\
    4\pi \tilde{\vec{q}}, & \text{if $\vec{x} \in \tilde{\Gamma}$}, \\
    8\pi \tilde{\vec{q}}, & \text{if $\vec{x}$ is inside the
    domain enclosed by $\tilde{\Gamma}$},
  \end{cases}
\end{equation}
holds. We will use this identity as a test case for the special
quadrature method in sections~\ref{sec:parameters} and
\ref{sec:res1-special-quadrature}. In
appendix~\ref{app:stresslet-identity} we show that a variant
of \eqref{eq:stresslet-identity} holds also for the wall
geometries that we consider, despite them not being closed
surfaces.

The double layer potential $\vec{u}(\vec{x}) = \vec{\Dp}[\Gamma,
\vec{q}](\vec{x})$ is a solution to \eqref{eq:stokes-1}--\eqref{eq:stokes-2}
in $\Omega$ for any continuous vector field $\vec{q}$. However,
not every solution to \eqref{eq:stokes-1}--\eqref{eq:stokes-2}
can be represented by a double layer potential alone; for
instance, as noted in
\cite{power-miranda}\cite[p.~119]{pozrikidis92}, the force and
torque exerted on any particle by the flow from a double layer
potential will always be zero, whereas a Stokes flow in general
can exert a nonzero force and torque on the particles (which is a
central point of the resistance and mobility problems mentioned in
section~\ref{sec:introduction}). This is related to the presence
of a nontrivial nullspace of the operator $\vec{q} \mapsto
\vec{\Dp}[\Gamma,\vec{q}]$ for external flows, which can be
immediately seen from the stresslet identity
\eqref{eq:stresslet-identity}.

To remove the nontrivial nullspace and allow for nonzero forces
and torques on the particles, we add a \emph{completion flow}
$\vec{\Vp}$, first introduced by Power and Miranda \cite{power-miranda}.
The completion flow is also a solution to
\eqref{eq:stokes-1}--\eqref{eq:stokes-2} and can be identified as
the flow from a point force $\vec{F}$ and a point torque
$\vec{\tau}$ located at $\vec{y}$. It is given by
\begin{equation}
  \label{eq:completion}
  \Vp_i[\vec{F}, \vec{\tau}, \vec{y}](\vec{x}) =
  \frac{1}{8\pi\mu} \gp{
    S_{ij}(\vec{x}-\vec{y}) F_j +
    R_{ij}(\vec{x}-\vec{y}) \tau_j
  }, \qquad \vec{x} \in \Omega,
\end{equation}
where the \emph{stokeslet} $\mat{S}$ and the \emph{rotlet}
$\mat{R}$ are given by%
\footnote{Here, $\delta_{ij}$ denotes the Kronecker delta, and
$\epsilon_{ijk}$ is the Levi-Civita symbol.}
\begin{equation}
  S_{ij}(\vec{r}) = \frac{\delta_{ij}}{\absi{\vec{r}}} +
  \frac{r_i r_j}{\absi{\vec{r}}^3} \qquad \text{and} \qquad
  R_{ij}(\vec{r}) = \epsilon_{ijk} \frac{r_k}{\absi{\vec{r}}^3},
\end{equation}
respectively. We call a pair $(\vec{F}, \vec{\tau})$ a
\emph{completion source}. Such completion sources are placed in the
interior of every particle. Mathematically, one completion
source per particle is sufficient, but this may lead to numerical
problems in some cases. In this paper, we allow for multiple
completion sources to be distributed along a line segment within
the particle; as we show in section~\ref{sec:res2-convergence-Nsrc},
this is important for elongated particles.
For the particle with index~$\alpha$, let
$\vec{F}^{(\alpha)}$ and $\vec{\tau}^{(\alpha)}$ be the net force
and torque, respectively, exerted on the fluid by the particle,
and let $\vec{y}_\text{c}^{(\alpha)}$ be the centre of mass of
the particle. (For a noninertial particle, $\vec{F}^{(\alpha)}$ and
$\vec{\tau}^{(\alpha)}$ would be equal to the net external force
and torque, respectively, acting on the particle.) The completion
flow associated with particle~$\alpha$ is then given by
\begin{equation}
  \label{eq:completion-particle}
  \vec{\Vp}^{(\alpha)}[\vec{F}^{(\alpha)}, \vec{\tau}^{(\alpha)}](\vec{x})
  = \frac{1}{N_\text{src}} \sum_{s=1}^{N_\text{src}}
  \vec{\Vp}[\vec{F}^{(\alpha)}, \vec{\tau}^{(\alpha)},
  \vec{y}_\text{c}^{(\alpha)} + C(s,N_\text{src})
  \vec{a}^{(\alpha)}](\vec{x}),
\end{equation}
where $N_\text{src}$ is the number of completion sources per
particle, $\vec{\Vp}$ is given by \eqref{eq:completion}, and
$\vec{a}^{(\alpha)}$ is a vector specifying the line segment
along which completion sources are placed. The function $C$ is
given by
\begin{equation}
  C(s, N_\text{src}) = \begin{cases}
    0, & \text{if $N_\text{src} = 1$}, \\
    -1 + 2\dfrac{s-1}{N_\text{src}-1}, & \text{if $N_\text{src} > 1$}.
  \end{cases}
\end{equation}

Both the double layer potential $\vec{\Dp}$ and the completion
flow $\vec{\Vp}$ have the property that they decay to zero as
$\vec{x} \to \infty$. To be able to represent flows which do not
decay, we add a \emph{background flow}
$\vec{u}_\text{bg}$, which is a known solution to
\eqref{eq:stokes-1}--\eqref{eq:stokes-2} in the whole physical space,
ignoring all particles and walls. The total flow~$\vec{u}$ in the
presence of particles and walls is thus written as
\begin{equation}
  \label{eq:total-flow}
  \vec{u}(\vec{x}) = \vec{u}_\text{bg}(\vec{x}) +
  \vec{u}_\text{d}(\vec{x}),
\end{equation}
where $\vec{u}_\text{d}$ is a \emph{disturbance flow} which is
responsible for enforcing the no-slip boundary conditions on the
solid boundary $\Gamma$. As $\vec{x}$ moves away from $\Gamma$,
the disturbance flow $\vec{u}_\text{d}$ should decay to zero, and
the total flow should therefore approach the background flow
$\vec{u}_\text{bg}$. The disturbance flow is written as
\begin{equation}
  \label{eq:dist-flow}
  \vec{u}_\text{d}(\vec{x}) = \vec{\Dp}[\Gamma,\vec{q}](\vec{x})
  + \sum_{\alpha=1}^M \vec{\Vp}^{(\alpha)}[\vec{F}^{(\alpha)},
  \vec{\tau}^{(\alpha)}](\vec{x}),
\end{equation}
where $\vec{\Vp}^{(\alpha)}$ is as in \eqref{eq:completion-particle},
and the double layer density $\vec{q}$ must be determined through
the boundary conditions. Note that $\vec{u}_\text{d}$ as given by
\eqref{eq:dist-flow} decays as $\vec{x} \to \infty$, and by the
superposition principle it satisfies \eqref{eq:stokes-1}--\eqref{eq:stokes-2}.
Also note that completion sources are
placed inside the particles since the flow is external to the
particles, but not inside the walls since the flow is internal to
the walls (for details we refer to
\cite[sec.~4.5]{pozrikidis92}).
On the other hand, the double layer density $\vec{q}$
is defined on the surfaces of both the particles and walls.
The formulation \eqref{eq:dist-flow} is complete, meaning that
any flow which satisfies \eqref{eq:stokes-1}--\eqref{eq:stokes-2} and
decays as $\vec{x} \to \infty$ can be represented in this way.

To derive the fundamental boundary integral equation, which is
used to determine the double layer density $\vec{q}$ in
\eqref{eq:dist-flow}, we insert \eqref{eq:dist-flow} into
\eqref{eq:total-flow} and then let $\vec{x} \in \Omega$ approach the
solid boundary $\Gamma$. Enforcing no-slip boundary conditions on
$\Gamma$ yields, recalling the jump condition \eqref{eq:double-jump},
\begin{equation}
  \label{eq:proto-bie}
  \vec{u}(\vec{x}) = \vec{u}_\text{bg}(\vec{x}) + \vec{\Dp}[\Gamma,\vec{q}](\vec{x})
  - 4\pi\vec{q}(\vec{x})
  + \sum_{\alpha=1}^M \vec{\Vp}^{(\alpha)}[\vec{F}^{(\alpha)},
  \vec{\tau}^{(\alpha)}](\vec{x})
  = \vec{U}_\Gamma(\vec{x}), \qquad \vec{x} \in \Gamma.
\end{equation}
The presence of the term $-4\pi \vec{q}(\vec{x})$, which is due
to the jump condition, makes the boundary integral equation
\eqref{eq:proto-bie} a Fredholm integral equation of the second
kind. The right-hand side $\vec{U}_\Gamma$ is the pointwise
velocity of the boundary $\Gamma$. We assume the walls to be
stationary and the particles to move as rigid bodies. This means
that, if we let $\Gw$ be the union of all wall
surfaces and $\Gp^{(\alpha)}$ the surface of particle
$\alpha$,
\begin{equation}
  \label{eq:rigid-body-motion}
  \vec{U}_\Gamma(\vec{x}) = \begin{cases}
    \vec{0}, & \vec{x} \in \Gw, \\
    \vec{U}_\text{RBM}^{(\alpha)} +
    \vec{\Omega}_\text{RBM}^{(\alpha)} \times
    (\vec{x} - \vec{y}_\text{c}^{(\alpha)}),
    & \vec{x} \in \Gp^{(\alpha)},
  \end{cases}
\end{equation}
where $\vec{U}_\text{RBM}^{(\alpha)}$ and
$\vec{\Omega}_\text{RBM}^{(\alpha)}$ are the translational and
angular velocity, respectively, of particle $\alpha$ (with RBM
denoting rigid body motion).

As mentioned in section~\ref{sec:introduction}, the viscous
resistance that the particles experience from the fluid is
related to their velocities. In the \emph{resistance problem},
the velocities (i.e.\ $\vec{U}_\text{RBM}^{(\alpha)}$ and
$\vec{\Omega}_\text{RBM}^{(\alpha)}$ for each particle) are
specified in \eqref{eq:proto-bie}--\eqref{eq:rigid-body-motion},
while in the \emph{mobility problem}, the viscous forces and
torques (i.e.\ $\vec{F}^{(\alpha)}$ and $\vec{\tau}^{(\alpha)}$ for
each particle) are specified \cite[p.~129]{pozrikidis92}. The
boundary integral equations resulting from these two problems are
described in more detail below. In both cases, the resulting
integral equation is discretized using the Nyström method, as
described in section~\ref{sec:disc-quadrature}.

\subsection{The resistance problem}

In this case, the velocities $\vec{U}_\text{RBM}^{(\alpha)}$ and
$\vec{\Omega}_\text{RBM}^{(\alpha)}$ of all particles are known,
while the corresponding forces~$\vec{F}^{(\alpha)}$ and
torques~$\vec{\tau}^{(\alpha)}$ are to be computed.
Following \cite[p.~130]{pozrikidis92}, the forces and torques
are related to the unknown double layer density $\vec{q}$ by
stipulating
\begin{equation}
  \label{eq:resistance-forces}
  \vec{F}^{(\alpha)}[\vec{q}] = \int_{\Gp^{(\alpha)}}
  \vec{q}(\vec{y}) \, \D S_{\vec{y}} \qquad \text{and} \qquad
  \vec{\tau}^{(\alpha)}[\vec{q}] = \int_{\Gp^{(\alpha)}}
  (\vec{y} - \vec{y}_\text{c}^{(\alpha)}) \times \vec{q}(\vec{y})
  \, \D S_{\vec{y}}.
\end{equation}
These relations are inserted into \eqref{eq:proto-bie}, which can
then be rearranged as
\begin{equation}
  \label{eq:resistance-bie}
  \vec{\Dp}[\Gamma,\vec{q}](\vec{x})
  - 4\pi\vec{q}(\vec{x})
  + \sum_{\alpha=1}^M \vec{\Vp}^{(\alpha)}[\vec{F}^{(\alpha)}[\vec{q}],
  \vec{\tau}^{(\alpha)}[\vec{q}]](\vec{x})
  = \vec{U}_\Gamma(\vec{x}) - \vec{u}_\text{bg}(\vec{x}),
  \qquad \vec{x} \in \Gamma.
\end{equation}
After solving this integral equation for $\vec{q}$, the forces
and torques can be computed using \eqref{eq:resistance-forces},
and the flow field can then be computed using
\eqref{eq:total-flow}--\eqref{eq:dist-flow}.

\subsection{The mobility problem}

In this case, the force~$\vec{F}^{(\alpha)}$ and
torque~$\vec{\tau}^{(\alpha)}$ exerted on the fluid by each
particle (which for a noninertial particle are equal to
the net external force and torque acting on the particle) are
known, but not the particle velocities $\vec{U}_\text{RBM}^{(\alpha)}$
and $\vec{\Omega}_\text{RBM}^{(\alpha)}$. Following
\cite[p.~135]{pozrikidis92}, the velocities are related to the double
layer density~$\vec{q}$ by
\begin{align}
  \label{eq:relation-URBM-q}
  \vec{U}_\text{RBM}^{(\alpha)}[\vec{q}]
  &= -\frac{4\pi}{\absi{\Gp^{(\alpha)}}}
  \int_{\Gp^{(\alpha)}}
  \vec{q}(\vec{y}) \, \D S_{\vec{y}}, \\[3pt]
  \label{eq:relation-OmegaRBM-q}
  \vec{\Omega}_\text{RBM}^{(\alpha)}[\vec{q}]
  &= -4\pi \sum_{n=1}^3 \frac{\vec{\omega}_n^{(\alpha)}}{A_n^{(\alpha)}}
  \gps{\bigg}{\vec{\omega}_n^{(\alpha)} \cdot
  \int_{\Gp^{(\alpha)}} (\vec{y} -
  \vec{y}_\text{c}^{(\alpha)}) \times \vec{q}(\vec{y}) \, \D S_{\vec{y}}}.
\end{align}
Here, $\absi{\Gp^{(\alpha)}}$ is the surface area of
$\Gp^{(\alpha)}$, and
\begin{equation}
  A_n^{(\alpha)} = \int_{\Gp^{(\alpha)}} \abs{
    \vec{\omega}_n^{(\alpha)} \times (\vec{y} - \vec{y}_\text{c}^{(\alpha)})
  }^2 \, \D S_{\vec{y}},
\end{equation}
while $\vec{\omega}_n^{(\alpha)}$ are three linearly independent
unit vectors which must satisfy
\begin{equation}
  \frac{1}{\sqrt{A_m^{(\alpha)} A_n^{(\alpha)}}}
  \int_{\Gp^{(\alpha)}} \fb{\vec{\omega}_m^{(\alpha)}
  \times (\vec{y} - \vec{y}_\text{c}^{(\alpha)})} \cdot
  \fb{\vec{\omega}_n^{(\alpha)}
  \times (\vec{y} - \vec{y}_\text{c}^{(\alpha)})} \, \D S_{\vec{y}} = \delta_{mn},
  \qquad m, n = 1,2,3.
\end{equation}
The boundary integral equation \eqref{eq:proto-bie} can then be
rearranged as
\begin{equation}
  \label{eq:mobility-bie}
  \vec{\Dp}[\Gamma,\vec{q}](\vec{x})
  - 4\pi\vec{q}(\vec{x})
  - \vec{U}_\Gamma[\vec{q}](\vec{x})
  = -\vec{u}_\text{bg}(\vec{x})
  - \sum_{\alpha=1}^M \vec{\Vp}^{(\alpha)}[\vec{F}^{(\alpha)},
  \vec{\tau}^{(\alpha)}](\vec{x}),
  \qquad \vec{x} \in \Gamma,
\end{equation}
where $\vec{U}_\Gamma[\vec{q}]$ is given by
\eqref{eq:rigid-body-motion} but with
$\vec{U}_\text{RBM}^{(\alpha)}$ and
$\vec{\Omega}_\text{RBM}^{(\alpha)}$ replaced by the expressions
in \eqref{eq:relation-URBM-q} and \eqref{eq:relation-OmegaRBM-q},
respectively. After solving \eqref{eq:mobility-bie} for
$\vec{q}$, the velocities can be computed using
\eqref{eq:relation-URBM-q}--\eqref{eq:relation-OmegaRBM-q}, and
the flow field can be computed using
\eqref{eq:total-flow}--\eqref{eq:dist-flow}.

\FloatBarrier
\chapter{Discretization and quadrature}
\label{sec:disc-quadrature}

In order to solve the boundary integral equation
\eqref{eq:resistance-bie} associated with the resistance problem,
or the boundary integral equation \eqref{eq:mobility-bie}
associated with the mobility problem, the integral operators in
these equations must be discretized. This amounts to discretizing
the double layer potential $\vec{\Dp}$, as well as the integrals
occurring in relation \eqref{eq:resistance-forces} for the resistance problem,
or relation \eqref{eq:relation-URBM-q}--\eqref{eq:relation-OmegaRBM-q}
for the mobility problem.
Following \cite{klinteberg16b}, we introduce the notation
\begin{equation}
  \label{eq:integral}
  \mathrm{I}[f] = \int_\Gamma f(\vec{y}) \, \D S(\vec{y})
\end{equation}
for the integral of the arbitrary function $f$ over the surface
$\Gamma$. We introduce a quadrature rule $Q_N$ called the
\emph{direct quadrature rule}, defined by a set
of $N$ nodes $\vec{x}_i \in \Gamma$ and weights $w_i \in \Real$,
$i=1,\ldots,N$. The
details of this quadrature rule is specified in
sections~\ref{sec:direct-quad-particles} and
\ref{sec:direct-quad-walls}. Using the direct quadrature rule $Q_N$, the integral in
\eqref{eq:integral} can be approximated as
\begin{equation}
  \label{eq:quadrature-rule}
  \mathrm{I}[f] \approx Q_N[f] = \sum_{i=1}^N f(\vec{x}_i) w_i.
\end{equation}
We denote an integral quantity approximated by $Q_N$ with a
superscript $h$, for example the double layer potential
\begin{equation}
  \label{eq:direct-quad-double-layer}
  \Dp^h_i[\Gamma, \vec{q}](\vec{x}) = Q_N[T_{ijk}(\vec{x} -
  {\cdot}) q_j(\cdot) n_k(\cdot)].
\end{equation}
We then discretize \eqref{eq:resistance-bie} or
\eqref{eq:mobility-bie} using the \emph{Nyström
method}~\cite[ch.~4]{atkinson97}\cite[sec.~12.2]{kress14}, in which the integral equation is enforced in the quadrature
nodes. For the resistance problem, \eqref{eq:resistance-bie} then becomes
\begin{equation}
  \label{eq:resistance-disc}
  \vec{\Dp}^h[\Gamma,\vec{q}](\vec{x}_i)
  - 4\pi\vec{q}(\vec{x}_i)
  + \sum_{\alpha=1}^M \vec{\Vp}^{(\alpha),h}[\vec{F}^{(\alpha)}[\vec{q}],
  \vec{\tau}^{(\alpha)}[\vec{q}]](\vec{x}_i)
  = \vec{U}_\Gamma(\vec{x}_i) - \vec{u}_\text{bg}(\vec{x}_i),
  \qquad i=1,\ldots,N.
\end{equation}
For the mobility problem, \eqref{eq:mobility-bie} becomes
\begin{equation}
  \label{eq:mobility-disc}
  \vec{\Dp}^h[\Gamma,\vec{q}](\vec{x}_i)
  - 4\pi\vec{q}(\vec{x}_i)
  - \vec{U}^h_\Gamma[\vec{q}](\vec{x}_i)
  = -\vec{u}_\text{bg}(\vec{x}_i)
  - \sum_{\alpha=1}^M \vec{\Vp}^{(\alpha)}[\vec{F}^{(\alpha)},
  \vec{\tau}^{(\alpha)}](\vec{x}_i),
  \qquad i=1,\ldots,N.
\end{equation}
The superscript $h$ on $\vec{\Vp}^{(\alpha),h}$ in
\eqref{eq:resistance-disc} and $\vec{U}_\Gamma^h$ in
\eqref{eq:mobility-disc} signifies that these quantities, while
not integrals themselves, contain integrals -- namely
\eqref{eq:resistance-forces} or
\eqref{eq:relation-URBM-q}--\eqref{eq:relation-OmegaRBM-q} --
which are approximated using the direct quadrature rule $Q_N$. In both
cases, the resulting linear system is solved iteratively using
the generalized minimal residual method (GMRES).

\begin{figure}[h!]
  \centering
  \begin{tikzpicture}
    \node at (-0.1,-1.55) {\includegraphics[scale=0.7,trim=1cm 1.5cm 1cm 1.5cm,clip]{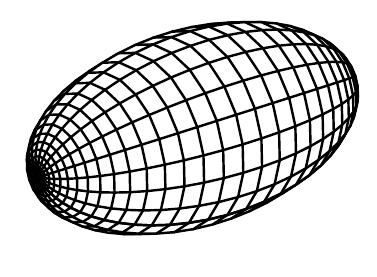}};
      \node at (-0.2,-4.8) {\includegraphics[scale=0.7,trim=5mm 1cm 5mm 1cm,clip]{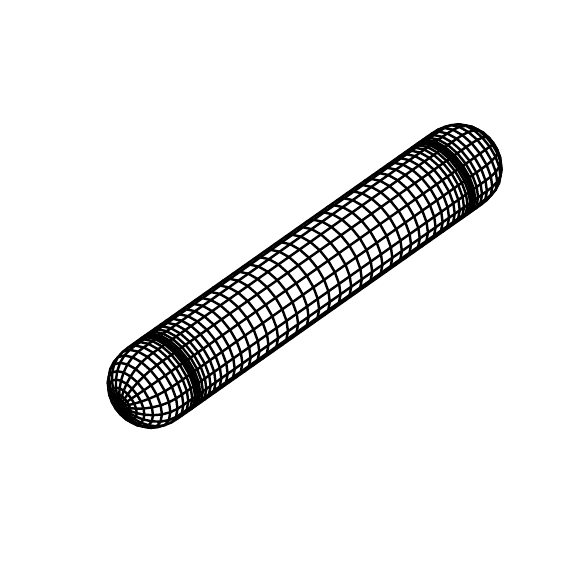}};
      \draw [line width=1.3pt, rounded corners=5pt, myRed] (-2.5,-0.3) rectangle (2.5,-6.5);
      \node [font=\bfseries] at (0,-2.8) {Spheroid};
      \node [font=\bfseries] at (0,-6.1) {Rod};
      \node [font=\bfseries,myRed] at (0,-6.9) {(a) Particles};
      \node [myRed] at (0,-7.4) {Particle-global quadrature};
    \begin{scope}[xshift=6cm]
      \node at (0,-1.5) {\includegraphics[scale=0.9]{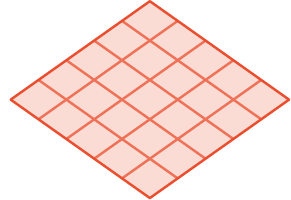}};
      \node at (0,-4.7) {\includegraphics[scale=0.9]{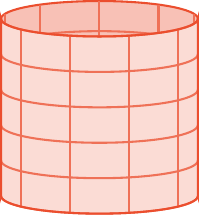}};
      \draw [line width=1.3pt, rounded corners=5pt, myRed]
      (-2.5,-0.3) rectangle (2.5,-6.5);
      \node [font=\bfseries] at (0,-2.8) {Plane wall};
      \node [font=\bfseries] at (0,-6.1) {Pipe};
      \node [font=\bfseries,myRed] at (0,-6.9) {(b) Walls};
      \node [myRed] at (0,-7.4) {Local patch-based quadrature};
    \end{scope}
  \end{tikzpicture}
  \caption{The geometrical objects considered in this paper are
  (a) particles (spheroids and rods) and (b) walls (plane walls and
  pipes). The quadrature rule is different for particles and
  walls.}
  \label{fig:geometry}
\end{figure}
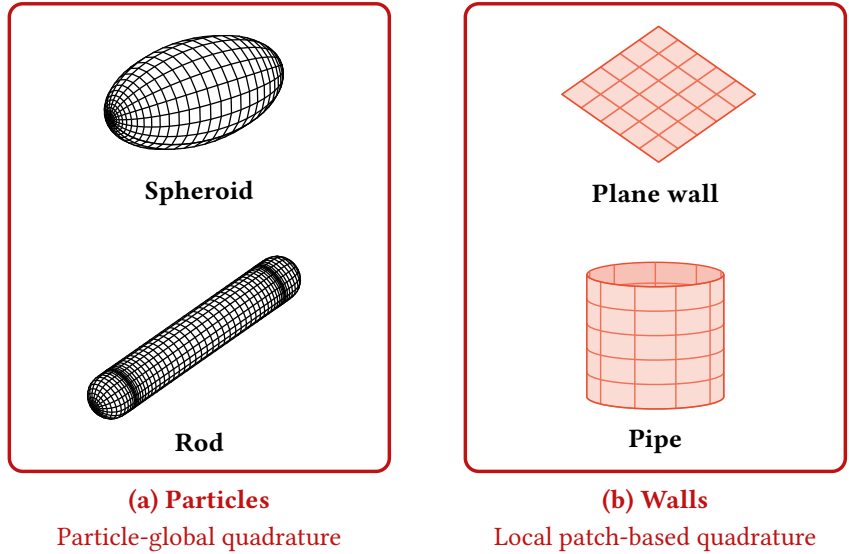

In this paper, we consider two distinct types of geometrical
objects, namely particles and walls, as indicated by
Figure~\ref{fig:geometry}. Particles are mobile rigid bodies
immersed in the fluid, while walls are stationary and surround
the fluid domain. We consider two types of particles: spheroids,
which are given by a surface
\begin{equation}
  \label{eq:spheroid-surface}
  \frac{x_1^2 + x_2^2}{a^2} + \frac{x_3^2}{c^2} = 1
\end{equation}
in local coordinates; and rods, which consist of a cylinder
with rounded caps, described in appendix~\ref{app:smooth-rod}.
We also consider two types of walls, namely
plane walls and pipes with circular cross section. Both wall
geometries extend to infinity in the periodic setting, but we
discretize only the part of each object that lies inside the
primary cell.

The nature of the direct quadrature rule $Q_N$ is different for
particles and walls: for particles, it is a particle-global
quadrature rule described in
section~\ref{sec:direct-quad-particles}, while for walls it is a
local patch-based quadrature rule described in
section~\ref{sec:direct-quad-walls}. The special quadrature
method for particles and walls is introduced in
section~\ref{sec:upsamp-special-quad}.

It should be noted that all geometrical objects shown in
Figure~\ref{fig:geometry} are smooth, i.e.\ of class $C^\infty$.
The construction of the smooth rod particles is described in
appendix~\ref{app:smooth-rod}. In section~\ref{sec:smoothness},
we consider the effect that a nonsmooth object would have on the
convergence of the special quadrature method.

\section{Direct quadrature for particles}
\label{sec:direct-quad-particles}

The discretization and direct quadrature rule of the spheroids are
exactly the same as in \cite{klinteberg16b}, while for the rods
they are a slight variation of the former. Both kinds of particles are
axisymmetric, and their parametrizations take this into account,
with one parameter $\varphi \in [0, 2\pi)$ varying in the
azimuthal direction and the other parameter $\theta \in [0,\pi]$
varying in the polar direction. For instance, the spheroid
\eqref{eq:spheroid-surface} is parametrized using spherical
coordinates as
\begin{equation}
  \label{eq:spheroid-parametrization}
  \begin{cases}
    x_1 = a \sin\theta \cos\varphi, \\
    x_2 = a \sin\theta \sin\varphi, \\
    x_3 = c \cos\theta.
  \end{cases}
\end{equation}
It is discretized using a tensorial grid with $n_\theta \times
n_\varphi$ grid points. For the polar direction, let $(\theta_i,
\lambda^\theta_i)$, $i=1,\ldots,n_\theta$, be the nodes and
weights of an $n_\theta$-point Gauss--Legendre quadrature rule
\cite[\href{https://dlmf.nist.gov/3.5.v}{sec.~3.5(v)}]{dlmf} on
the interval $[0,\pi]$. For the azimuthal direction, let
$(\varphi_j, \lambda^\varphi_j)$, $j=1,\ldots,n_\varphi$, be the
nodes and weights of the trapezoidal rule on the interval
$[0,2\pi)$. Since the integrand is periodic on this interval, the
trapezoidal rule has spectral accuracy in this case \cite{trefethen14}.
The resulting tensorial quadrature rule, called the direct
quadrature rule of the spheroid, is
\begin{equation}
  \label{eq:spheroidal-rule}
  Q_{n_\theta n_\varphi}[f] = \sum_{i=1}^{n_\theta}
  \sum_{j=1}^{n_\varphi} f(\vec{x}(\theta_i,\varphi_j))
  W_\text{sph}(\theta_i,\varphi_j) \lambda_i^\theta \lambda_j^\varphi,
\end{equation}
where $W_\text{sph}(\theta,\varphi)$ is the area element associated with the
parametrization \eqref{eq:spheroid-parametrization}.

\begin{figure}[h!]
  \centering
  \includegraphics[trim=0mm 2mm 0mm 0mm,clip]{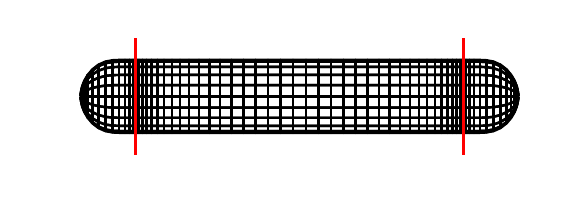}
  \caption{The grid on the rods consists of three parts: two caps
  and a middle cylinder.}
  \label{fig:rod-parts}
\end{figure}

The rod consists of a cylinder with rounded caps. While the
surface is smooth everywhere, the grid is divided into three
parts as shown in Figure~\ref{fig:rod-parts}. The reason for this
is to be able to increase the resolution at the caps
independently of the resolution at the middle of the rod.%
\footnote{We also tested a discretization of the rod using a grid
spanning the whole rod without dividing it into parts. We did not
find any significant improvement in the quadrature error or
computational cost from using such a grid rather than the one
shown in Figure~\ref{fig:rod-parts}.}

The rod is parametrized as
\begin{equation}
  \label{eq:rod-parametrization}
  \begin{cases}
    x_1 = \varrho(\theta; L, R) \cos\varphi, \\
    x_2 = \varrho(\theta; L, R) \sin\varphi, \\
    x_3 = \beta(\theta; L, R),
  \end{cases}
\end{equation}
where $L$ is the length of the rod and $R$ its radius. The shape
functions $\varrho(\cdot\,;L,R) : [0,\pi] \to [0,R]$ and
$\beta(\cdot\,;L,R) : [0,\pi] \to [-\tfrac{1}{2}L, \tfrac{1}{2}L]$ are described in
appendix~\ref{app:smooth-rod}. They are chosen such that $\theta \in [0,
\pi/3] = I_1$ and $\theta \in [2\pi/3, \pi] = I_3$ correspond to the two caps,
while $\theta \in [\pi/3, 2\pi/3] = I_2$ corresponds to the middle
part. Each cap is discretized using $n_1 \times n_\varphi$ grid
points, and the middle cylinder is discretized using $n_2 \times
n_\varphi$ grid points, so the total grid has $(2n_1+n_2) \times
n_\varphi$ grid points. The trapezoidal rule is again used in the
azimuthal direction. In the polar direction, a separate
Gauss--Legendre quadrature rule is used for each of the three
parts. The tensorial direct quadrature rule of the rod is thus (with
$n_3=n_1$)
\begin{equation}
  \label{eq:rod-rule}
  Q_{(2n_1+n_2)n_\varphi}[f] = \sum_{k=1}^3 \sum_{i=1}^{n_k}
  \sum_{j=1}^{n_\varphi} f(\vec{x}(\theta_i^k,\varphi_j))
  W_\text{rod}(\theta_i^k,\varphi_j) \lambda_i^k \lambda_j^\varphi,
\end{equation}
where $(\theta^k_i, \lambda^k_i)$, $i=1,\ldots,n_k$, are the
nodes and weights of an $n_k$-point Gauss--Legendre quadrature on
the interval $I_k$, and $W_\text{rod}(\theta, \varphi)$ is the
area element associated with \eqref{eq:rod-parametrization}.

The direct quadrature rules \eqref{eq:spheroidal-rule} and
\eqref{eq:rod-rule} are both particle-global in the sense
that each particle is treated as a single unit, and the
quadrature rule is applied to the particle as a whole. The
quadrature rules has spectral accuracy for smooth integrands,
i.e.\ it converges exponentially as the number of grid points
increases.

\section{Direct quadrature for walls}
\label{sec:direct-quad-walls}

\begin{figure}[ht!]
  \centering\small
  \begin{minipage}[b]{0.5\textwidth}%
    \centering
    \includegraphics[scale=0.5,trim=0mm 8mm 0mm 8mm,clip]{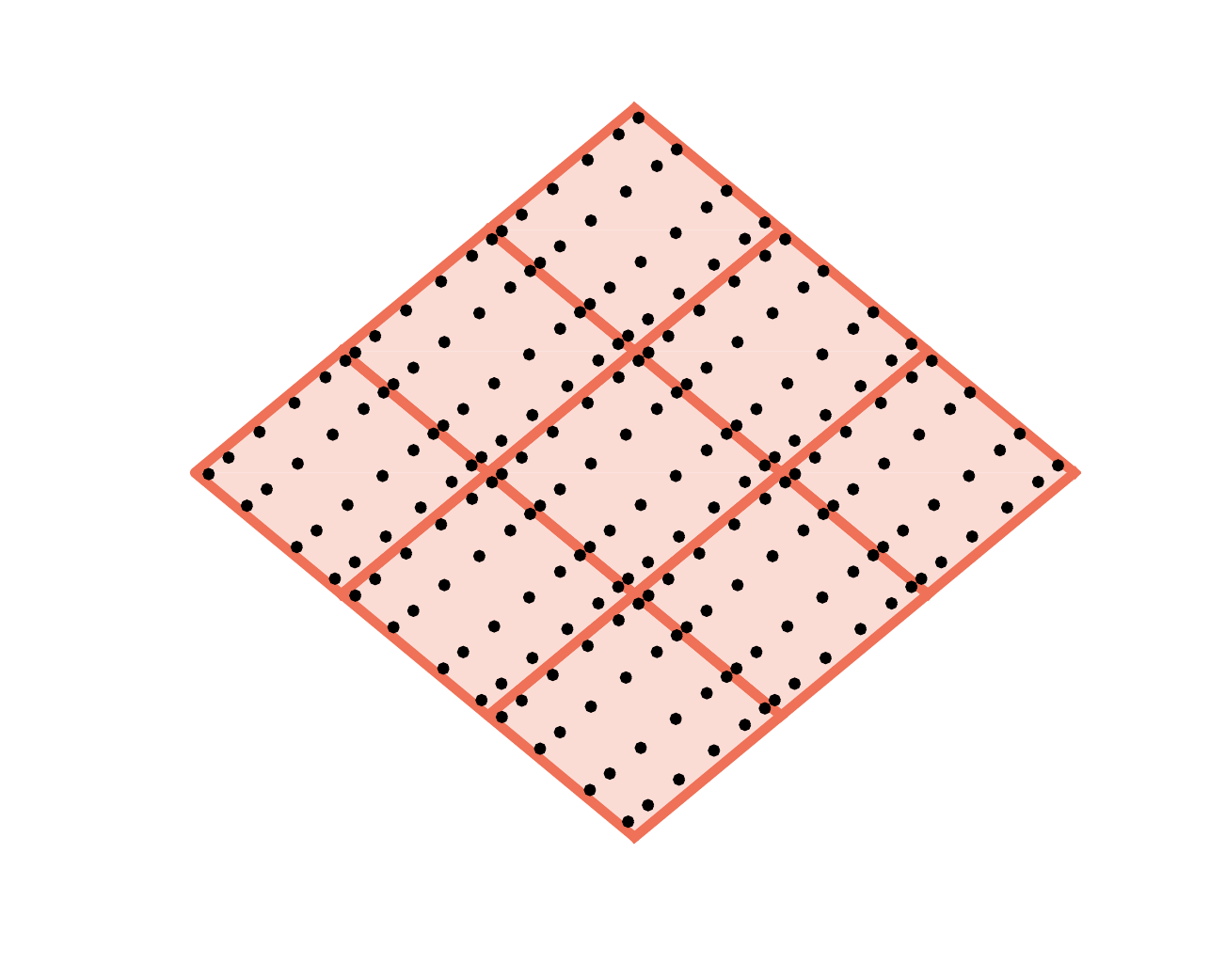}\\
    (a) Plane wall
  \end{minipage}%
  \begin{minipage}[b]{0.5\textwidth}%
    \centering
    \includegraphics[scale=0.8,trim=0mm 5mm 0mm 4mm,clip]{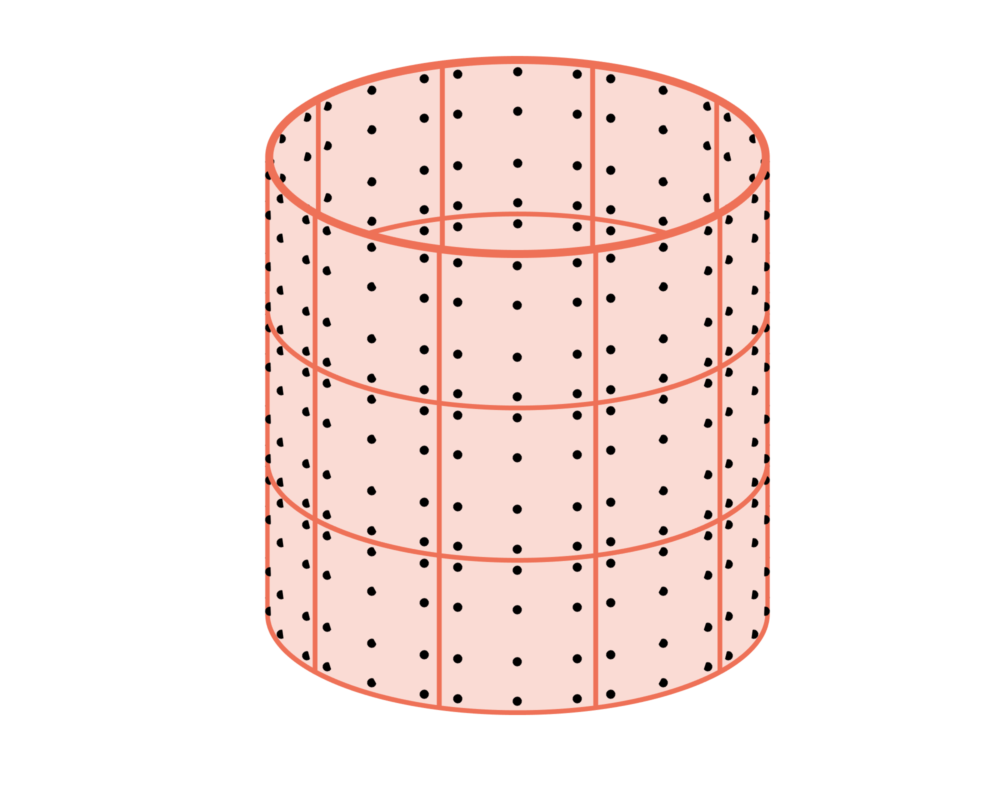}\\
    (b) Pipe
  \end{minipage}%
  \caption{(a)~A plane wall divided into $3 \times 3$ patches,
  each patch discretized using $6 \times 4$ grid points. (b)~A
  pipe divided into $3 \times 10$ patches, each patch discretized
  using $4 \times 3$ grid points.}
  \label{fig:wall-patches}
\end{figure}

The wall geometries are present only in the periodic setting, and
then only the part inside the primary cell needs to be
discretized. For the plane wall, this part consists of a flat
rectangle of size $L_1 \times L_2$, which is divided into $P_1
\times P_2$ flat subrectangles, called patches. Each patch is
discretized using a tensorial grid with $n_1 \times n_2$
Gauss--Legendre grid points, as shown in
Figure~\ref{fig:wall-patches}~(a).
In each direction of the patch, an $n_d$-point Gauss--Legendre quadrature
rule is used with nodes and weights $(s_i^d, \lambda_i^d)$, $i=1,\ldots,n_d$,
$d=1,2$. The resulting tensorial direct quadrature rule of the patch is
\begin{equation}
  \label{eq:wall-rule}
  Q_{n_1 n_2}[f] = \sum_{i=1}^{n_1}
  \sum_{j=1}^{n_2} f(\vec{x}(s_i^1,s_j^2))
  W_\text{wall}(s_i^1,s_j^2) \lambda_i^1 \lambda_j^2,
\end{equation}
where $W_\text{wall}(s^1, s^2)$ is the area element of the wall.

The part of the pipe in the primary cell consists of a cylinder
with radius $R_\text{c}$ and length $L_\text{c}$. Like the plane
wall, it is divided into rectangular patches, but these are
curved, as seen in Figure~\ref{fig:wall-patches}~(b). Apart from that,
the discretization and quadrature rule are the same as for the
plane wall. The direct quadrature rule of the pipe is thus also
given by \eqref{eq:wall-rule}, the only difference compared to
the plane wall being the area element $W_\text{wall}$ and the
parametrization $(s^1,s^2) \mapsto \vec{x}$.

The direct quadrature rule \eqref{eq:wall-rule} is local in the
sense that the wall is subdivided into smaller patches, and the
quadrature rule is applied to each patch separately. The grid can
be refined in two different ways: by adding more grid points to
each patch (which we call $n$-refinement), or by reducing the
size of the patches and thus increasing their number (which we
call $P$-refinement). Under $n$-refinement, the quadrature rule
has spectral accuracy like the direct quadrature rule of the
particles, while under $P$-refinement the quadrature rule is
polynomially accurate with order determined by $n_1$ and $n_2$.

\section{Special quadrature: upsampled quadrature and quadrature
by expansion (QBX)}
\label{sec:upsamp-special-quad}

The double layer potential $\vec{\Dp}$ given by
\eqref{eq:double-layer} is challenging to compute using direct
quadrature in two different situations, in both cases due to its
kernel $\mat{T}$. Firstly, when the evaluation point $\vec{x}$ is
on $\Gamma$ itself, $\mat{T}$ becomes singular at the point
$\vec{y} = \vec{x}$ (we refer to this as the \emph{singular case}
or the \emph{onsurface evaluation case}). The integral exists as
an improper integral as long as $\Gamma$ is a Lyapunov surface
\cite[p.~37]{pozrikidis92}, but clearly a special quadrature
method of some sort is needed to compute it.
Secondly, when $\vec{x} \in \Omega$ is close to $\Gamma$, but not
on $\Gamma$, $\mat{T}$ becomes very peaked and thus hard to resolve
using the direct quadrature rule (we refer to this as the
\emph{nearly singular case} or the \emph{offsurface evaluation
case}).

\begin{figure}[b!]
  \centering\small
  \begin{minipage}{0.45\textwidth}%
    \centering
    \includegraphics[scale=0.5,trim=0mm 22mm 0mm 22mm,clip]%
    {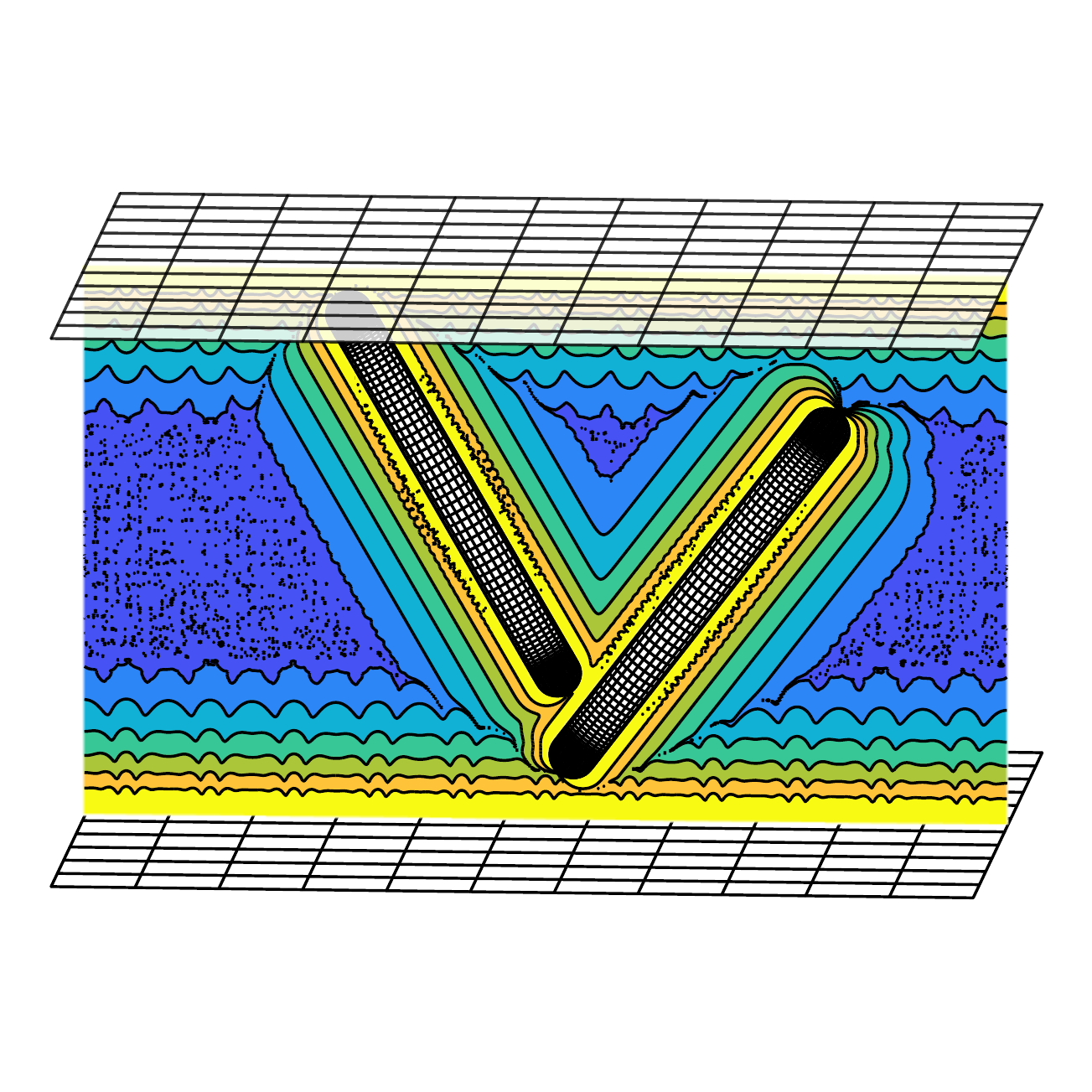}\\
    (a) Direct quadrature error
  \end{minipage}%
  \begin{minipage}{0.45\textwidth}%
    \centering
    \includegraphics[scale=0.5,trim=0mm 22mm 0mm 22mm,clip]%
    {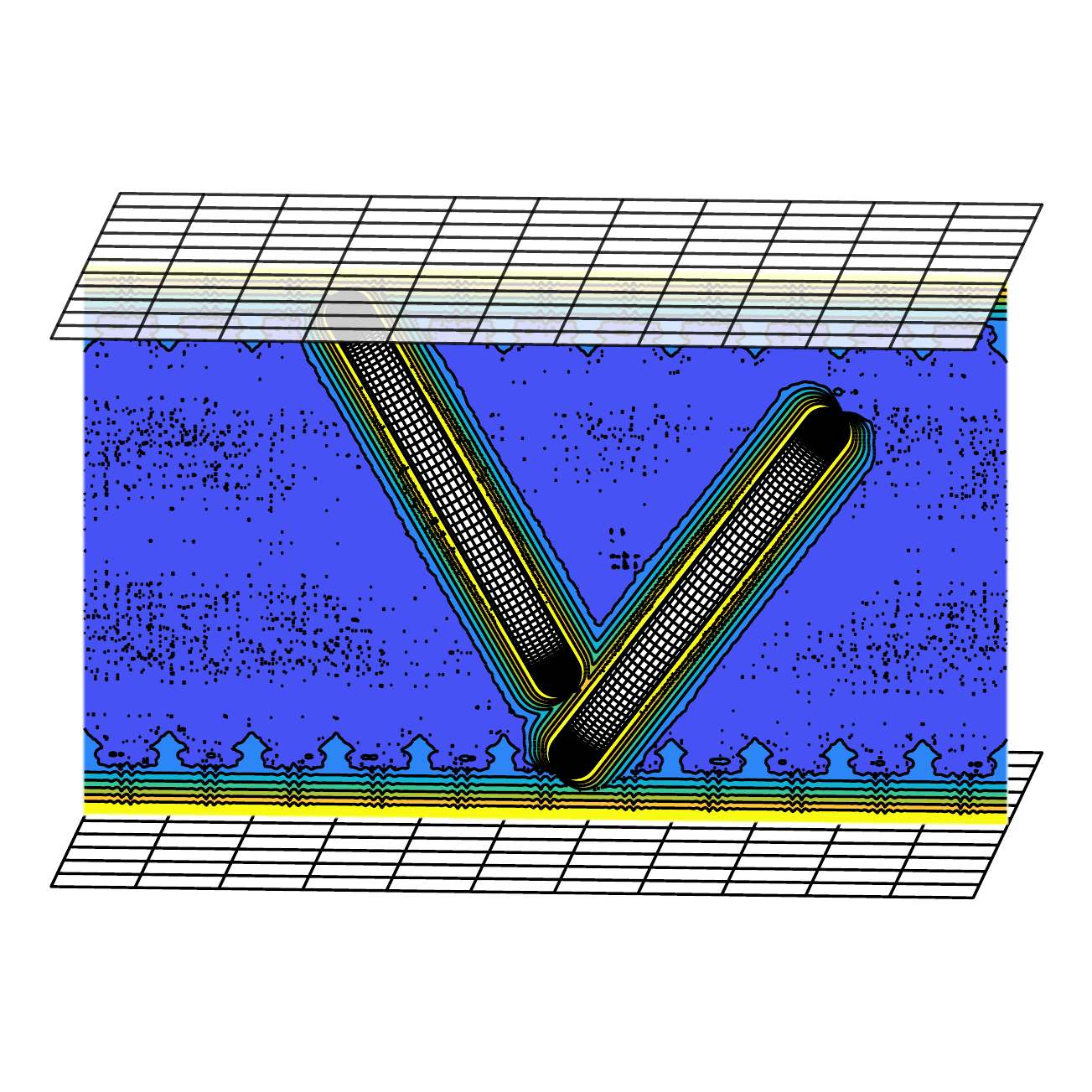}\\
    (b) Upsampled quadrature error, $\kappa = 2$
  \end{minipage}%
  \begin{minipage}{0.1\textwidth}%
    \centering
    \includegraphics[scale=1]{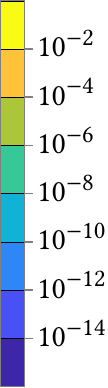}\\
    \vphantom{(c)}
  \end{minipage}%
  \caption{Relative error in the centre plane when evaluating the
  stresslet identity \eqref{eq:stresslet-identity} for two rod
  particles between a pair of parallel horizontal plane walls in a
  periodic setting, using the direct quadrature rule in~(a) and
  the upsampled quadrature rule with upsampling factor $\kappa=2$
  in~(b). Note that the error is still large very close to the
  particles and walls in (b). The density $\vec{q}$ is constant
  in this example, so the error in (a) and (b) comes entirely
  from the nearly singular behaviour of the stresslet~$\mat{T}$.}
  \label{fig:direct-quad}
\end{figure}%

The singular case is always present when solving the boundary
integral equation \eqref{eq:proto-bie}, while the
nearly singular case occurs when particles are close to each
other or close to a wall, and also if the flow field
\eqref{eq:total-flow}--\eqref{eq:dist-flow} is to be
computed close to a particle or wall. The latter situation is
illustrated in Figure~\ref{fig:direct-quad}~(a),
where the error grows exponentially as the evaluation
point~$\vec{x}$ approaches the boundary~$\Gamma$. This
behaviour is well-known, and in two dimensions there are error
estimates available for the Laplace and Helmholtz potentials in
\cite{klinteberg17} and for the Stokes potential in
\cite{palsson19}. To compute the double layer potential
accurately close to a particle or wall, \emph{special quadrature}
is needed. Here, we consider two types of special quadrature:
upsampled quadrature and quadrature by expansion (QBX).

\begin{figure}[h!]
  \centering
  \includegraphics{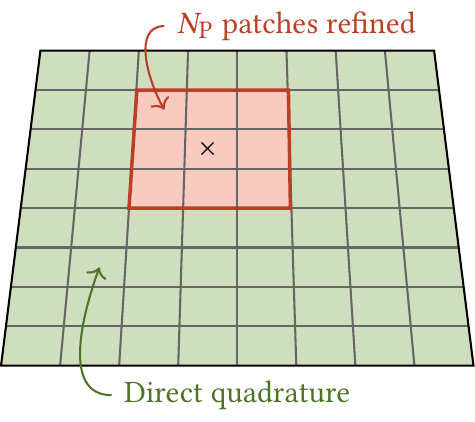}
  \caption{Upsampled quadrature for a plane wall: the
  $N_\text{P}$ patches closest to the evaluation point (marked
  with $\times$) are refined, here for $N_\text{P}=9$. The other
  patches are treated using direct quadrature without refinement.}
  \label{fig:wall-shells}
\end{figure}

Assuming that the density $\vec{q}$ itself is well-resolved on
the grid, \emph{upsampled quadrature} provides a partial solution
for the nearly singular case. In upsampled quadrature,
the double layer density $\vec{q}$ is interpolated onto a grid
refined by a factor $\kappa$ in both directions, and the integral
is then evaluated using direct quadrature on the finer grid. For
the particle-global quadrature rules in section~\ref{sec:direct-quad-particles},
the grid of the whole particle is refined (increasing the number
of grid points of each individual quadrature rule). The density is
interpolated onto the finer grid using trigonometric
interpolation in the azimuthal direction and barycentric Lagrange
interpolation \cite{berrut04} in the polar direction. For the
local patch-based quadrature rules in
section~\ref{sec:direct-quad-walls}, only the $N_\text{P}$
patches closest to the evaluation point $\vec{x}$ are refined,
using $n$-refinement (thus increasing the number of grid points
on them); other patches are sufficiently far away from the singularity
that direct quadrature can be used. This is illustrated in
Figure~\ref{fig:wall-shells} for $N_\text{P}=9$. The refinement
has spectral accuracy for both particles and walls. Since all
geometrical objects are rigid, interpolation matrices can be
precomputed.

\begin{figure}[h!]
  \centering\small
  \begin{minipage}[b]{0.35\textwidth}%
    \centering
    \includegraphics{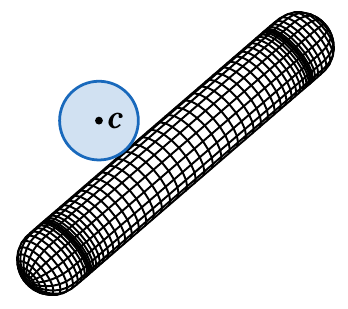}\\
    (a)
  \end{minipage}%
  \begin{minipage}[b]{0.35\textwidth}%
    \centering
    \includegraphics{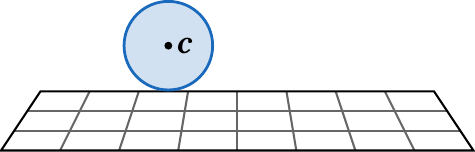}\\
    (b)
  \end{minipage}%
  \caption{The idea behind QBX is to make a series expansion of
  the potential close to a particle (a) or wall (b), valid inside
  a ball of convergence shown as a blue disc. The expansion is
  also valid at the point where the ball of convergence touches
  the boundary.}
  \label{fig:qbx-idea}
\end{figure}

As Figure~\ref{fig:direct-quad}~(b) shows,
upsampled quadrature pushes the region where the error is large
closer to the boundary $\Gamma$. However, the error will always
be large very close to $\Gamma$ no matter how large the
upsampling factor $\kappa$ is. To be able to achieve small errors
arbitrarily close to $\Gamma$, we use a special quadrature rule
specifically designed for layer potentials with singular kernels, namely
\emph{quadrature by expansion} (QBX) \cite{klockner13,barnett14}.
The idea behind QBX is to make a local series expansion of the
potential $\vec{\Dp}$ in the fluid domain, which converges
rapidly since $\vec{\Dp}$ is smooth all the way up to the
boundary $\Gamma$. The expansion is made
around a point $\vec{c}$, called the \emph{expansion centre},
which is inside the fluid domain (i.e.\ not on $\Gamma$), and it can
be used to evaluate the potential inside a ball around $\vec{c}$
called the \emph{ball of convergence}, as shown in
Figure~\ref{fig:qbx-idea}.
The expansion is valid even at the
point where the ball touches $\Gamma$ \cite{epstein13}, and can therefore
be used in the singular case as well as the nearly singular case.
The application of QBX to the Stokes double layer potential
$\vec{\Dp}$ will be described in detail in section~\ref{sec:qbx-detail}.

\begin{figure}[h!]
  \centering
  \begin{minipage}[b]{0.4\textwidth}
    \centering
    \includegraphics{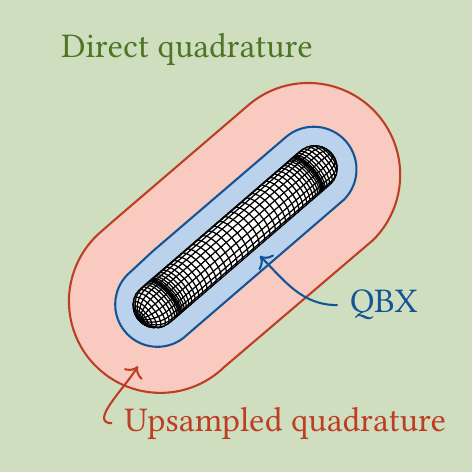}\\
    (a)
  \end{minipage}%
  \begin{minipage}[b]{0.4\textwidth}
    \centering
    \includegraphics{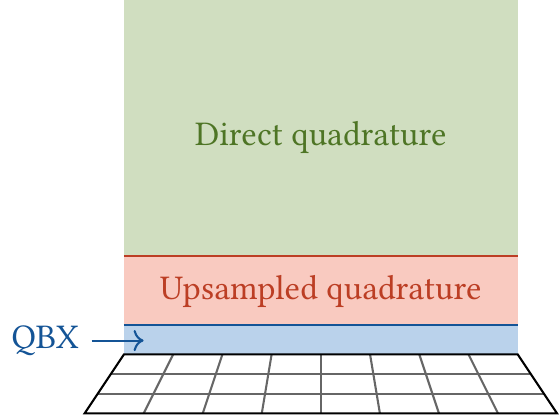}\\
    \hspace*{9mm}(b)
  \end{minipage}%
  \caption{The regions of the combined quadrature strategy, shown
  here for a rod particle in (a) and a plane wall in (b).
  Depending on the location of the evaluation point $\vec{x}$, it is
  treated using direct quadrature, upsampled quadrature or QBX.}
  \label{fig:combined-strategy}
\end{figure}

In this paper, we use a \emph{combined quadrature} strategy, where
direct quadrature is used far away from the boundary,
upsampled quadrature is used in an intermediate region, and QBX
is used in a small region closest to the boundary, as illustrated
in Figure~\ref{fig:combined-strategy}. For each particle and
wall in the geometry, the evaluation point~$\vec{x}$ is
classified into one of these three regions, and the contribution
to the double layer potential $\vec{\Dp}$ from that particle or
wall is computed as follows:
\begin{itemize}
  \item
    If $\vec{x}$ is in the direct quadrature region, the double
    layer potential \eqref{eq:double-layer} is computed using
    direct quadrature \eqref{eq:direct-quad-double-layer} over
    the whole particle or wall, as described in
    sections~\ref{sec:direct-quad-particles} and
    \ref{sec:direct-quad-walls}.
  \item
    If $\vec{x}$ is in the upsampled quadrature region, the
    behaviour is different for particles and walls, as described
    above. For a particle, the density is upsampled globally on
    the whole particle surface and then integrated using direct
    quadrature on the fine grid. For a wall, the density is
    upsampled only on the $N_\text{P}$ patches closest to the
    evaluation point, while direct quadrature without upsampling
    is used on the other patches (as in Figure~\ref{fig:wall-shells}).
  \item
    If $\vec{x}$ is in the QBX region, the behaviour is similar
    to the upsampled quadrature region. For a particle, the
    density on the whole particle surface is used when computing
    the coefficients of the local expansion which is then used at
    the evaluation point (particle-global QBX). For a wall, only the density on the
    $N_\text{P}$ patches closest to the expansion centre $\vec{c}$
    is used to compute the expansion (local QBX), while the
    contribution from other patches is computed using direct
    quadrature. In other words, the expansion is computed using a
    truncated wall, with $N_\text{P}$ determining the number of
    patches in the truncated wall. The difference between
    particle-global and local QBX is described in more detail in
    section~\ref{sec:global-local-qbx}.
\end{itemize}
The total double layer potential at $\vec{x}$ is then retrieved
using superposition, i.e.\ by summing the contributions from all
particles and walls.

When using local QBX, the convergence rate of the local expansion
will depend on the ratio between the distance from $\vec{c}$ to
the wall and the distance from $\vec{c}$ to the edge of the
truncated wall \cite{siegel18}. We have observed that
$N_\text{P}=1$ is too low for the wall QBX region in our case,
since the expansion centre may then be too close to the edge of the
truncated wall. Setting $N_\text{P}=9$ seems to be sufficient to
remedy this, and increasing $N_\text{P}$ further has no effect.
We therefore fix $N_\text{P}=9$ both for the QBX region and the
upsampled quadrature region of plane walls and pipes, for the
rest of this paper.

The distances from the surface at which to switch from one
quadrature region to the next (i.e.\ direct quadrature, upsampled
quadrature, QBX) are parameters to be set, and these will be
discussed in section~\ref{sec:parameters}.

\FloatBarrier
\chapter{Quadrature by expansion for the Stokes double layer
potential}
\label{sec:qbx-detail}

In order to apply QBX to the Stokes double layer potential
$\vec{\Dp}$ given by \eqref{eq:double-layer}, we need to be able
to write down a local expansion of the potential. We use the same
approach as in \cite{klinteberg16b}, which is summarized in
section~\ref{sec:expansion-double-layer}. The differences between
particles (for which particle-global QBX is used) and walls (for
which local QBX is used) are summarized in section~\ref{sec:global-local-qbx}.
Finally, the precomputation scheme which is crucial for accelerating the
method is described in section~\ref{sec:qbx-precomp}.

\section{Local expansion of the double layer potential}
\label{sec:expansion-double-layer}

\begin{figure}[h!]
  \centering
  \includegraphics[scale=1]{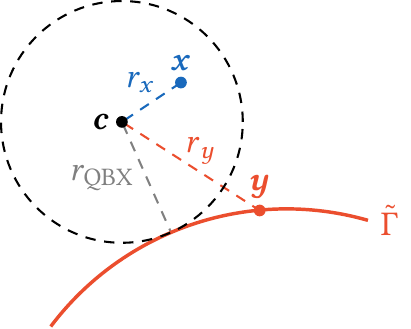}
  \caption{Illustration of the points $\vec{x}$, $\vec{y}$ and
  $\vec{c}$, and the ball of convergence of
  \eqref{eq:laplace-expansion}.}
  \label{fig:expansion-centre}
\end{figure}

Instead of expanding the double layer potential $\vec{\Dp}$
itself directly, we use the fact that $\vec{\Dp}$ can be
expressed in terms of the so-called dipole potential $\Lp$
using the relation \cite{tornberg08,klinteberg16b}
\begin{equation}
  \label{eq:dipole-relation}
  \Dp_i[\tilde{\Gamma}, \vec{q}](\vec{x}) = \gp{x_j \pd{}{x_i} -
  \delta_{ij}} \Lp[\tilde{\Gamma}, q_j \vec{n} + n_j \vec{q}](\vec{x})
  - \pd{}{x_i} \Lp[\tilde{\Gamma}, y_k q_k \vec{n} + y_k n_k \vec{q}](\vec{x}),
\end{equation}
where $\tilde{\Gamma}$ is any subset of $\Gamma$.
The dipole potential is the double layer potential of the Laplace
equation and is defined as
\begin{equation}
  \label{eq:dipole-potential}
  \Lp[\tilde{\Gamma}, \vec{\rho}](\vec{x}) = \int_{\tilde{\Gamma}}
  \vec{\rho}(\vec{y}) \cdot \grad_{\vec{y}}
  \frac{1}{\absi{\vec{x} - \vec{y}}} \, \D S(\vec{y}).
\end{equation}
The kernel of the dipole potential has a natural expansion based
on the so-called Laplace expansion
\begin{equation}
  \label{eq:laplace-expansion}
  \frac{1}{\absi{\vec{x}-\vec{y}}} = \sum_{l=0}^\infty
  \frac{4\pi}{2l+1} \sum_{m=-l}^l r_x^l Y_l^{-m}(\theta_x,
  \varphi_x) \frac{1}{r_y^{l+1}} Y_l^m(\theta_y, \varphi_y),
\end{equation}
where $Y_l^m$ is the spherical harmonics function of degree $l$
and order $m$ (defined as in \cite[Eq.~(3.5)]{epstein13}),
while $(r_x, \theta_x, \varphi_x)$ and $(r_y,
\theta_y, \varphi_y)$ are spherical coordinates of the points
$\vec{x}$ and $\vec{y}$ respectively, with respect to a chosen
expansion centre $\vec{c}$, as shown in
Figure~\ref{fig:expansion-centre}. The expansion
\eqref{eq:laplace-expansion} is valid as long as $r_x < r_y$,
i.e., it can be used for all $\vec{x}$ within the ball of radius
$r_\text{QBX} = \min_{\vec{y} \in \tilde{\Gamma}} r_y$ centred at
$\vec{c}$.

Inserting \eqref{eq:laplace-expansion} into
\eqref{eq:dipole-potential} leads to the expansion
\begin{equation}
  \label{eq:dipole-expansion}
  \Lp[\tilde{\Gamma}, \vec{\rho}](\vec{x}) = \sum_{l=0}^\infty
  \sum_{m=-l}^l r_x^l Y_l^{-m} (\theta_x, \varphi_x)
  z_{lm}[\vec{\rho}]
\end{equation}
of the dipole potential, where the coefficients
$z_{lm}[\vec{\rho}]$ are given by
\begin{equation}
  \label{eq:dipole-qbx-coeffs}
  z_{lm}[\vec{\rho}] = \frac{4\pi}{2l+1} \int_{\tilde{\Gamma}} \vec{\rho}(\vec{y}) \cdot
  \grad_{\vec{y}} \frac{1}{r_y^{l+1}} Y_l^m (\theta_y, \varphi_y)
  \, \D S(\vec{y}).
\end{equation}
These coefficients are complex-valued due to $Y_l^m$, but the
dipole potential $\Lp$ itself is real.
Since the spherical harmonics functions satisfy $Y_l^{-m} = (Y_l^m)^*$,
the coefficients also satisfy $z_{l,-m} = (z_{lm})^*$,
where the asterisk denotes the complex conjugate.
It is therefore enough to compute the coefficients for $m \geq
0$.
The expansion \eqref{eq:dipole-expansion} is in fact valid also
at the point where the ball in Figure~\ref{fig:expansion-centre}
touches $\tilde{\Gamma}$ (where $r_x = r_\text{QBX}$), as
established in \cite{epstein13}.
This means that the expansion can be used both
for offsurface evaluation (in the interior of the ball,
where the double layer potential is nearly singular) as
well as onsurface evaluation (at the point on $\tilde{\Gamma}$
closest to $\vec{c}$, where the double layer potential is singular).

Relation~\eqref{eq:dipole-relation} allows us to express
$\vec{\Dp}$ using four dipole potentials with densities
\begin{equation}
  \label{eq:density-conversion}
  \begin{array}{r@{\:}l}
    \vec{\rho}^{(j)} &= q_j \vec{n} + n_j \vec{q}, \qquad j=1,2,3, \\[4pt]
    \vec{\rho}^{(4)} &= y_k q_k \vec{n} + y_k n_k \vec{q}.
  \end{array}
\end{equation}
Each of the four dipole potentials is expanded using
\eqref{eq:dipole-expansion}, with coefficients given by
\eqref{eq:dipole-qbx-coeffs}, which together with
\eqref{eq:dipole-relation} provides a local expansion of
the Stokes double layer potential.

In practice, the expansion \eqref{eq:dipole-expansion} must be
truncated, which is done at $l = l_\text{max} = p_\text{QBX}$.
This results in the approximation
\begin{equation}
  \label{eq:truncated-expansion}
  \Lp[\tilde{\Gamma}, \vec{\rho}^{(j)}](\vec{x}) \approx
  \Lp^\text{QBX}[\tilde{\Gamma}, \vec{\rho}^{(j)}](\vec{x}) = \sum_{l=0}^{p_\text{QBX}}
  \sum_{m=-l}^l r_x^l Y_l^{-m} (\theta_x, \varphi_x) z^h_{lm,j},
  \qquad j=1,2,3,4.
\end{equation}
The coefficients $z^h_{lm,j} = z_{lm}^h[\vec{\rho}^{(j)}]$ are
here computed using the upsampled quadrature rule introduced in
section~\ref{sec:upsamp-special-quad}, with upsampling factor
$\kappa = \kappa_\text{QBX}$. Upsampling is needed since the
integrand in \eqref{eq:dipole-qbx-coeffs} becomes quite
peaked for large $l$. However, the cost of upsampling can be
entirely hidden in a precomputation step, as explained in
section~\ref{sec:qbx-precomp}. The number of coefficients that
needs to be computed in \eqref{eq:truncated-expansion} for each
$j$ is
\begin{equation}
  \label{eq:num-coeffs}
  N_\text{QBX} = \frac{(p_\text{QBX}+1)(p_\text{QBX}+2)}{2},
\end{equation}
which takes into account that only coefficients with $m \geq 0$
need to be computed directly.

If the expansion \eqref{eq:dipole-expansion} is absolutely
convergent, the terms must decay in magnitude as $l \to \infty$. The
size of the terms can be estimated using the bound
\begin{equation}
  \label{eq:qbx-term-bound}
  \abs{ \sum_{m=-l}^l r_x^l Y_l^{-m} (\theta_x, \varphi_x)
  z_{lm}} \leq
  r_\text{QBX}^l \sqrt{\frac{2l+1}{4\pi}}
  \gp{\sum_{m=-l}^l \abs{z_{lm}}^2}^{1/2},
\end{equation}
where we have used the fact that \cite[Eq.~(3.36)]{epstein13}
\begin{equation}
  \sum_{m=-l}^l \absi{Y_l^m(\theta,\varphi)}^2 =
  \frac{2l+1}{4\pi}.
\end{equation}
For a single dipole expansion such as
\eqref{eq:truncated-expansion} with a fixed $j$, the bound
\eqref{eq:qbx-term-bound} with $z_{lm} = z_{lm,j}^h$ provides a
way to estimate the decay of the terms and thus the truncation
error of the truncated expansion. It is however not directly
applicable to the Stokes double layer potential, which involves
derivatives of dipole potentials as seen in
\eqref{eq:dipole-relation}. To estimate the truncation error for
the Stokes double layer potential, we instead evaluate the error
directly in the grid points, as explained in section~\ref{sec:param1-qbx}.

In summary, to compute
$\vec{\Dp}[\tilde{\Gamma},\vec{q}](\vec{x})$ using QBX, the
density $\vec{q}$ is first upsampled to a finer grid with
upsampling factor $\kappa_\text{QBX}$ and then converted into
four dipole densities using \eqref{eq:density-conversion}. From
these, four sets of dipole coefficients $z^h_{lm,j}$ are computed
using the direct quadrature rule on the refined grid. The coefficients
are used to evaluate the dipole potentials
\eqref{eq:truncated-expansion}, from which the Stokes double
layer potential $\vec{\Dp}$ can be computed using
\eqref{eq:dipole-relation}. Note that the derivatives with
respect to $\vec{x}$ in \eqref{eq:dipole-relation} can be computed
analytically.

Since QBX can be used for both onsurface and offsurface
evaluation, it is useful to introduce one expansion for
each grid point. For each grid point $\vec{x}_i$ on the boundary,
an expansion centre $\vec{c}_i^{+}$ is thus placed at a
distance~$r_\text{QBX}$ away from the boundary in the normal
direction (i.e.\ in the fluid domain). This expansion centre can
be used to evaluate the double layer potential in a ball touching
that grid point. In practice
the balls of convergence of neighbouring expansion centres will
overlap, and for a given evaluation point the closest expansion
centre is used to evaluate the QBX potential.

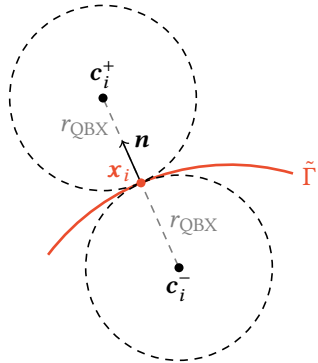
\begin{figure}[b!]
  \centering
  \begin{tikzpicture}
    \node [inner sep=0pt] (c) at (-2.5,1.0) {};
    \node [inner sep=0pt] (cm) at (-1.5,-1.25) {};
    \node [inner sep=0pt] (y) at (-2,-0.125) {};

    \draw [gray, line width=0.6pt, dashed, dash pattern=on 3pt off 2.5pt]
        (c) -- (y) node [midway,left,yshift=1.5mm,xshift=-1pt] {$r_\text{QBX}$};
    \draw [gray, line width=0.6pt, dashed, dash pattern=on 3pt off 2.5pt]
        (cm) -- (y) node [midway,right] {$r_\text{QBX}$};
    \draw [->, line width=0.6pt] (y) -- (-2.25, 0.4375) node [right=1pt] {$\vec{n}$};
    \draw [fred, line width=1pt] (0,0) arc [radius=3cm, start angle=74, end angle=143];
    \node [right, fred] at (0,0) {$\tilde{\Gamma}$};
    \draw [line width=0.6pt, dashed, dash pattern=on 3pt off 2.5pt]
        (c) circle [radius=1.228cm];
    \draw [line width=0.6pt, dashed, dash pattern=on 3pt off 2.5pt]
        (cm) circle [radius=1.228cm];

    \node [fill=black, inner sep=1.2pt, circle] at (c) {};
    \node [above] at (c) {$\vec{c}_i^{+}$};

    \node [fill=black, inner sep=1.2pt, circle] at (cm) {};
    \node [below] at (cm) {$\vec{c}_i^{-}$};

    \node [fill=fred, inner sep=1.2pt, circle] at (y) {};
    \node [left,fred,yshift=1.3mm] at (y) {$\vec{x}_i$};
  \end{tikzpicture}
  \caption{The two expansion centres $\vec{c}_i^{+}$ and
  $\vec{c}_i^{-}$ used for onsurface evaluation in a grid point
  $\vec{x}_i$.}
  \label{fig:double-sided}
\end{figure}

For onsurface evaluation (but not offsurface evaluation), we also
use a second expansion centre~$\vec{c}_i^{-}$ for each grid
point, placed at a
distance $r_\text{QBX}$ away from the boundary in the
\emph{negative} normal direction
(i.e.\ outside the fluid domain), as
shown by Figure~\ref{fig:double-sided}. The reason for this is
that it significantly improves the convergence when solving the boundary
integral equation using GMRES, since the spectrum of the discrete operator
better matches that of the continuous operator, as was noted in
\cite{klockner13,rahimian18,klinteberg16b}.
Note that due to the jump condition \eqref{eq:double-jump},
the correct value of the potential on $\tilde{\Gamma}$ is the
average of the values from the two sides:
\begin{equation}
  \label{eq:average-onsurface-evaluation}
  \vec{\Dp}[\tilde{\Gamma}, \vec{q}](\vec{x}_i) =
  \frac{\vec{\Dp}^{+}[\tilde{\Gamma}, \vec{q}](\vec{x}_i)
  + \vec{\Dp}^{-}[\tilde{\Gamma}, \vec{q}](\vec{x}_i)}{2},
\end{equation}
where $\vec{\Dp}^{+}$ is the limit from the fluid domain and
$\vec{\Dp}^{-}$ is the limit from the other side of
$\tilde{\Gamma}$. While using two expansions may seem to double
the computational cost, the extra cost appears only in the
precomputation step, as described in
section~\ref{sec:qbx-precomp}, and thus does not affect the
cost of evaluation itself.

There are two sources of error in the QBX approximation: the
\emph{truncation error} due to the fact that the expansion in
\eqref{eq:truncated-expansion} is truncated at $l=p_\text{QBX}$,
and the \emph{coefficient error} (called the ``quadrature error'' in
\cite{epstein13,klinteberg16b,klinteberg17}) due to the fact that the coefficients
\eqref{eq:dipole-qbx-coeffs} are computed using a quadrature
rule with finite precision. These two errors are controlled by
the following three QBX parameters:
\begin{itemize}
  \item
    The expansion radius $r_\text{QBX}$, which is the distance
    from the expansion centre to $\tilde{\Gamma}$ and also the radius of the
    ball in which the expansion is valid. Increasing
    $r_\text{QBX}$ makes the truncation error grow since the ball
    of convergence (and hence $r_x$) becomes larger, but the
    coefficient error decreases since the integrand in
    \eqref{eq:dipole-qbx-coeffs} becomes easier to resolve as
    $r_y$ becomes larger.
  \item
    The expansion order $p_\text{QBX}$, which governs the number
    of terms to be included in the sum in
    \eqref{eq:truncated-expansion}. Increasing $p_\text{QBX}$
    makes the truncation error decrease since more terms are
    included, but the coefficient error grows since the integrand
    in \eqref{eq:dipole-qbx-coeffs} is harder to resolve for
    larger $l$.
  \item
    The upsampling factor $\kappa_\text{QBX}$, which governs the
    amount of grid refinement when computing the dipole
    coefficients \eqref{eq:dipole-qbx-coeffs}. Increasing
    $\kappa_\text{QBX}$ makes the coefficient error decrease
    since the resolution of the underlying quadrature rule
    increases.
\end{itemize}
A simple way to decrease both the truncation error and
coefficient error is to increase $p_\text{QBX}$ and
$\kappa_\text{QBX}$ simultaneously while keeping $r_\text{QBX}$
fixed. We will continue to discuss how the QBX parameters should
be selected to achieve a small overall error in
section~\ref{sec:parameters}. For a more in-depth analysis, we
refer to \cite{epstein13} for the truncation error,
\cite{klinteberg17} for the coefficient error, as well as the
summary in \cite[sec.~3.5]{klinteberg16b}.

\section{Global and local QBX}
\label{sec:global-local-qbx}

As was mentioned in section~\ref{sec:intro-related-work}, a QBX
method can be either \emph{(fully) global}, \emph{particle-global} or
\emph{local}, the difference being which part of the boundary
(i.e.\ which source points) to include when forming the local
expansion. Here we use particle-global QBX for particles and local QBX
for walls. In essence, the difference between the three variants is
what $\tilde{\Gamma}$ in section~\ref{sec:expansion-double-layer}
is taken to be:
\begin{itemize}
  \item
    For a fully global QBX method, all grid points on the whole
    boundary are used to form the local expansion, i.e.\
    $\tilde{\Gamma} = \Gamma$.
  \item
    For a particle-global QBX method, all grid points on a single
    particle are used to form the expansion, i.e.\
    $\tilde{\Gamma} = \Gp^{(\alpha)}$, where $\Gp^{(\alpha)}$ is
    the surface of the particle with index~$\alpha$.
  \item
    For a local QBX method, only the grid points which are close
    to the expansion centre are used to form the expansion. In
    our case, we choose $\tilde{\Gamma}$ to be the $N_\text{P}$
    patches of the wall which are closest to the expansion
    centre, as shown in Figure~\ref{fig:wall-shells}. (The
    contribution from patches further away is not included in the
    expansion but computed using direct quadrature.)
\end{itemize}
Note that $\tilde{\Gamma}$ may depend on the location of the
expansion centre to be used, which in turn depends on the
evaluation point. In principle, it is sufficient to let
$\tilde{\Gamma}$ consist of the grid points close to the
expansion centre (i.e.\ local QBX), since that is where the
integrand becomes nearly singular; for grid points further away,
direct quadrature can be used. The reason to extend
$\tilde{\Gamma}$ further is to improve the regularity of the
layer potential that is being expanded, so that
the expansion converges more rapidly. Indeed, in
local QBX, the expanded layer potential consists of the contribution
from a truncated part of the boundary, and may not be very smooth since
$\tilde{\Gamma}$ ends abruptly. However, the larger
$\tilde{\Gamma}$ is, the further away from the ball of
convergence will the edge of $\tilde{\Gamma}$ be, and the less
will it affect the convergence of the expansion. We have observed
that $N_\text{P}=9$ is sufficient for $\tilde{\Gamma}$ for the
walls.

In particle-global QBX, the expanded layer potential has the
contribution from a whole particle, which consists of a closed
and smooth surface, so the potential from it should be smooth. In
a fully global QBX, the expanded layer potential is the global potential,
which is smooth if $\Gamma$ is regular enough. Unlike the
particle-global QBX, the fully global QBX quickly becomes expensive
unless a fast method (such as the FMM) is used to compute the
far-field contribution. Therefore the fully global QBX variant is
not used in this paper.

An advantage of the local and particle-global QBX variants over
the fully global QBX is that, if the individual particles and
walls are rigid, $\tilde{\Gamma}$ is the same (in local
coordinates) for all particles or patches of the same shape, even
if they have different orientations. This makes precomputation
possible, which we shall return to in
section~\ref{sec:qbx-precomp}. Another advantage is that
expansion centres can be placed without regards to other
particles or walls, since each expansion contains only the
contribution from a single particle or wall segment. In a fully
global QBX method, each ball of convergence must be completely
outside \emph{all} particles and walls, which would complicate the
placement of the expansion centres.

\section{Precomputation for QBX}
\label{sec:qbx-precomp}

The mapping given by \eqref{eq:density-conversion} and the
discrete version of \eqref{eq:dipole-qbx-coeffs}, which takes the
double layer density~$\vec{q}$ on~$\tilde{\Gamma}$ and returns
the dipole coefficients \smash{$z^h_{lm,j}$} for a single expansion
centre $\vec{c}_i$, is a linear function of $\vec{q}$ and can therefore be
represented by a matrix $\mat{M}_i$. This matrix is of size
$4N_\text{QBX} \times 3\tilde{N}$, where $N_\text{QBX}$
is given by \eqref{eq:num-coeffs} and $\tilde{N}$ is the number
of grid points on $\tilde{\Gamma}$ (before upsampling). There is
one such matrix $\mat{M}_i$ for every expansion centre, and it
depends only on the geometry~$\tilde{\Gamma}$, its discretization
and the location of the expansion centre in the local coordinates
of $\tilde{\Gamma}$. For a rigid geometry~$\tilde{\Gamma}$, such
as in our case, the matrix $\mat{M}_i$ can therefore be precomputed
and stored.

Note that the upsampling factor $\kappa_\text{QBX}$ is effectively
``hidden'' in this precomputation step: upsampling influences
the computation of $\mat{M}_i$ since $\vec{q}$ is upsampled before
being inserted into \eqref{eq:density-conversion}, but it has no
effect on the size of $\mat{M}_i$, which is set by the
discretization of $\tilde{\Gamma}$ prior to upsampling.
Therefore, upsampling does not affect the computational
complexity of the method once $\mat{M}_i$ has been precomputed.

The matrix $\mat{M}_i$ which computes the coefficients $z^h_{lm,j}$
is used for \emph{offsurface} evaluation, when the evaluation
point is not known beforehand; the coefficients can then be used to
evaluate the expansion at any evaluation point within the ball of
convergence. For \emph{onsurface} evaluation, i.e.\ evaluation at
one of the grid points of the boundary, the evaluation point
itself is known beforehand and precomputation can be taken even
further. In fact, the mapping that takes the expansion
coefficients to the value of the potential $\vec{\Dp}[\tilde{\Gamma},
\vec{q}](\vec{x}_i)$, given by \eqref{eq:truncated-expansion} and
\eqref{eq:dipole-relation}, is also linear and can therefore be
represented by a matrix $\mat{S}_i$. This allows us to compute a
matrix $\mat{R}_i = \mat{S}_i \mat{M}_i$ which maps the density
$\vec{q}$ on $\tilde{\Gamma}$ directly to the value of the double
layer potential $\vec{\Dp}$ at one of the grid points --
effectively representing a set of target-specific quadrature
weights for every grid point. The matrix $\mat{R}_i$ is of size $3
\times 3\tilde{N}$ and there is one such matrix for each grid
point $\vec{x}_i$ on the boundary. Precomputing the $\mat{R}_i$ matrix hides
not only $\kappa_\text{QBX}$ but also $p_\text{QBX}$.

Since two expansion centres are used for onsurface evaluation, as
the reader may recall from Figure~\ref{fig:double-sided}, there
are actually two $\mat{R}_i$ matrices for each grid point:
$\mat{R}_i^{+}$ and $\mat{R}_i^{-}$, associated with $\vec{c}_i^{+}$
and $\vec{c}_i^{-}$, respectively. From
\eqref{eq:average-onsurface-evaluation}, it is clear that these
matrices can be combined as
\begin{equation}
  \mat{R}_i = \frac{\mat{R}_i^{+} + \mat{R}_i^{-}}{2}
\end{equation}
to form a single matrix $\mat{R}_i$ for each grid point. This way,
the extra cost of using two expansions is completely hidden in
the precomputation step.

\begin{figure}[h!]
  \centering
  \includegraphics[scale=0.9]{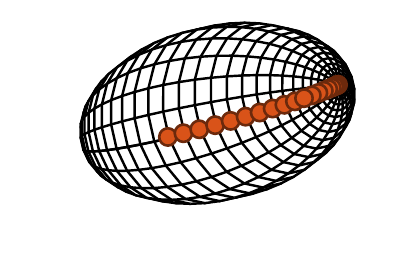}%
  \includegraphics[scale=0.9]{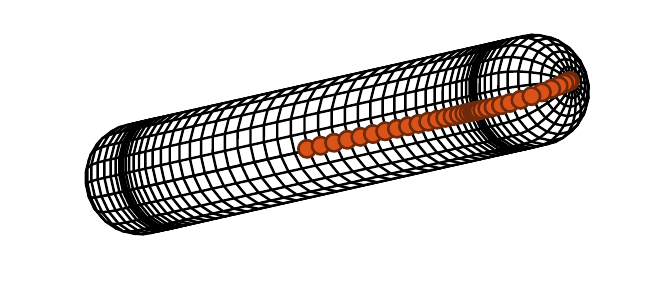}%
  \caption{The matrices $\mat{M}_i$ and $\mat{R}_i$ need only be
  stored for the grid points along half a line of longitude,
  here indicated with red dots.}
  \label{fig:precomp-symm-particle}
\end{figure}

For the particles, the axisymmetry can be used to vastly reduce
the amount of computations and storage needed to precompute the
matrices $\mat{M}_i$ and $\mat{R}_i$. In fact, due to reflective
symmetry, it suffices to compute $\mat{R}_i$ for the
$n_\theta/2$ grid points ($n_1 + n_2/2$ grid points for rod
particles) shown in Figure~\ref{fig:precomp-symm-particle}, and $\mat{M}_i$
for the corresponding expansion centres.
The matrices for all other grid points and their expansion
centres are then calculated using rotations and
reflections, as in \cite{klinteberg16b}. Note that if several
particles of the same shape appear in a simulation, the
precomputation only needs to be done for one such particle.

For a wall geometry with uniform patch size, as in our case, the
geometry has a discrete translational symmetry for offsets equal
to the patch size, due to periodicity. This means that the geometry
``looks'' exactly the same seen from any patch of the wall, and
it is therefore enough to precompute the $\mat{M}_i$ and
$\mat{R}_i$ matrices for the $n_1 n_2$ grid points and
corresponding expansion centres of a single patch of the wall. In
this case $\tilde{\Gamma}$ consists of that patch and its
$N_\text{P}-1$ closest neighbours, as indicated in
Figure~\ref{fig:wall-shells} for $N_\text{P}=9$.

\pagebreak
\vspace*{-1.2cm}
\chapter{Periodicity and fast methods}
\label{sec:periodicity}

Up to this point we have not taken periodicity into account in
the description of the mathematical formulation and its
discretization; it is now time to remedy this. We will here give
the details of the periodic formulation indicated in
Figure~\ref{fig:periodicity}~(b), and in particular focus on how the
special quadrature methods are combined with the fast summation
method used for the periodic problem.

Consider a primary cell with side lengths $\vec{B} = (B_1, B_2,
B_3)$ which is replicated periodically in all three spatial
directions. The flow field is then periodic, i.e.\
$\vec{u}(\vec{x}) = \vec{u}(\vec{x} + \vec{k} \cdot \vec{B})$ for
any $\vec{k} \in \Integer^3$. This changes the boundary integral
formulation introduced in section~\ref{sec:bif} in the following
way: The layer potential~$\vec{\Dp}$ and completion
flow~$\vec{\Vp}^{(\alpha)}$ which appear in the flow field
\eqref{eq:dist-flow} and in the fundamental boundary integral
equation \eqref{eq:proto-bie} are replaced by their periodic
counterparts $\vec{\Dp}^\text{3P}$ and
$\vec{\Vp}^{(\alpha),\text{3P}}$.
These are defined as infinite sums over the periodic lattice,
i.e.
\begin{equation}
  \label{eq:periodic-sums}
  \vec{\Dp}^\text{3P}[\Gamma, \vec{q}](\vec{x}) =
  \sum_{\vec{k} \in \Integer^3} \vec{\Dp}[\Gamma,
  \vec{q}](\vec{x} + \vec{k} \cdot \vec{B}
  ), \qquad
  \vec{\Vp}^{(\alpha),\text{3P}}[\vec{F}, \vec{\tau}](\vec{x}) =
  \sum_{\vec{k} \in \Integer^3} \vec{\Vp}^{(\alpha)}[\vec{F}, \vec{\tau}]
  (\vec{x} + \vec{k} \cdot \vec{B}).
\end{equation}
These sums converge slowly, and their value depends on the order
of summation, so they cannot be computed using direct summation.
We compute them using the \emph{Spectral Ewald (SE)} method
\cite{lindbo10,lindbo11}, a fast Ewald summation method based on
the fast Fourier transform (FFT). The SE method is described in
detail for the stokeslet in \cite{lindbo10}, for the stresslet in
\cite{klinteberg14} and for the rotlet in \cite{klinteberg16a},
and has been combined with QBX previously in
\cite{klinteberg16b}. In the SE method, each of
the periodic sums in \eqref{eq:periodic-sums} is split into two
parts: the \emph{real-space part}, which decays fast and can
therefore be summed directly in real space; and the
\emph{Fourier-space part}, which is smooth and therefore decays
fast in Fourier space.

No special treatment is needed for the completion flow
$\vec{\Vp}^{(\alpha),\text{3P}}$ since the evaluation point is
never close to the singular points (which are inside the particle),
so the SE method as described in \cite{lindbo10,klinteberg16a} is
used without modification. For the double layer potential
$\vec{\Dp}^\text{3P}$, special quadrature is needed so SE must be
combined with QBX and the upsampled quadrature rule. How this is
done is described below.

The periodic sum for the double layer potential can be written
explicitly as
\begin{equation}
  \label{eq:stresslet-sum}
  \Dp_i^\text{3P}[\Gamma, \vec{q}](\vec{x}) =
  \sum_{\vec{k} \in \Integer^3} \int_\Gamma
  T_{ijl}(\vec{x}+\vec{k}\cdot\vec{B}-\vec{y}) q_j(\vec{y}) n_l(\vec{y})
  \, \D S(\vec{y}).
\end{equation}
The stresslet $\mat{T}$ is split into two parts
\begin{equation}
  \label{eq:stresslet-split}
  T_{ijl}(\vec{r}) = T_{ijl}^\text{R}(\vec{r}; \xi) +
  T_{ijl}^\text{F}(\vec{r}; \xi),
\end{equation}
where $\mat{T}^\text{R}$ is the real-space part and
$\mat{T}^\text{F}$ is the Fourier-space part. The Ewald parameter
$\xi$ is a positive number which is used to balance the decay of
$\mat{T}^\text{R}$ in real space and the decay of the Fourier
coefficients of $\mat{T}^\text{F}$ (a larger value of $\xi$ makes
the real-space part decay faster and the Fourier-space part decay
slower, thus shifting computational work into Fourier space). In the
split that we use, $\mat{T}^\text{R}$ is given by \cite{klinteberg14,klinteberg16b}
\begin{align}
  \notag
  T_{ijl}^\text{R}(\vec{r}; \xi) &= -\frac{2}{r^4} \gp{
    \frac{3}{r} \op{erfc}(\xi r)
    + \frac{2\xi}{\sqrt{\pi}}(3 + 2\xi^2r^2 - 4\xi^4r^4)
    \E^{-\xi^2r^2}
  } r_i r_j r_l \\
  \label{eq:explicit-realspace-stresslet}
  &\phantom{=\:\,}
  + \frac{8\xi^3}{\sqrt{\pi}} (2-\xi^2r^2)
  \E^{-\xi^2r^2}
  (\delta_{ij}r_l + \delta_{jl} r_i + \delta_{li} r_j),
\end{align}
where $r = \absi{\vec{r}}$. The Fourier-space part is simply given
by $\mat{T}^\text{F} = \mat{T} - \mat{T}^\text{R}$.
Inserting \eqref{eq:stresslet-split} into
\eqref{eq:stresslet-sum} splits the periodic double layer potential into two parts
$\vec{\Dp}^\text{3P} = \vec{\Dp}^\text{3P,R} +
\vec{\Dp}^\text{3P,F}$, where
\begin{align}
  \label{eq:realspace-sum}
  \Dp_i^\text{3P,R}[\Gamma, \vec{q}](\vec{x}; \xi) &=
  \sum_{\vec{k} \in \Integer^3} \int_\Gamma
  T_{ijl}^\text{R}(\vec{x}+\vec{k}\cdot\vec{B}-\vec{y}; \xi) q_j(\vec{y}) n_l(\vec{y})
  \, \D S(\vec{y}), \\
  \label{eq:fourier-sum}
  \Dp_i^\text{3P,F}[\Gamma, \vec{q}](\vec{x}; \xi) &=
  \sum_{\vec{k} \in \Integer^3} \int_\Gamma
  T_{ijl}^\text{F}(\vec{x}+\vec{k}\cdot\vec{B}-\vec{y}; \xi) q_j(\vec{y}) n_l(\vec{y})
  \, \D S(\vec{y}).
\end{align}
The singularity of the stresslet is completely
transferred to $\mat{T}^\text{R}$, while $\mat{T}^\text{F}$ is
nonsingular~\cite{klinteberg16b}. The Fourier-space part
\eqref{eq:fourier-sum} is computed using FFTs as described in
appendix~\ref{app:streamlines} and \cite{klinteberg14}. The
real-space potential \eqref{eq:realspace-sum} is evaluated in
real space, and requires special quadrature due to the singularity of
$\mat{T}^\text{R}$, much as in the free-space setting. Note that since
$\mat{T}^\text{R}(\vec{r};\xi)$ decays fast as $\absi{\vec{r}}
\to \infty$ it can be neglected for $\absi{\vec{r}} > r_\text{c}$,
where $r_\text{c}$ is called the \emph{cutoff radius}. We can
thus change the integration domain in \eqref{eq:realspace-sum} to
$\Gamma^\star = \Gamma^\star(\vec{x}, \vec{k}; r_\text{c})
= \seti{\vec{y} \in \Gamma : \absi{\vec{x} + \vec{k} \cdot
\vec{B} - \vec{y}} \leq r_\text{c}}$ and approximate
\begin{equation}
  \label{eq:truncated-realspace-sum}
  \Dp_i^\text{3P,R}[\Gamma,\vec{q}](\vec{x}; \xi) \approx
  \Dp_i^{\text{3P,R}\star}[\Gamma,\vec{q}](\vec{x}; \xi) =
  \sum_{\vec{k} \in \Integer^3} \int_{\Gamma^\star}
  T_{ijl}^\text{R}(\vec{x}+\vec{k}\cdot\vec{B}-\vec{y}; \xi) q_j(\vec{y}) n_l(\vec{y})
  \, \D S(\vec{y}).
\end{equation}
The error of this approximation is determined by the product $\xi
r_\text{c}$ as described in \cite{klinteberg14}. Rather than
deriving a new QBX expansion from scratch for the real-space part
$\vec{\Dp}^{\text{3P,R}\star}$, we reuse the expansion of the
total layer potential $\vec{\Dp}$ from section~\ref{sec:qbx-detail}.
To be able to do this, we insert $\mat{T}^\text{R} = \mat{T} -
\mat{T}^\text{F}$ into \eqref{eq:truncated-realspace-sum} to get
\begin{align}
  \label{eq:qbx-realspace-sum-1}
  \Dp_i^{\text{3P,R}\star}[\Gamma, \vec{q}](\vec{x}; \xi) &=
  \sum_{\vec{k} \in \Integer^3} \int_{\Gamma^\star}
  T_{ijl}(\vec{x}+\vec{k}\cdot\vec{B}-\vec{y}; \xi) q_j(\vec{y}) n_l(\vec{y})
  \, \D S(\vec{y}) \\
  \label{eq:qbx-realspace-sum-2}
  &\phantom{=} -
  \sum_{\vec{k} \in \Integer^3} \int_{\Gamma^\star}
  T_{ijl}^\text{F}(\vec{x}+\vec{k}\cdot\vec{B}-\vec{y}; \xi) q_j(\vec{y}) n_l(\vec{y})
  \, \D S(\vec{y}).
\end{align}
Note that the integration domain $\Gamma^\star$ ensures that both
of these sums have few terms since $r_\text{c}$ should be small
-- typically smaller than the size of the periodic cell. The integral in
\eqref{eq:qbx-realspace-sum-1} represents the total layer
potential from $\Gamma^\star$ and can thus be computed using the
combined special quadrature method from
section~\ref{sec:upsamp-special-quad}, with truncation at
$r_\text{c}$. The integral in \eqref{eq:qbx-realspace-sum-2} is
computed using direct quadrature, which is possible since
$\mat{T}^\text{F}$ is nonsingular.

As in the free-space setting in section~\ref{sec:upsamp-special-quad},
the evaluation point~$\vec{x}$ is classified into one of three
regions (see Figure~\ref{fig:combined-strategy}). The potential
is evaluated using \eqref{eq:truncated-realspace-sum} in the
direct quadrature region and
\eqref{eq:qbx-realspace-sum-1}--\eqref{eq:qbx-realspace-sum-2} in
the other two regions.
The reader may wonder why we in the upsampled quadrature region
do not simply evaluate \eqref{eq:truncated-realspace-sum} using
upsampled quadrature. The reason is that \eqref{eq:qbx-realspace-sum-2}
evaluated using the same quadrature method as in the
Fourier-space part---i.e.\ direct quadrature---is needed to
cancel discretization errors in the latter.%
\footnote{
  Due to the nonlocal nature of the Fourier transform, the
  Fourier-space part must be computed using the same quadrature
  method everywhere; here we use direct quadrature. Another
  possibility would be to use upsampled quadrature, but then
  upsampling would need to be done for \emph{all} evaluation
  points, not only those in the upsampled quadrature region. In
  that case one might want to remove the direct quadrature region
  altogether and use only upsampled quadrature and QBX.
}
These discretization errors may be larger than the SE error
tolerance and are caused by the fact that
$\mat{T}^\text{F}(\vec{r};\xi)$, while nonsingular, tends to
become slightly peaked for small $\absi{\vec{r}}$, i.e.\ close to
$\Gamma^\star$. Cancellation prevents these errors from
influencing the error of the full method.

Another important point to note is that $r_\text{c}$ must be
chosen large enough so that no special quadrature is needed for
the total potential when $\absi{\vec{r}} > r_\text{c}$. This is
because, for $\absi{\vec{r}} > r_\text{c}$, the total potential
is equal to the Fourier-space part, which is always computed
using direct quadrature. Thus, $r_\text{c}$ must be at least as
large as the distance from $\Gamma$ to the direct quadrature
region.

\FloatBarrier
\chapter{Parameter selection}
\label{sec:parameters}

In this section, we develop our strategy for selecting the
parameters of the combined special quadrature, i.e.\ upsampled
quadrature and QBX, when evaluating the Stokes double layer
potential $\vec{\Dp}$. We assume that a discretization of the
geometry is given, with sufficient resolution for the density to
be well-resolved and the direct quadrature to achieve a given
error tolerance $\etol$ at a given distance (sufficiently far
away) from all surfaces.
The goal is to select quadrature parameters for each particle and
wall so that the error tolerance $\etol$ is achieved also in the
upsampled quadrature region and QBX region. Of course, there are
many different ways to choose the parameters, some resulting in
higher computational efficiency than others. Here, we do not aim
to optimize the efficiency; instead, our focus is on achieving
the given error tolerance at an acceptable (albeit not optimal)
computational cost.

\begin{figure}[ht!]
  \centering
  \includegraphics{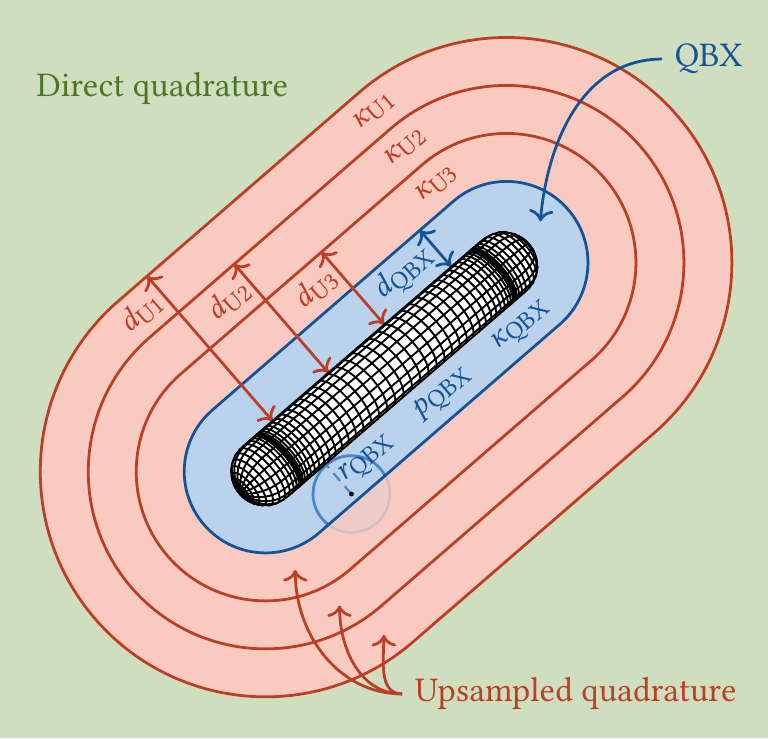}
  \caption{Parameters of the combined special quadrature, here
  for a rod with $N_\text{U}=3$ upsampled quadrature regions
  (each geometrical object has its own set of corresponding
  parameters). Note that while the ball of convergence for QBX
  can be larger than the QBX region, the expansion is used only
  inside the QBX region.}
  \label{fig:parameters}
\end{figure}

The parameters that must be selected are shown in
Figure~\ref{fig:parameters}. Note that we allow for multiple
upsampled quadrature regions with different upsampling factors
$\kappa_\text{U$i$}$, in order to gradually increase the
upsampling closer to the surface. Due to the precomputation
scheme for QBX, using QBX may in fact be faster than using
upsampled quadrature with the same upsampling factor. Therefore,
the QBX region may extend further away from the surface than the
expansion centre (i.e.\ $d_\text{QBX}$ may be larger than
$r_\text{QBX}$, but of course not larger than $2r_\text{QBX}$).

The parameters to be selected are as follows:
\begin{itemize}
  \item
    The threshold distances $d_\text{U$i$}$ for the upsampled
    quadrature regions, $i=1,2,\ldots, N_\text{U}$, and the
    threshold distance $d_\text{QBX}$ for the QBX region. These
    determine at what distance from the surface each region
    starts. If $d_\Gamma(\vec{x})$ is the distance from the
    evaluation point $\vec{x}$ to the surface $\Gamma$, then
    $\vec{x}$ belongs to the $i$th upsampled quadrature region if
    \begin{equation}
      d_\text{U$i$} \geq d_\Gamma(\vec{x}) \geq \begin{cases}
        d_\text{U$(i+1)$} & \text{if $i < N_\text{U}$} \\
        d_\text{QBX} & \text{if $i=N_\text{U}$},
      \end{cases}
    \end{equation}
    and $\vec{x}$ belongs to the QBX region if $d_\text{QBX} \geq
    d_\Gamma(\vec{x}) \geq 0$. Each of these distances should be
    chosen so that the error does not exceed the tolerance in the
    region further away from the surface (for example,
    $d_\text{U1}$ is selected based on the direct quadrature
    error).

  \item
    The upsampling factors $\kappa_\text{U$i$}$ for the upsampled
    quadrature regions, $i=1,2,\ldots,N_\text{U}$. These should
    be increasing, i.e.\ $\kappa_\text{U1} < \kappa_\text{U2} <
    \ldots < \kappa_\text{U$N_\text{U}$}$. The upsampling factor
    $\kappa_\text{U$i$}$ determines the distance
    $d_\text{U$(i+1)$}$ at which the next region must begin,
    which we will come back to in section~\ref{sec:param-example1}.

  \item
    The QBX upsampling factor $\kappa_\text{QBX}$, which controls the
    amount of upsampling used when computing the coefficients in
    \eqref{eq:dipole-qbx-coeffs}, and thus the QBX coefficient
    error as mentioned in section~\ref{sec:expansion-double-layer}.
    It should be chosen large enough so that the coefficient
    error and the truncation error are balanced. The upsampling
    factor influences the QBX precomputation time (which grows
    like $O(\kappa_\text{QBX}^2)$), but not the
    size of the precomputed $\mat{M}_i$ and $\mat{R}_i$ matrices, and
    thus not the evaluation time.

  \item
    The QBX expansion order $p_\text{QBX}$, which controls the number of
    terms included in the expansion in~\eqref{eq:truncated-expansion},
    and thus the QBX truncation
    error as mentioned in section~\ref{sec:expansion-double-layer}.
    It should be chosen so that the truncation error is below the
    error tolerance everywhere in the QBX region. The expansion
    order affects the size of the $\mat{M}_i$ matrix used for
    offsurface evaluation (which grows like $O(p_\text{QBX}^2)$),
    but not that of the $\mat{R}_i$ matrix
    used for onsurface evaluation. It should be noted that as
    $p_\text{QBX}$ increases, the upsampling factor
    $\kappa_\text{QBX}$ must also increase since higher-order
    coefficients are harder to resolve.

  \item
    The QBX expansion radius $r_\text{QBX}$, which affects both the
    coefficient error and the truncation error, but neither the
    precomputation time nor the evaluation time
    directly. It should typically be chosen as small as possible,
    since this speeds up the convergence of the expansion in
    \eqref{eq:truncated-expansion} so that $p_\text{QBX}$ can be
    chosen small. On the other hand, a very small $r_\text{QBX}$
    means that the upsampling factor $\kappa_\text{QBX}$ must be
    large, since the expansion centre moves closer to the
    surface.

    However, in our implementation the primary restriction on
    $r_\text{QBX}$ is that it must be large enough for the balls
    of convergence to cover the QBX region sufficiently well. In
    general, $r_\text{QBX}$ should not be smaller than the
    distance from one expansion centre to the next, to ensure a
    good coverage. Since we have one expansion centre per grid
    point, we require that $r_\text{QBX}$ should not be smaller
    than the grid spacing.%
    \footnote{An alternative would be to introduce more expansion
    centres to maintain the coverage of the QBX region as
    $r_\text{QBX}$ decreases below the grid spacing. Doing so
    would also increase the amount of work and storage needed in
    the precomputation step.}
    Letting $h$ be some measure of the grid spacing (for example
    the largest distance between neighbouring grid points on the
    surface), it is useful to consider the ratio $r_\text{QBX}/h$
    when selecting parameters. As noted in \cite{klinteberg16b},
    this has the advantage that if $r_\text{QBX}/h$ is kept fixed
    during refinement of the original grid, then the coefficient
    error is constant, assuming that the upsampling factor
    $\kappa_\text{QBX}$ is also fixed. We will therefore consider
    $r_\text{QBX}/h$ in the rest of this section.
\end{itemize}

Unfortunately, no general error estimates are available in
three dimensions for the quadrature rules that we use here. The
parameters must therefore be selected based on numerical
experiments, and we present a strategy for doing so here. The
idea is to start from the outermost upsampled quadrature region
(U1) and then proceed inwards towards the surface of the
particle or wall, determining the parameters in the following
order:
\begin{enumerate}
  \item
    Threshold distances and upsampling factors for the upsampled
    quadrature regions,
  \item
    The QBX parameters $r_\text{QBX}/h$ and $d_\text{QBX}$,
  \item
    The QBX parameters $p_\text{QBX}$ and $\kappa_\text{QBX}$.
\end{enumerate}
The process must be repeated for each type of particle and wall
to be used. We develop the strategy in the context of a specific
rod particle in section~\ref{sec:param-example1}; a summary of
the parameter selection strategy in the general case follows in
section~\ref{sec:param-general-summary}. In
sections~\ref{sec:param-example2} and \ref{sec:param-example3} we
apply the strategy to two more examples (a rod with a higher
aspect ratio and a plane wall).

In order to estimate the error during the parameter selection
process, we apply a constant density~$\tilde{\vec{q}}$ such that
$\absi{\tilde{\vec{q}}} = 1$ to the surface and evaluate the
stresslet identity~\eqref{eq:stresslet-identity}.%
\footnote{
  Since the computation of the layer potential is a linear
  function of $\vec{q}$, the error will scale with
  $Q = \max_{\vec{x} \in \Gamma} \absi{\vec{q}(\vec{x})}$.
  In particular, if the maximum error is $\etol$ when $Q = 1$,
  the maximum error will be $\alpha \etol$ when the density is
  multiplied by a constant $\alpha$.
}
This may seem like an overly simple test case since both the
density and the solution are constant.
However, the QBX expansions are not of the constant double layer
potential $\vec{\Dp}$ itself, but of the four dipole potentials
$\Lp$ defined by \eqref{eq:dipole-relation}, and these are
\emph{not} constant. In practice, the stresslet identity seems to
provide a decent test case for both direct quadrature, upsampled
quadrature and QBX, as shown by the results in
section~\ref{sec:res2-convergence} where the density is not
constant.

The Spectral Ewald parameters $\xi$ and $r_\text{c}$ will not be
discussed at length here, but we note that the requirement that
no special quadrature be needed for $\absi{\vec{r}} > r_\text{c}$
implies that $r_\text{c}$ must be at least as large as
$d_\text{U1}$. The Spectral Ewald error is determined by the
product $\xi r_\text{c}$ in real space and $\xi h_\text{F}$ in
Fourier space, where $h_\text{F}$ is the grid spacing of the
uniform grid used for the Fourier-space part (see
appendix~\ref{app:streamlines}). Given a tolerance $\etol$, the
parameters $\xi$, $r_\text{c}$ and $h_\text{F}$ must satisfy the
system $\xi r_\text{c} = A(\etol)$, $\xi h_\text{F} = B(\etol)$,
where $A$ and $B$ are known functions. This leaves one degree of
freedom which can be used to minimize the computational cost,
albeit under the constraint $r_\text{c} \geq d_\text{U1}$. For a
general discussion on the selection of Spectral Ewald parameters,
including the functions $A$ and $B$, we refer to
\cite{klinteberg14} for the stresslet, \cite{lindbo10} for the
stokeslet, and \cite{klinteberg16a} for the rotlet.

\section{Introductory example: a rod particle with low aspect ratio}
\label{sec:param-example1}

\begin{figure}[h!]
  \centering\small
  \begin{minipage}{0.38\textwidth}%
    \centering
    \includegraphics[scale=0.4]{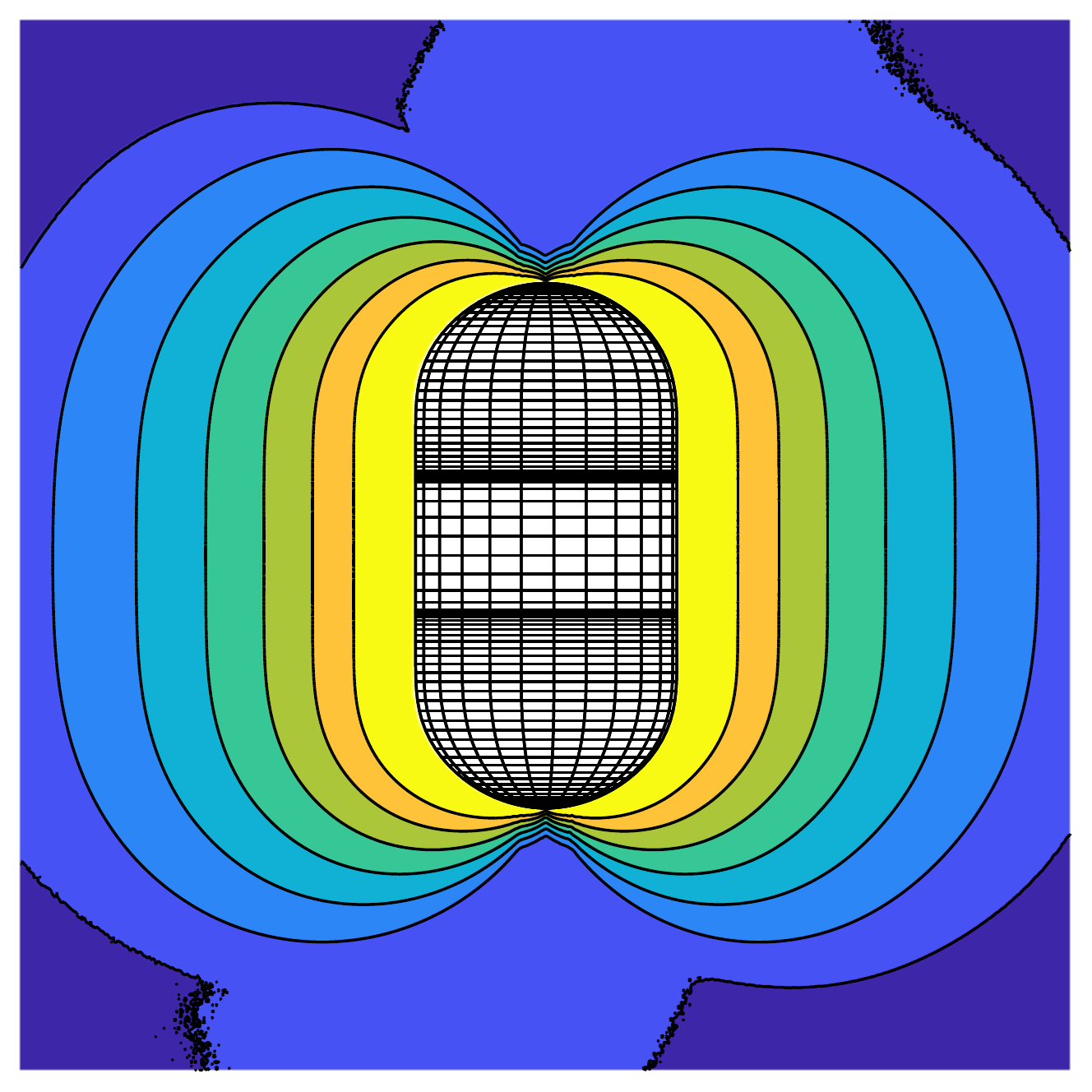}\\%
    (a) Direct quadrature error, slice~1
  \end{minipage}%
  \begin{minipage}{0.38\textwidth}%
    \centering
    \includegraphics[scale=0.4]{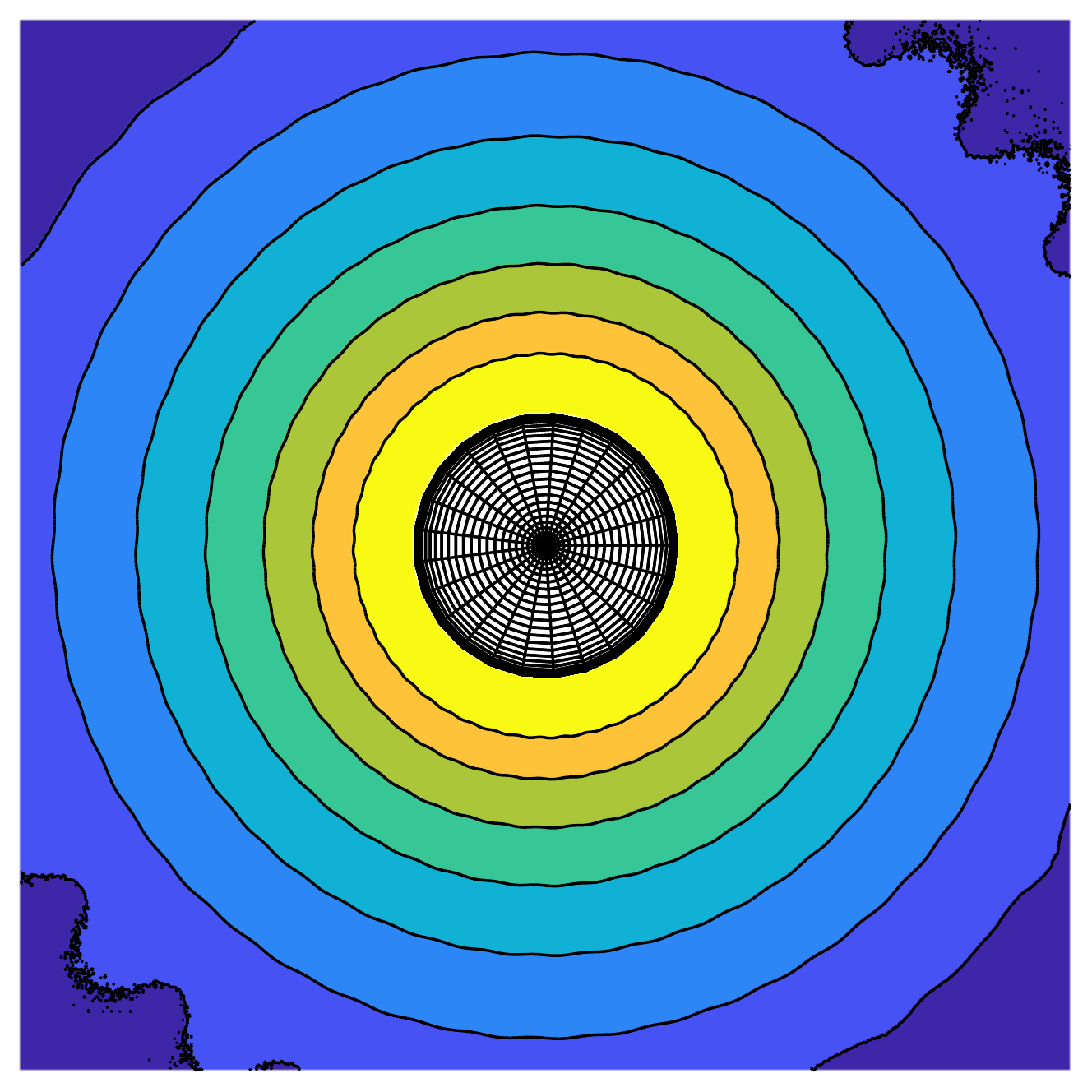}\\%
    (b) Direct quadrature error, slice~2
  \end{minipage}%
  \begin{minipage}{0.1\textwidth}%
    \centering
    \includegraphics[scale=1]{fig/colorbar.pdf}\\%
    \vphantom{(c)}
  \end{minipage}%
  \caption{Error in two perpendicular slices (a) and (b) through the rod
  particle, when evaluating the stresslet identity~\eqref{eq:stresslet-identity}
  using direct quadrature in free space.}
  \label{fig:param1-direct}
\end{figure}

In this first example, we consider a rod particle of length~$L=2$
and radius~$R=0.5$ (i.e.\ aspect ratio~2), shown in
Figure~\ref{fig:param1-direct}. The grid used for the direct
quadrature has parameters $n_1 = 40$, $n_2 = 10$ and $n_\varphi =
25$ (introduced in section~\ref{sec:direct-quad-particles}),
for a total of 2250 grid points. To give an idea of the
error associated with the direct quadrature, we apply the constant
density $\tilde{\vec{q}} = (1, 1, 1) / \sqrt{3}$ to the particle
surface and compute the stresslet identity~\eqref{eq:stresslet-identity}
using direct quadrature in two planes intersecting the particle.
The absolute error in these planes is shown in
Figure~\ref{fig:param1-direct}.

\begin{figure}[h!]
  \centering\small
  \begin{minipage}[b]{0.5\textwidth}%
    \centering
    \includegraphics[scale=0.8,trim=0mm 9mm 0mm 5mm,clip]%
    {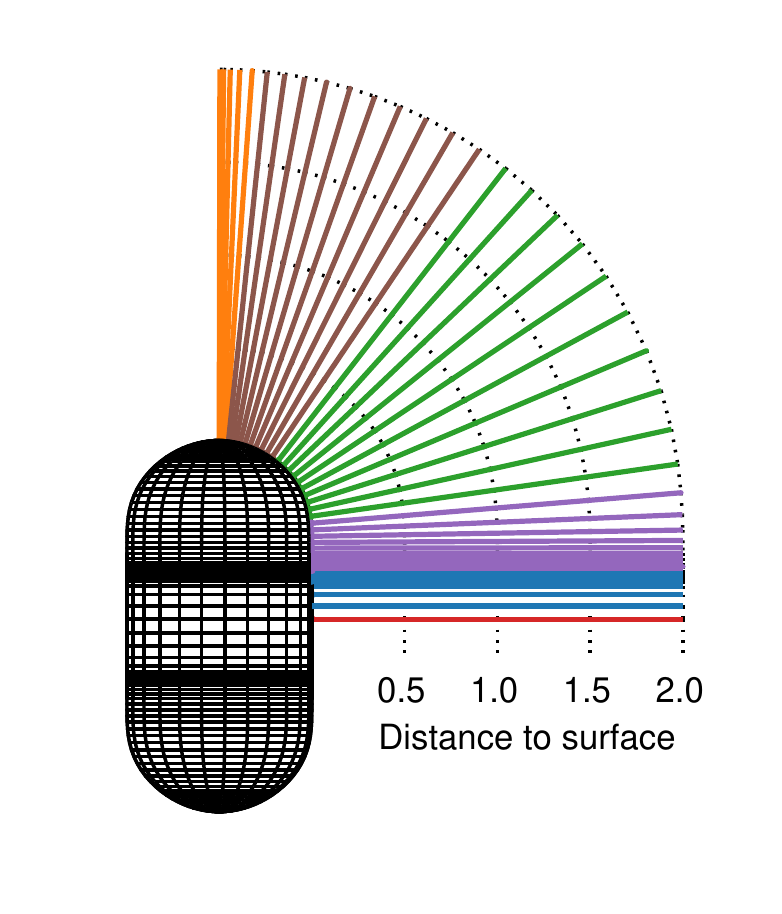}\\
    (a) Lines along which the error is plotted
  \end{minipage}%
  \begin{minipage}[b]{0.5\textwidth}%
    \centering
    \includegraphics[scale=0.8]{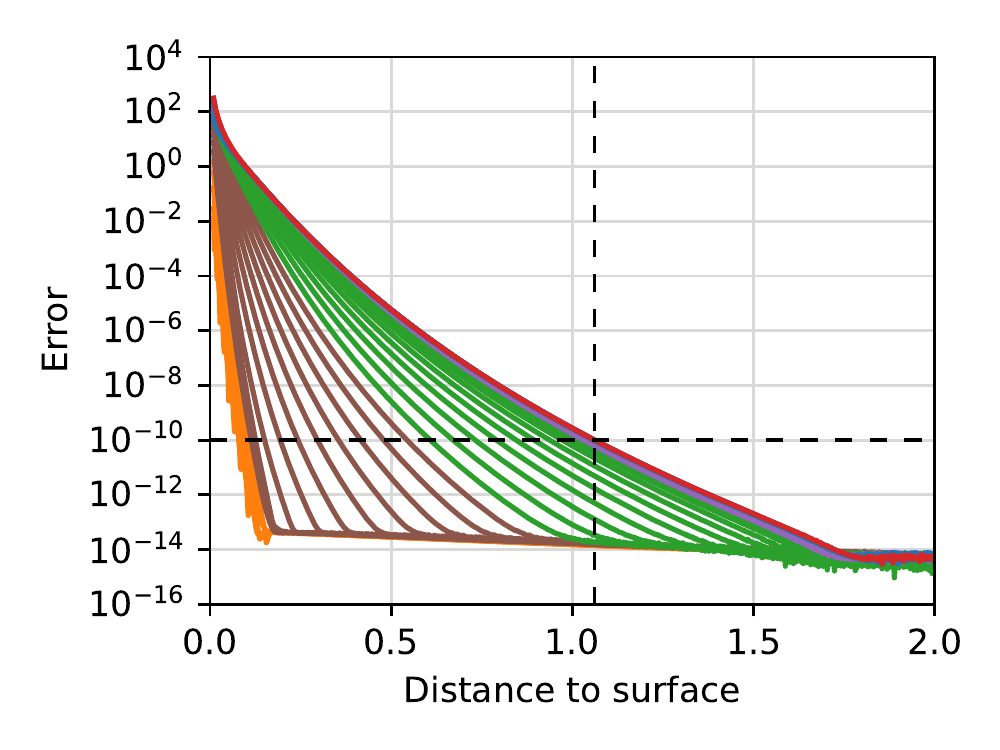}\\
    \hspace*{6mm}(b) Direct quadrature error along the lines
  \end{minipage}%
  \caption{Error of the direct quadrature along 45 normal lines
  centred on grid points of the particle surface. In (a), the
  normal lines are shown coloured in groups of ten. In (b), the
  error along each line is shown with the same colours (here, the
  blue and purple curves are obscured by the red curve). The line
  with the greatest error is coloured red; the error along this
  line reaches the value \num{e-10} at distance 1.061 from the
  surface. (The smallest distance to surface included here is
  0.01.)}
  \label{fig:param1-direct-line}
\end{figure}

To determine how the error varies with the distance to the
surface, we evaluate \eqref{eq:stresslet-identity} along
several normal lines centred on grid points of the particle; due
to the symmetry of the error it is enough to consider the
$n_1+n_2/2 = 45$ lines shown in Figure~\ref{fig:param1-direct-line}~(a).
The error along these lines is shown in
Figure~\ref{fig:param1-direct-line}~(b). Given an error tolerance
$\etol$, the smallest distance at which the error does not exceed
$\etol$ can be determined numerically. This distance is
taken as $d_\text{U1}$. In this example, we will use the error
tolerance $\etol = \num{e-10}$. As indicated in
Figure~\ref{fig:param1-direct-line}~(b), the error reaches
\num{e-10} at $d_\text{U1} = 1.061$; special quadrature must be used
within this distance to the surface. Having established the first
threshold distance $d_\text{U1}$, we now proceed to determine the
rest of the parameters for the upsampled quadrature regions, in
section~\ref{sec:param1-upsamp}.

\subsection{Parameters for the upsampled quadrature regions}
\label{sec:param1-upsamp}

\begin{figure}[h!]
  \centering
  \begin{minipage}{0.5\textwidth}%
    \centering
    \includegraphics[scale=0.8]{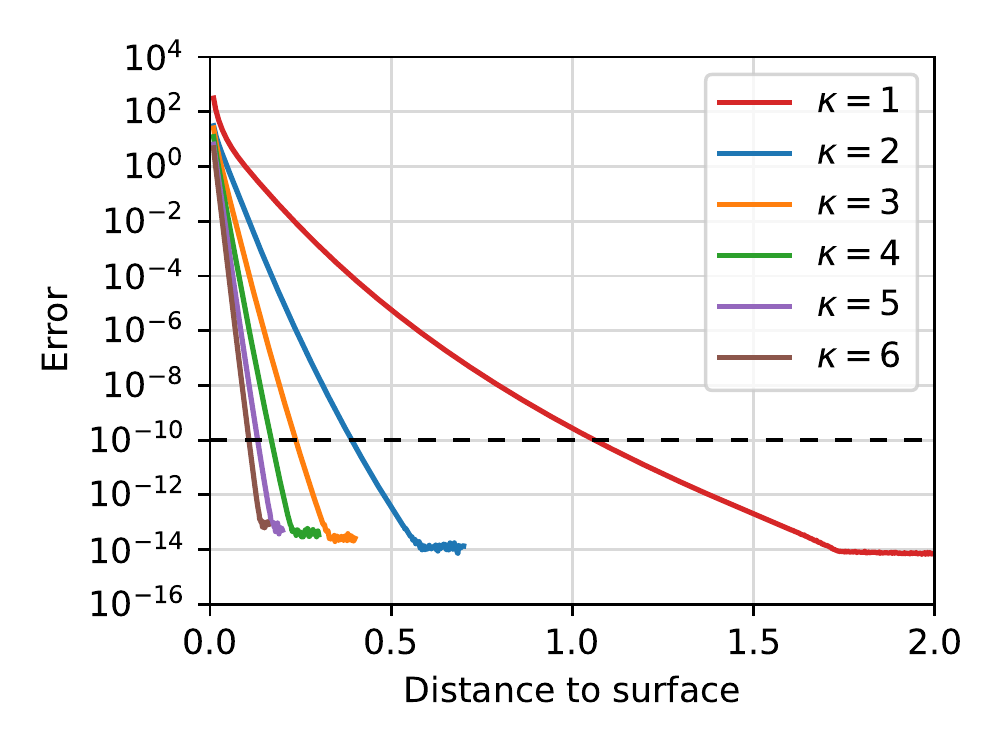}%
  \end{minipage}%
  \caption{Maximal stresslet identity error along any of the
  lines shown in Figure~\ref{fig:param1-direct-line}~(a), for
  upsampled quadrature with different upsampling factors
  $\kappa$. (The smallest distance to surface included here is
  0.01.)}
  \label{fig:param1-upsamp-line}
\end{figure}

For the sake of simplicity we will always choose the upsampling
factors to be $\kappa_\text{U$i$} = i+1$, meaning that the first
upsampling factor will be $\kappa_\text{U1} = 2$, the next will
be $\kappa_\text{U2} = 3$ and so on (this may not be the optimal
strategy with regard to computational cost, but recall that our
goal is not to optimize for computational efficiency).
In order to determine the threshold distance $d_\text{U$i$}$ of
every upsampled quadrature region, we repeat the investigation from
Figure~\ref{fig:param1-direct-line} for different upsampling
factors $\kappa = 1, 2, 3, \ldots$, computing the stresslet
identity error as a function of the distance to the surface for
each upsampling factor. The maximal error at each distance is
shown in Figure~\ref{fig:param1-upsamp-line} ($\kappa=1$
corresponds to Figure~\ref{fig:param1-direct-line}). The
threshold distance $d_\text{U$(i+1)$}$ is now taken as the
distance at which the error curve corresponding to $\kappa =
\kappa_\text{U$i$}$ intersects the error tolerance~$\etol$ ($i=1,
2, \ldots$). For instance, in this case the curve corresponding to
$\kappa=\kappa_\text{U1} = 2$ intersects $\etol = \num{e-10}$
around $0.391 = d_\text{U2}$.

This procedure sets all of the parameters for the upsampled
quadrature, as shown in Table~\ref{tab:param1-upsamp}. However,
at some point we must switch from upsampled quadrature to QBX,
which is determined by the QBX threshold distance~$d_\text{QBX}$.
Selecting $d_\text{QBX}$ will also fix the number of upsampled
quadrature regions~$N_\text{U}$. We will determine $d_\text{QBX}$
together with the other QBX parameters in
section~\ref{sec:param1-qbx}.

\begin{table}[h!]
  \centering
  \caption{Parameters for the upsampled
  quadrature regions, tolerance \num{e-10}.}
  \begin{tabular}{cccccccc}
    \toprule
    $i$ & 1 & 2 & 3 & 4 & 5 & 6 & $\cdots$ \\
    \midrule
    $\kappa_\text{U$i$}$ & 2 & 3 & 4 & 5 & 6 & 7 & $\cdots$ \\
    \midrule
    $d_\text{U$i$}$ & 1.061 & 0.391 & 0.237 & 0.169 & 0.132 &
    0.108 & $\cdots$ \\
    \bottomrule
  \end{tabular}
  \label{tab:param1-upsamp}
\end{table}%

\subsection{Parameters for the QBX region}
\label{sec:param1-qbx}

\begin{figure}[h!]
  \centering\small
  \begin{minipage}[b]{0.35\textwidth}%
    \centering
    \includegraphics[scale=0.4]{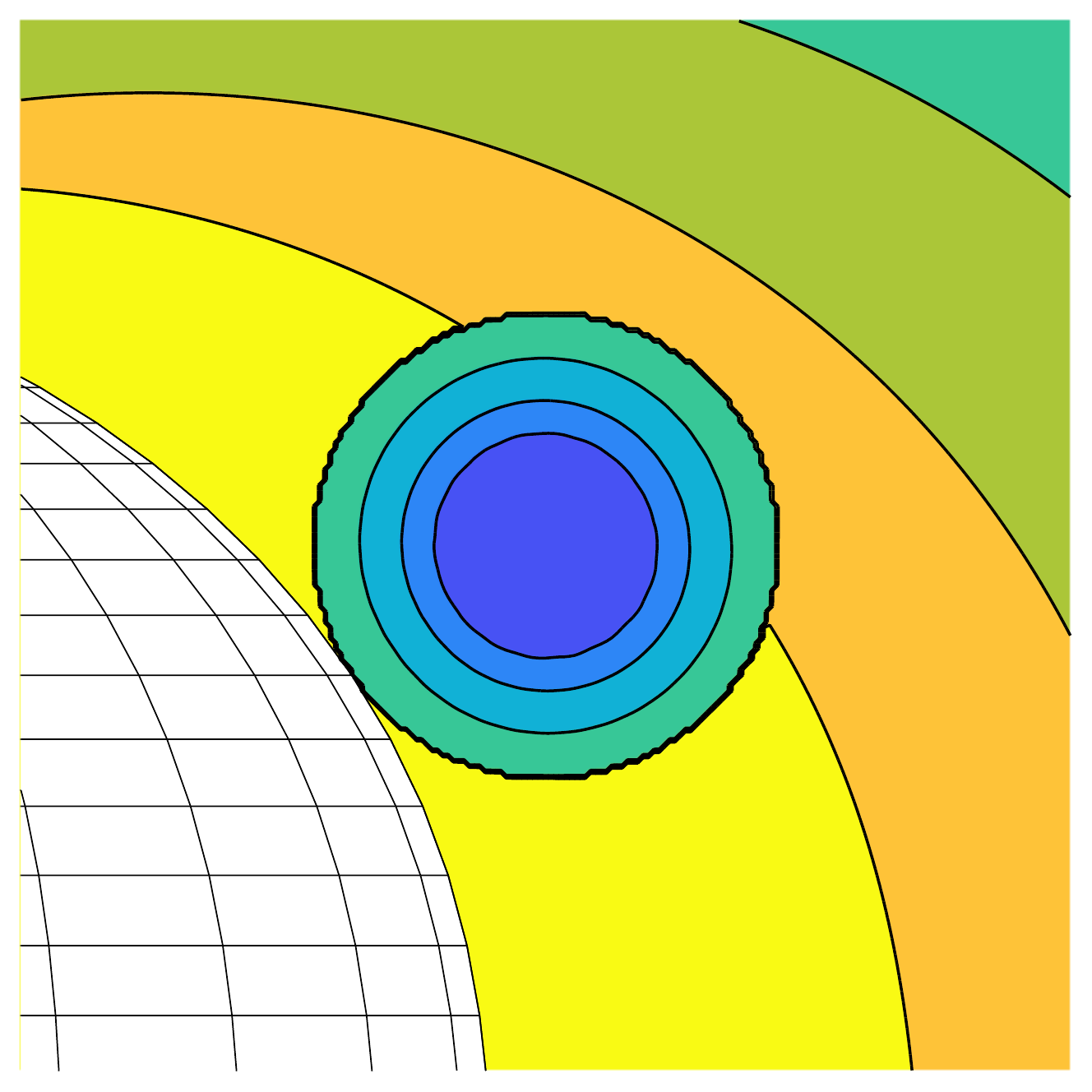}\\%
    (a) Offsurface error from a single QBX expansion
  \end{minipage}%
  \begin{minipage}[b]{0.35\textwidth}%
    \centering
    \includegraphics[scale=0.4,trim=0mm 5mm 0mm 0mm,clip]{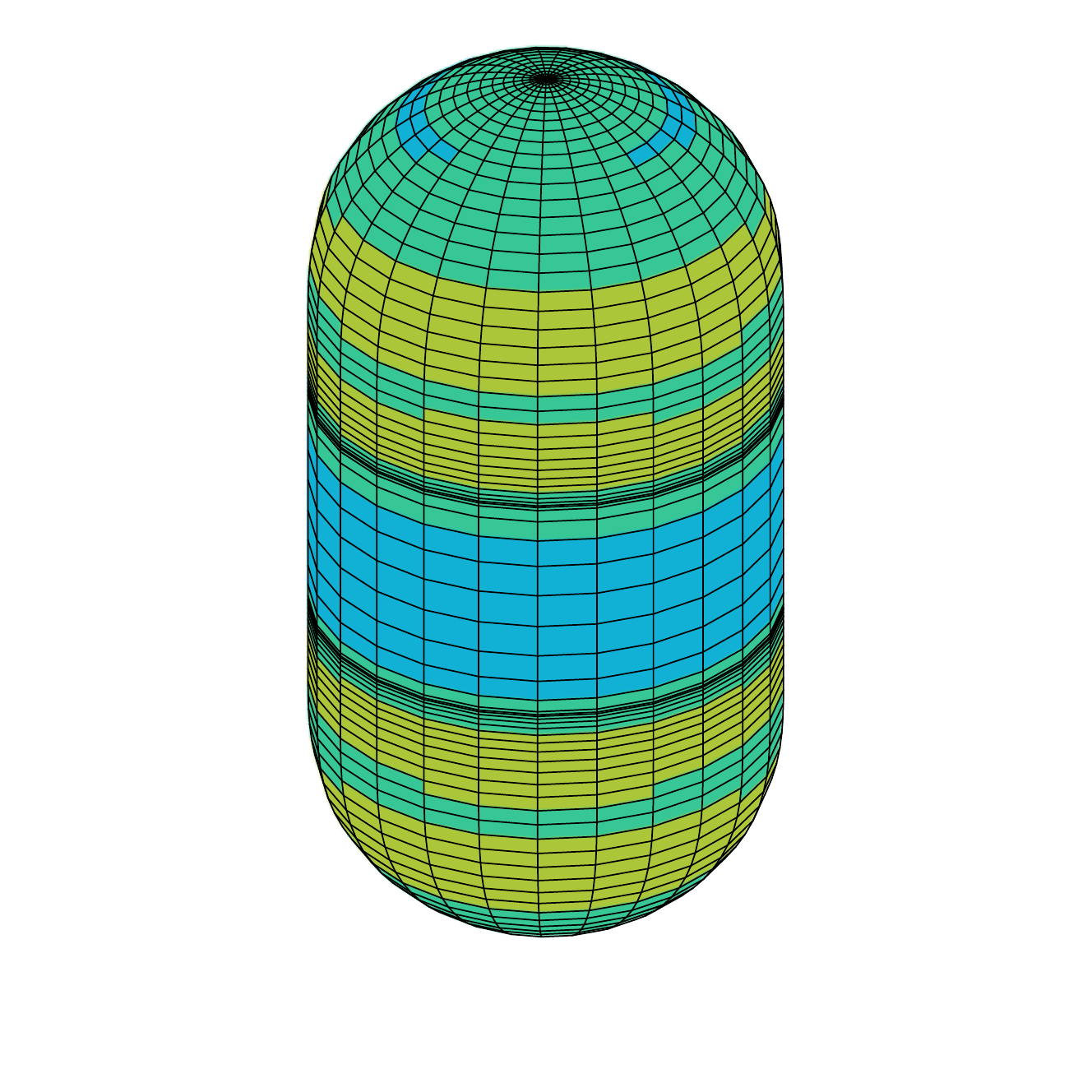}\\%
    (b) Onsurface error\\\vphantom{QBX expansion}
  \end{minipage}%
  \begin{minipage}[b]{0.1\textwidth}%
    \centering
    \includegraphics[scale=1]{fig/colorbar.pdf}\\%
    \vphantom{(c)}\hspace*{0.001pt}\\\vphantom{QBX expansion}\vspace{7mm}
  \end{minipage}%
  \caption[QBX error]{Error when evaluating the stresslet
  identity \eqref{eq:stresslet-identity} using quadrature by
  expansion with the pa\-ram\-e\-ters $r_\text{QBX} = h =
  \pi/25$, $p_\text{QBX} = 18$ and $\kappa_\text{QBX} = 15$:
  (a) in a slice through the particle
  centre (i.e.\ offsurface), using a single QBX expansion centre
  and direct quadrature outside the ball of convergence;
  (b) in the grid points of the particle (i.e.\ onsurface). Note
  in (b) that the error seems to be related to the curvature of
  the boundary; in particular the error is larger in areas
  where the curvature changes, viz.\ in the smooth transition
  from cylinder to cap. (Similar observations related to the
  convexity of the boundary were reported in \cite{barnett14} and
  \cite{klockner13}.)}
  \label{fig:param1-qbx-demonstration}
\end{figure}

To understand how the QBX error behaves, we plot the offsurface
error from a single expansion in
Figure~\ref{fig:param1-qbx-demonstration}~(a). Note that the QBX
parameters used in this figure are not yet selected to achieve the
error tolerance, but meant only to demonstrate the general
behaviour of the error. Since the QBX error is the largest at the
boundary of the ball of convergence (outside this ball the direct
quadrature error is shown in
Figure~\ref{fig:param1-qbx-demonstration}~(a)), it is sufficient
to measure the error at a point on this boundary, for example at
the point where the ball touches the particle. Thus, we measure
the QBX error at all the grid points of the rod -- the onsurface
error -- shown in Figure~\ref{fig:param1-qbx-demonstration}~(b)
for these particular QBX parameters.

For particles, we define the grid spacing $h$ as the distance
between grid points in the azimuthal direction at the equator of
the particle, i.e.\
\begin{equation}
  \label{eq:grid-spacing-particles}
  h = \frac{2\pi R}{n_\varphi}~\text{for rods},
  \hspace{8em}
  h = \frac{2\pi a}{n_\varphi}~\text{for spheroids},
\end{equation}
where $n_\varphi$ is the number of grid points in the azimuthal
direction (as defined in section~\ref{sec:direct-quad-particles}),
$R$ is the radius of the rod and $a$ is the equatorial semiaxis
of the spheroid (which appears in \eqref{eq:spheroid-surface}).
For the rod that we consider in this example, $R=0.5$ and
$n_\varphi = 25$, so $h = \pi/25 \approx 0.1257$.

\begin{figure}[t!]
  \centering\small
  \begin{minipage}[b]{0.5\textwidth}%
    \centering
    \includegraphics[scale=0.8]{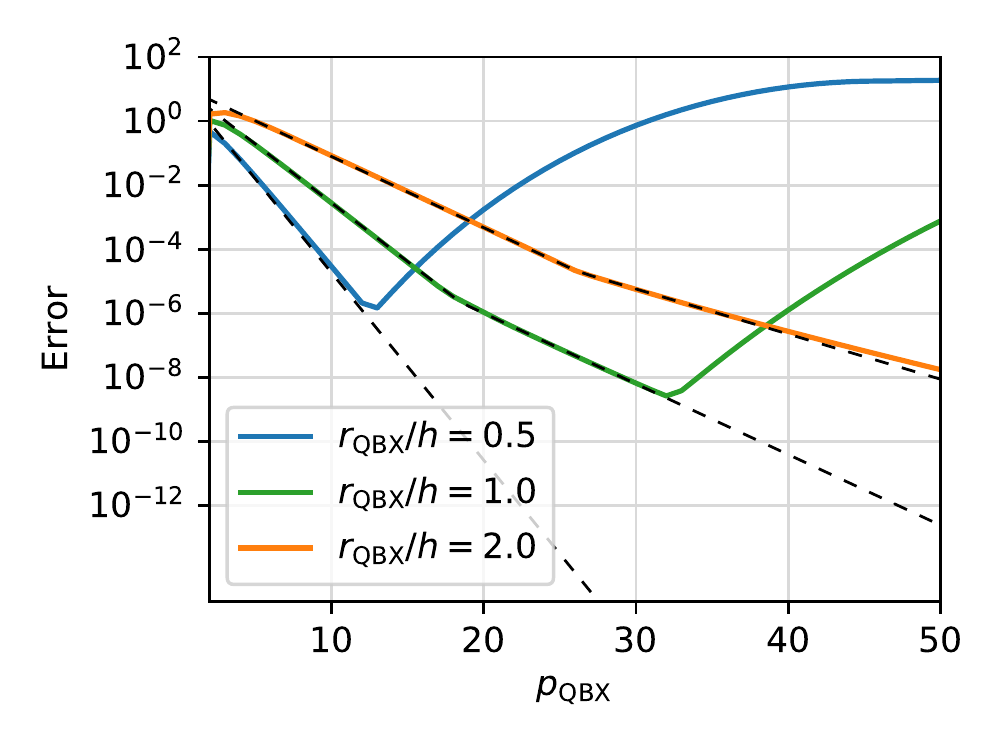}\\%
    \hspace*{6mm}(a) Onsurface QBX error for different $r_\text{QBX}/h$
  \end{minipage}%
  \begin{minipage}[b]{0.5\textwidth}%
    \centering
    \includegraphics[scale=0.8]{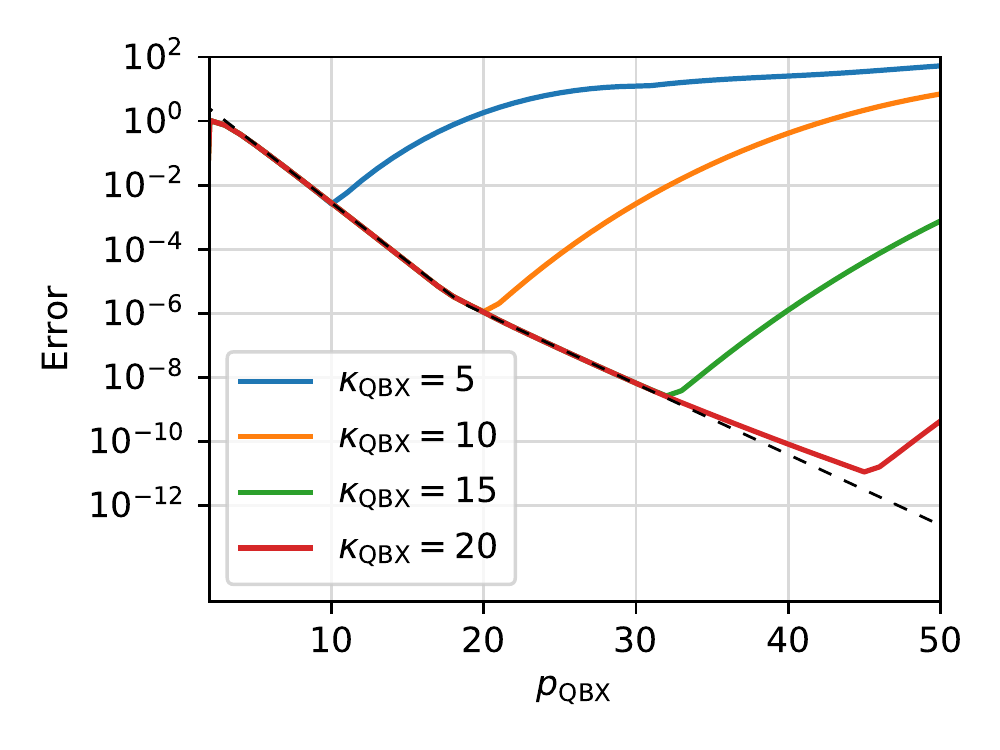}\\%
    \hspace*{6mm}(b) Onsurface QBX error for different
    $\kappa_\text{QBX}$
  \end{minipage}%
  \caption[Maximal onsurface error]%
  {Maximal onsurface error on the rod when evaluating the
  stresslet identity using QBX. In (a) for different
  $r_\text{QBX}/h$ with $\kappa_\text{QBX}=15$ fixed, and in (b)
  for different $\kappa_\text{QBX}$ with $r_\text{QBX}/h=1$
  fixed. Note that each curve has a minimum, which is where the
  truncation error (which decreases as $p_\text{QBX}$ grows) and
  the coefficient error (which increases as $p_\text{QBX}$ grows)
  balance. The dashed lines indicate the experimental truncation
  error estimate~\eqref{eq:param1-truncation-estimate}.}
  \label{fig:param1-qbx-p-kr}
\end{figure}

\begin{figure}[H]
  \centering
  \def\ffR{50}
  \def\ffr{1.3}
  \def\ffd{1.8}
  \def\ffV{3.5}
  \def\ffth{0.6}
  \def\pA{({\ffR*cos(90+\ffV)},{\ffR*sin(90+\ffV)})}
  \def\pB{({\ffR*cos(90-\ffV)},{\ffR*sin(90-\ffV)})}
  \def\pC{({(\ffR+\ffd)*cos(90+\ffV)},{(\ffR+\ffd)*sin(90+\ffV)})}
  \def\pD{({(\ffR+\ffd)*cos(90-\ffV)},{(\ffR+\ffd)*sin(90-\ffV)})}
  \def\pGa{({\ffR*cos(90+3*\ffth)},{\ffR*sin(90+3*\ffth)})}
  \def\pGb{({\ffR*cos(90+\ffth)},{\ffR*sin(90+\ffth)})}
  \def\pGc{({\ffR*cos(90-\ffth)},{\ffR*sin(90-\ffth)})}
  \def\pGd{({\ffR*cos(90-3*\ffth)},{\ffR*sin(90-3*\ffth)})}
  \def\pEa{({(\ffR+\ffr)*cos(90+3*\ffth)},{(\ffR+\ffr)*sin(90+3*\ffth)})}
  \def\pEb{({(\ffR+\ffr)*cos(90+\ffth)},{(\ffR+\ffr)*sin(90+\ffth)})}
  \def\pEc{({(\ffR+\ffr)*cos(90-\ffth)},{(\ffR+\ffr)*sin(90-\ffth)})}
  \def\pEd{({(\ffR+\ffr)*cos(90-3*\ffth)},{(\ffR+\ffr)*sin(90-3*\ffth)})}
  \begin{tikzpicture}[scale=0.8]
    \fill [fblue!20] \pA -- \pC
    arc [radius=\ffR+\ffd, start angle=90+\ffV, end angle=90-\ffV] -- \pB
    arc [radius=\ffR, start angle=90-\ffV, end angle=90+\ffV];

    \draw [line width=1pt, KTHred] \pA
    arc [radius=\ffR, start angle=90+\ffV, end angle=90-\ffV]
    node [anchor=west] {Surface $\Gamma$};
    \node [circle, fill=KTHred, inner sep=1pt] at \pGa {};
    \node [circle, fill=KTHred, inner sep=1pt] at \pGb {};
    \node [circle, fill=KTHred, inner sep=1pt] at \pGc {};
    \node [circle, fill=KTHred, inner sep=1pt] at \pGd {};

    \draw [line width=0.7pt, fblue] \pC
    arc [radius=\ffR+\ffd, start angle=90+\ffV, end angle=90-\ffV];

    \node [circle, fill=black, inner sep=1pt] at \pEa {};
    \node [circle, fill=black, inner sep=1pt] at \pEb {};
    \node [circle, fill=black, inner sep=1pt] at \pEc {};
    \node [circle, fill=black, inner sep=1pt] at \pEd {};
    \draw \pEa circle [radius=\ffr];
    \draw \pEb circle [radius=\ffr];
    \draw \pEc circle [radius=\ffr];
    \draw \pEd circle [radius=\ffr];
    \draw [dashed] \pGa -- \pEa
    node [right, midway, font=\normalsize, xshift=-1pt] {$r_\text{QBX}$};
    \draw [dashed] \pGb -- \pEb;
    \draw [dashed] \pGc -- \pEc;
    \draw [dashed] \pGd -- \pEd;

    \draw [<->, dashed, yshift=-5pt, KTHred] \pGb
    arc [radius=\ffR, start angle=90+\ffth, end angle=90-\ffth]
    node [below=-1pt, midway, font=\normalsize] {$h$};
    \draw [<->, dashed, fblue] \pA -- \pC
    node [midway, left, font=\normalsize, xshift=1pt] {$d_\text{QBX}$};
    \path [fblue] \pB -- \pD node [midway, right] {QBX region};
  \end{tikzpicture}
  \caption{Balls of convergence with the parameters $r_\text{QBX}$ and $d_\text{QBX}$
  and the grid spacing~$h$ marked.}
  \label{fig:param1-qbx-ball-sketch}
\end{figure}
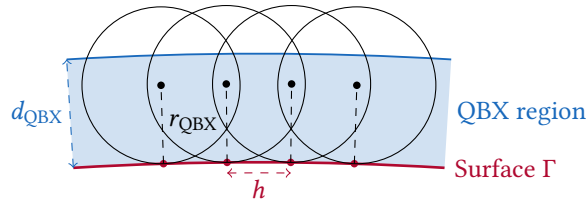

We now focus on selecting the parameters $r_\text{QBX}/h$,
$p_\text{QBX}$ and $\kappa_\text{QBX}$ such that the error is
bounded by $\etol$ in the whole ball of convergence, for all QBX
expansions of the particle. To do this, we consider the maximal
onsurface error as we vary these three parameters, shown in
Figure~\ref{fig:param1-qbx-p-kr}. As seen in
Figure~\ref{fig:param1-qbx-p-kr}~(a), $r_\text{QBX}/h$ should be
chosen as small as possible since this improves the decay
of the truncation error as $p_\text{QBX}$ grows; if
$r_\text{QBX}/h$ is small, $p_\text{QBX}$ can also be chosen
small, which is important since the offsurface evaluation time for
QBX grows as $O(p_\text{QBX}^2)$. On the other hand, as
Figure~\ref{fig:param1-qbx-ball-sketch} shows, $r_\text{QBX}$
must not be too small compared to $h$, since then the balls of
convergence would not overlap properly, and large areas of the
QBX region would not be covered by any ball of convergence.%
\footnote{Some areas of the QBX region will inevitably fall
outside every ball of convergence no matter how large
$r_\text{QBX}$ is. However, these areas are mainly very close to
the surface but not at the grid points, where it is typically not
necessary to evaluate the layer potential.}
For this reason we require that $r_\text{QBX} \geq h$. In fact,
since $r_\text{QBX}$ should be as small as possible, we will
always set $r_\text{QBX} = h$, so that $r_\text{QBX}/h = 1$.

To select $p_\text{QBX}$ and $\kappa_\text{QBX}$, we use the data
shown in Figure~\ref{fig:param1-qbx-p-kr}~(b), which is for
$r_\text{QBX}/h=1$. As can be seen
there, the truncation error is independent of $\kappa_\text{QBX}$
and depends only on $p_\text{QBX}$, so we simply select the
smallest $p_\text{QBX}$ such that the truncation error is below
the tolerance.%
\footnote{
  The dashed curves in Figure~\ref{fig:param1-qbx-p-kr} indicate
  the experimental truncation error estimate
  \begin{equation}
    \label{eq:param1-truncation-estimate}
    e_\text{trunc} \approx \max\!\Big(13 (0.245 \log(\rho) + 0.43)^{p_\text{QBX}},
    ~0.07 (\rho-0.63) (0.175 \log(\rho) + 0.602)^{p_\text{QBX}} \Big),
  \end{equation}
  where $\rho = r_\text{QBX}/h$. This estimate was constructed
  for the rod particle in this particular example by applying
  curve fitting to data from a parameter study similar to that
  shown in Figure~\ref{fig:param1-qbx-p-kr} itself.
  Unfortunately, this experimental estimate is of limited use in
  the parameter selection process since it would have to be
  reconstructed for every new geometrical object (such as a rod
  with a different aspect ratio), while the data used to
  construct it can just as well be used directly to select
  $p_\text{QBX}$.
}
Then we select the smallest $\kappa_\text{QBX}$ (restricted to
multiples of five for convenience) such that the coefficient
error is no larger than the truncation error (i.e.\ such that the
minimum point of the error curve is to the right of the selected
$p_\text{QBX}$). For example, for $\etol=\num{e-10}$, we must
choose $p_\text{QBX} = 40$ and $\kappa_\text{QBX} = 20$.

It remains to choose the threshold distance $d_\text{QBX}$, which
determines the extent of the QBX region as shown in
Figure~\ref{fig:param1-qbx-ball-sketch}. Clearly $d_\text{QBX}$
cannot be larger than $2 r_\text{QBX}$ since then the balls of
convergence would not reach the edge of the QBX region. Even with
$d_\text{QBX} = 2 r_\text{QBX}$, there would be areas in the QBX
region, close to its edge, that would not be inside any ball of
convergence. To mitigate this problem, we introduce a safety
factor $\gamma$, derived in appendix~\ref{app:safety-factor}, and
require that
\begin{equation}
  \label{eq:dQBX-essential-restriction}
  d_\text{QBX} \leq 2 \gamma r_\text{QBX},
\end{equation}
where $\gamma = 0.85$. As long as $d_\text{QBX}$ satisfies
\eqref{eq:dQBX-essential-restriction}, it can be chosen
arbitrarily, in the sense that its value will not affect\linebreak

\noindent
the conformance to the error tolerance, only the computational cost.
We introduce the somewhat arbitrary additional constraint that
$d_\text{QBX} \geq r_\text{QBX}$, and then select $d_\text{QBX}$
as follows: If the interval $[r_\text{QBX}, 2\gamma
r_\text{QBX}]$ contains any threshold distance $d_\text{U$i$}$
for the upsampled quadrature regions, set $d_\text{QBX}$ equal to
the largest $d_\text{U$i$}$ in the interval (i.e.\ the one with
the smallest $i$). Otherwise, set $d_\text{QBX} = r_\text{QBX}$.
In any case, this also sets the number of upsampled quadrature regions
$N_\text{U}$ since the last upsampled quadrature region ends
where the QBX region begins. Our choices here are motivated by
keeping $N_\text{U}$ as low as possible since this reduces the
computational cost, which we will return to in
section~\ref{sec:param1-timing}.

In our current example, $d_\text{QBX}$ should be in the interval
$[r_\text{QBX}, 2 \gamma r_\text{QBX}] \approx [0.1257, 0.2136]$.
As seen in Table~\ref{tab:param1-upsamp}, $d_\text{U4} = 0.169$
is the largest threshold distance in this interval, and thus we
select $d_\text{QBX} = 0.169$ which means that $N_\text{U} = 3$
upsampled quadrature regions are used.

\subsection{Verification of selected parameters}
\label{sec:param1-results}

\begin{figure}[h!]
  \centering\small
  \begin{minipage}[b]{0.5\textwidth}%
    \centering
    \includegraphics[scale=0.95]{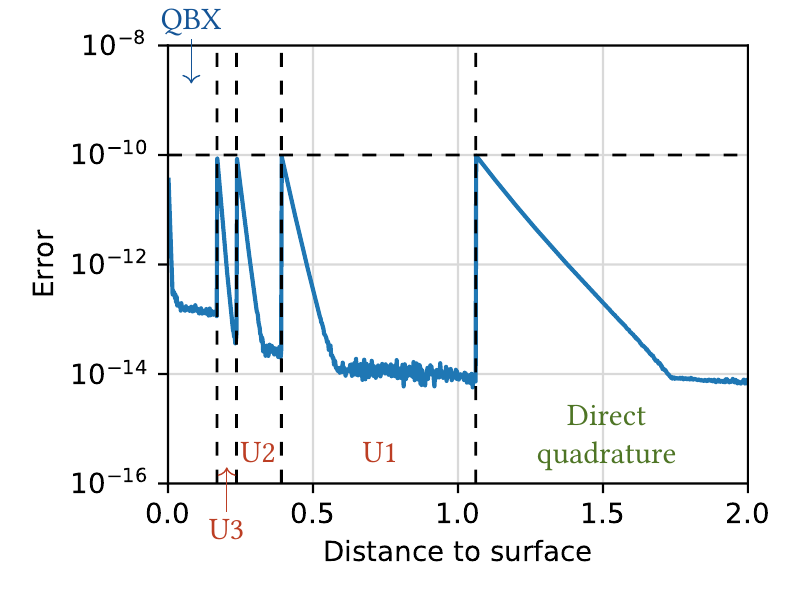}\\[-5pt]%
    \hspace*{6mm}(a) Special quadrature error, tolerance \num{e-10}
  \end{minipage}%
  \begin{minipage}[b]{0.5\textwidth}%
    \centering
    \includegraphics[scale=0.95]{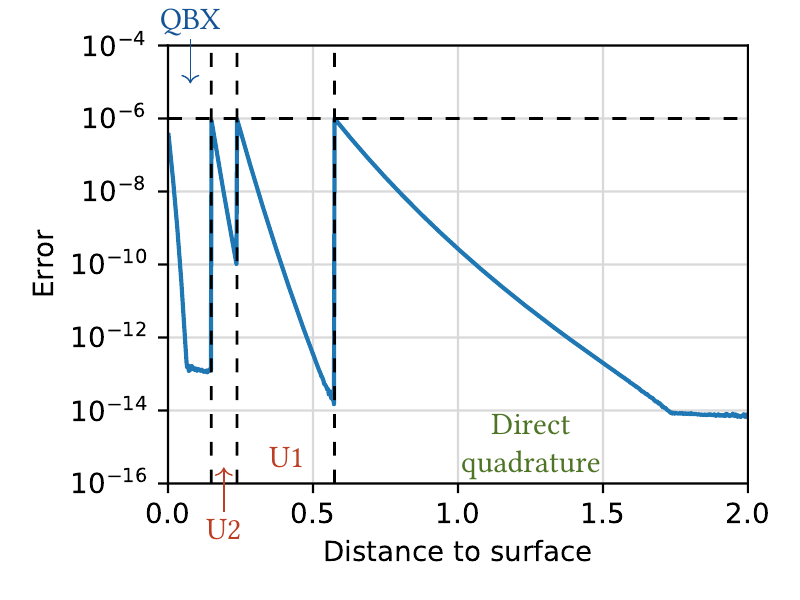}\\[-5pt]%
    \hspace*{6mm}(b) Special quadrature error, tolerance \num{e-6}
  \end{minipage}%
  \caption{Maximal stresslet identity error along any of the
  lines shown in Figure~\ref{fig:param1-direct-line}~(a) as a
  function of the distance to the surface (for 1000 equispaced
  distances in $[0,2]$), using the combined
  special quadrature with (a)~tolerance~\num{e-10} and
  (b)~tolerance~\num{e-6}. The different quadrature regions are
  marked. The largest error is \num{9.709e-11} in (a), and
  \num{9.995e-7} in (b).}
  \label{fig:param1-special-distance}
\end{figure}

\begin{figure}[h!]
  \centering\small
  \begin{minipage}{0.38\textwidth}%
    \centering
    \includegraphics[scale=0.35]{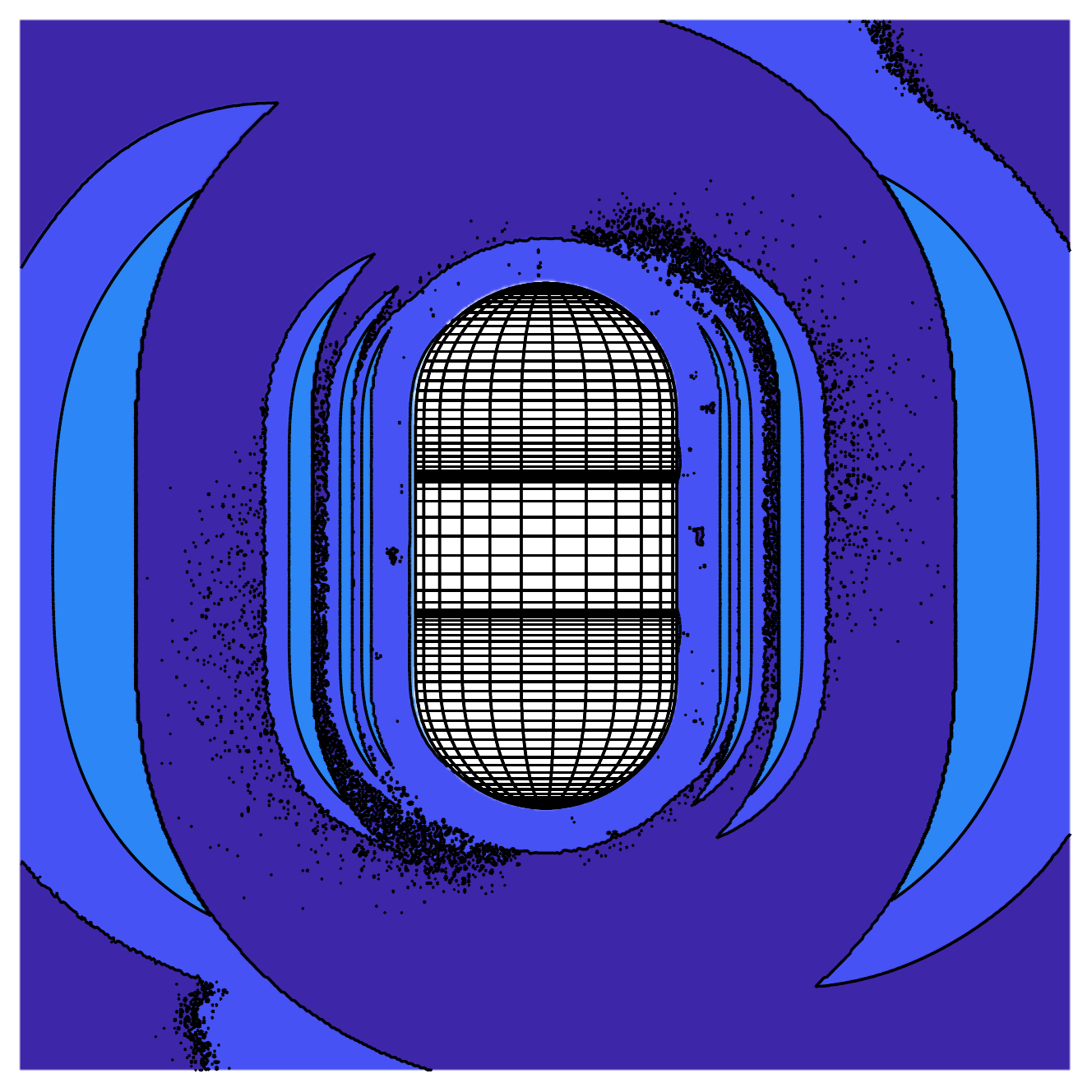}\\%
    (a) Special quadrature error,\\tolerance \num{e-10}, slice~1
  \end{minipage}%
  \begin{minipage}{0.38\textwidth}%
    \centering
    \includegraphics[scale=0.35]{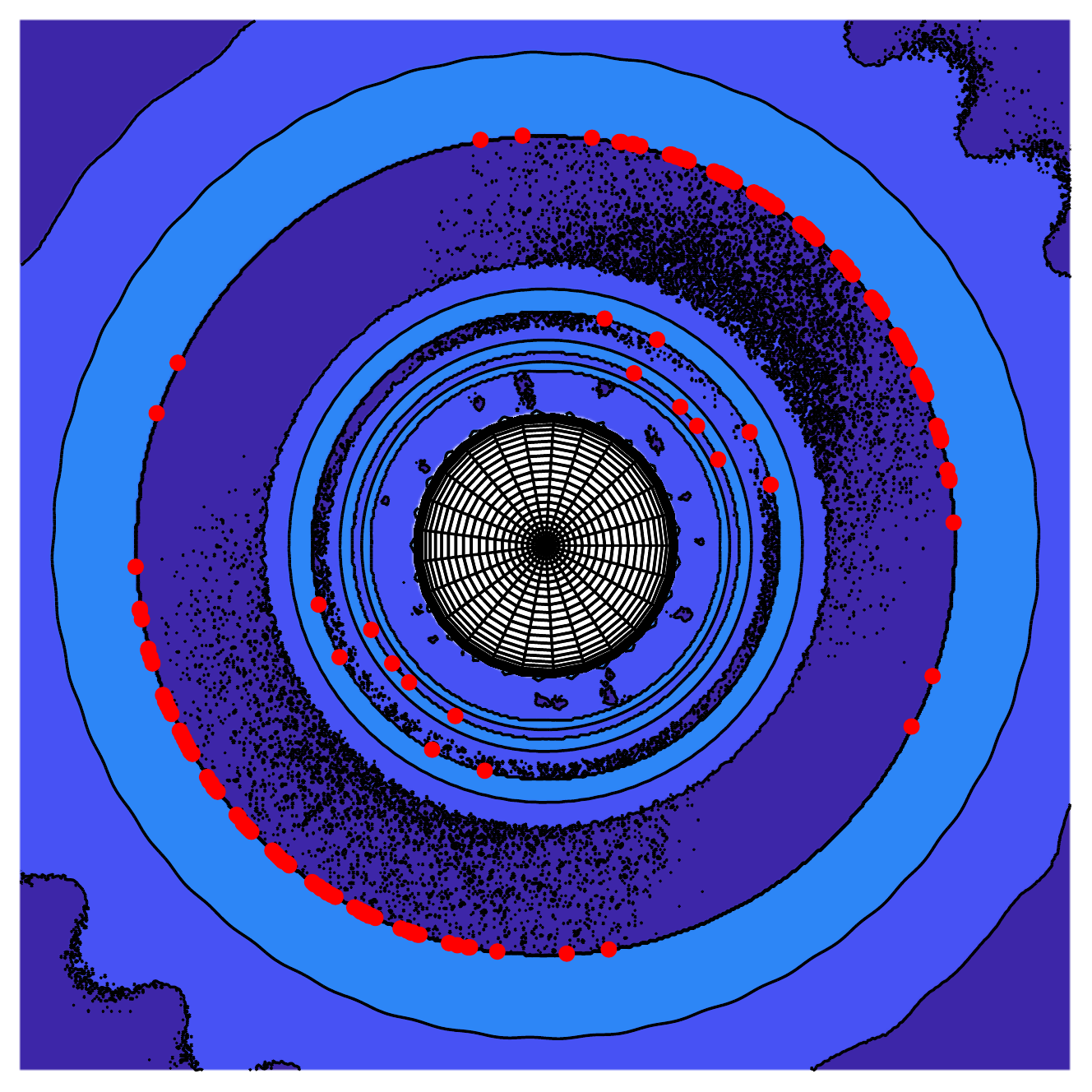}\\%
    (b) Special quadrature error,\\tolerance \num{e-10}, slice~2
  \end{minipage}%
  \begin{minipage}{0.1\textwidth}%
    \centering
    \includegraphics[scale=1]{fig/colorbar.pdf}\\%
    \vphantom{(c)}\hspace*{0.001pt}\\\vphantom{tolerance}
  \end{minipage}%
  \caption{Error in two perpendicular slices (a) and (b) through
  the rod particle, when evaluating the stresslet
  identity~\eqref{eq:stresslet-identity} using combined special quadrature
  with tolerance \num{e-10} in free space. Each slice consists of
  $500 \times 500$ evaluation points. All points in (a) are below the tolerance,
  but 166 points in (b), marked red, are above the tolerance. The largest
  error is \num{9.768e-11} in (a) and \num{1.074e-10} in (b).}
  \label{fig:param1-special-1}
\end{figure}

To summarize, the parameters that were selected above for the rod
in this example was, with error tolerance \num{e-10},
\begin{equation}
  \begin{array}{r@{\:}c@{\:}l@{\qquad}r@{\:}c@{\:}l@{\qquad}r@{\:}c@{\:}l}
    r_\text{QBX}/h &=& 1, & N_\text{U} &=& 3, \\
    d_\text{QBX} &=& 0.169, & d_\text{U1} &=& 1.061, & \kappa_\text{U1} &=& 2, \\
    p_\text{QBX} &=& 40, & d_\text{U2} &=& 0.391, & \kappa_\text{U2} &=& 3, \\
    \kappa_\text{QBX} &=& 20, & d_\text{U3} &=& 0.237, & \kappa_\text{U3} &=& 4. \\
  \end{array}
\end{equation}
If the selected QBX upsampling factor $\kappa_\text{QBX}$ seems
large, recall that this parameter is completely hidden in the QBX
precomputation step, as explained in
section~\ref{sec:qbx-precomp}. To verify that the selected
parameters keep the error below the tolerance, we plot in
Figure~\ref{fig:param1-special-distance}~(a) the maximum error
along the 45 lines that were used earlier (shown in
Figure~\ref{fig:param1-direct-line}~(a)). We also plot the error
in two slices in Figure~\ref{fig:param1-special-1}. The fact that
the error slightly exceeds the tolerance at some points in (b)
should come as no surprise, since we have used the error only
along certain lines to select the parameters, not in the whole
space. All the points where the tolerance is exceeded are close
to the boundary between different quadrature regions and could
thus be eliminated by adjusting the threshold distances slightly
upwards (which we will however not do here).

\begin{figure}[h!]
  \centering\small
  \begin{minipage}{0.38\textwidth}%
    \centering
    \includegraphics[scale=0.35]{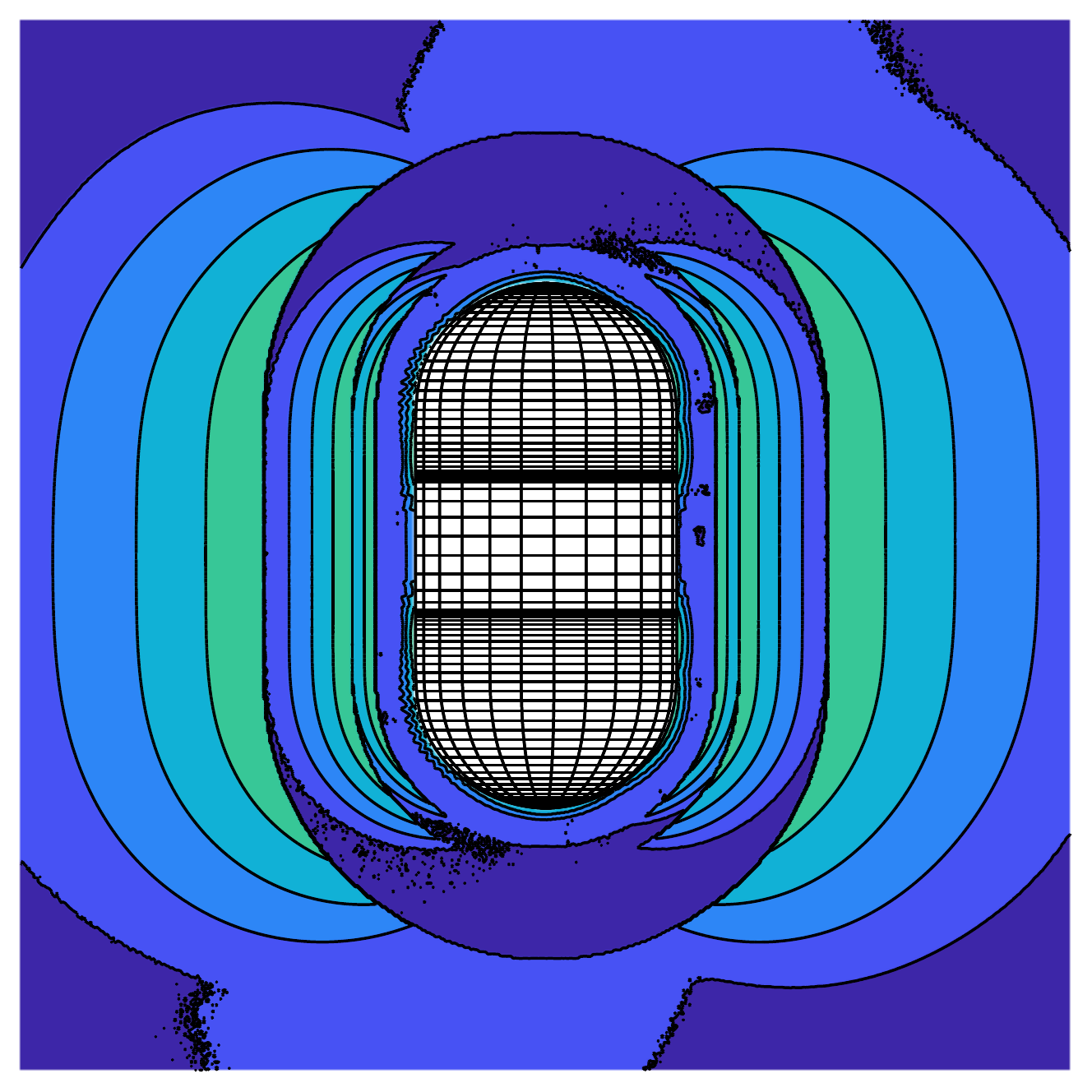}\\%
    (a) Special quadrature error,\\tolerance \num{e-6}, slice~1
  \end{minipage}%
  \begin{minipage}{0.38\textwidth}%
    \centering
    \includegraphics[scale=0.35]{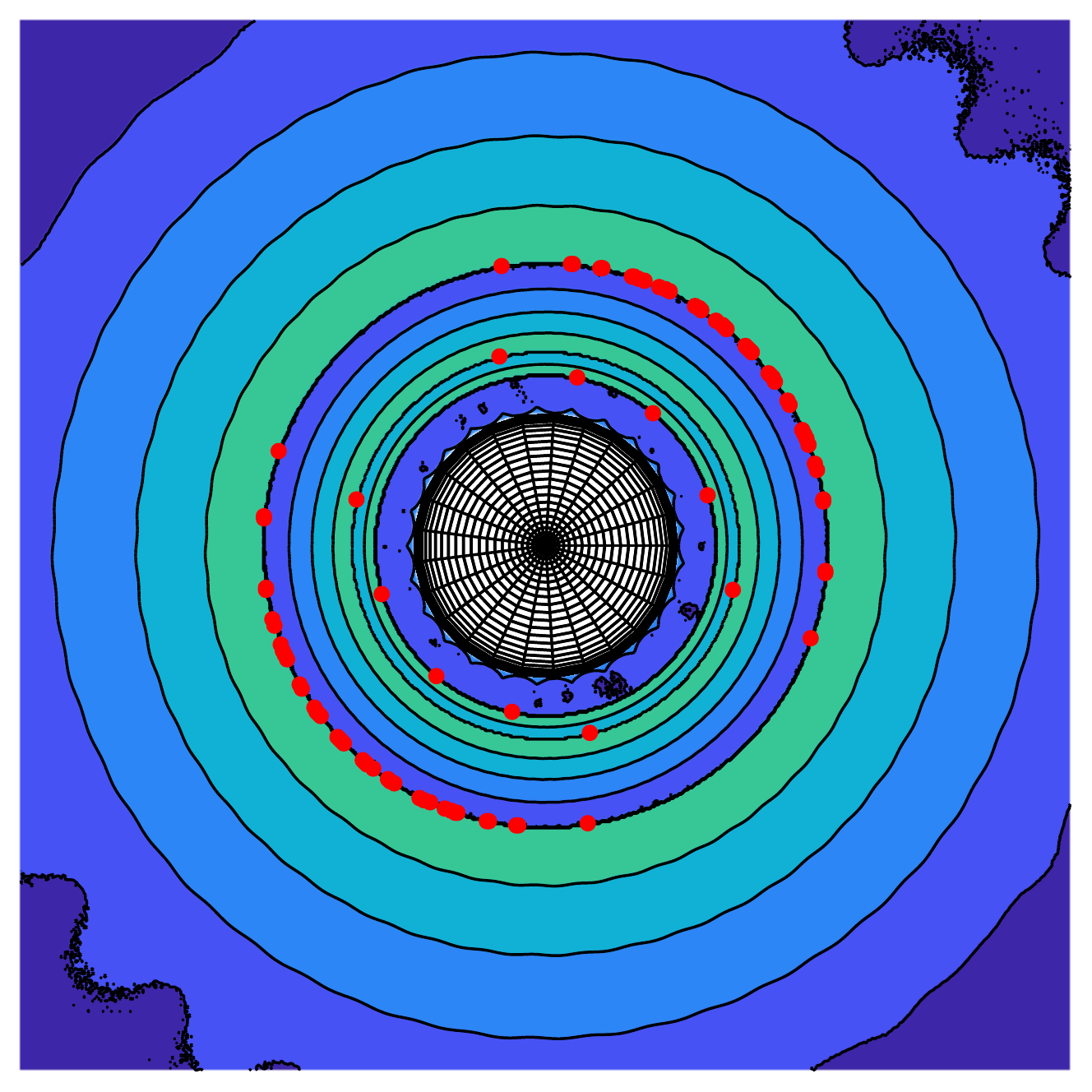}\\%
    (b) Special quadrature error,\\tolerance \num{e-6}, slice~2
  \end{minipage}%
  \begin{minipage}{0.1\textwidth}%
    \centering
    \includegraphics[scale=1]{fig/colorbar.pdf}\\%
    \vphantom{(c)}
  \end{minipage}%
  \caption{Error in two perpendicular slices (a) and (b) through
  the rod particle, when evaluating the stresslet
  identity~\eqref{eq:stresslet-identity} using combined special quadrature
  with tolerance \num{e-6} in free space. Each slice consists of
  $500 \times 500$ evaluation points. All points in (a) are below
  the tolerance, but 90 points in (b), marked red, are above the tolerance.
  The largest error is \num{9.214e-7} in (a) and \num{1.079e-6} in (b).}
  \label{fig:param1-special-2}
\end{figure}

The parameter selection procedure is repeated for the same
rod with the looser error tolerance \num{e-6}. The parameters for
tolerance \num{e-6} are
\begin{equation}
  \label{eq:param1-params-tol6}
  \begin{array}{r@{\:}c@{\:}l@{\qquad}r@{\:}c@{\:}l@{\qquad}r@{\:}c@{\:}l}
    r_\text{QBX}/h &=& 1, & N_\text{U} &=& 2, \\
    d_\text{QBX} &=& 0.149, & d_\text{U1} &=& 0.575, & \kappa_\text{U1} &=& 2, \\
    p_\text{QBX} &=& 21, & d_\text{U2} &=& 0.238, & \kappa_\text{U2} &=& 3, \\
    \kappa_\text{QBX} &=& 15, \\
  \end{array}
\end{equation}
and the error is shown in Figure~\ref{fig:param1-special-distance}~(b)
and Figure~\ref{fig:param1-special-2}.

\FloatBarrier
\subsection{A note on the computational cost}
\label{sec:param1-timing}

While our parameter selection strategy does not try to optimize
the computational cost, we naturally strive for a reasonably low
cost. We therefore comment on the computational cost for the
different quadrature methods considered here. The computational
complexity for evaluating the layer potential using each
quadrature method is shown in Table~\ref{tab:complexity-evaluation}.
The precomputation time (for constructing the interpolation matrices
and QBX matrices), which is naturally independent of the number of
evaluation points, is not included. Note that the total
evaluation time depends on the number of evaluation points in
each quadrature region, and also on the number of expansions that
are used for the evaluation points in the QBX region (recall that
the closest expansion centre is used for each evaluation point).

An example of evaluation times for a specific computer machine is
given in Table~\ref{tab:example-time}, again excluding
precomputation. The time required to find the closest expansion
centre for each evaluation point in the QBX region has been
omitted from Tables~\ref{tab:complexity-evaluation} and
\ref{tab:example-time} since it is negligible (around $\num{e-8}
\times N_\text{eval,QBX}$ seconds).

For a particle, $N_\text{grid}$ is the number of
grid points on the whole particle, i.e.\ $N_\text{grid} = 2250$
for the rod that we have considered so far. Let us study the
special case of a single evaluation point, relevant for example
when computing a streamline. Based on
Table~\ref{tab:example-time}, the evaluation time for this single
point can be computed, depending on which quadrature region the
point belongs to and the parameter $\kappa_\text{U$i$}$ or
$p_\text{QBX}$. This is shown in Table~\ref{tab:example-time-single}.
From this, it can for example be seen that the evaluation takes
roughly 1000 times longer for a point in the upsampled quadrature
region with $\kappa_\text{U$i$}=2$ compared to the direct
quadrature region. (The upsampled quadrature cost is in this case
completely dominated by interpolating the density, i.e.\
multiplying it by the precomputed interpolation matrix.)

\begin{table}[h!]
  \centering
  \caption{Computational complexities for the different
  quadrature methods.}
  \label{tab:complexity-evaluation}
  \begin{tabular}{p{4.5cm}@{}p{7cm}}
    \toprule
    Direct quadrature \\
    \quad Evaluate & $T_\text{DE} = O(N_\text{grid} N_\text{eval,D})$ \\
    \midrule
    Upsampled quadrature \\
    \quad Interpolate density & $T_\text{UI,$i$} = O(\kappa_\text{U$i$}^2 N_\text{grid}^2)$ \\
    \quad Evaluate & $T_\text{UE,$i$} = O(\kappa_\text{U$i$}^2 N_\text{grid}
    N_\text{eval,U$i$})$ \\
    \midrule
    QBX \\
    \quad Compute coefficients & $T_\text{QC} = O(p_\text{QBX}^2 N_\text{grid} N_\text{exp})$ \\
    \quad Evaluate expansion & $T_\text{QE} = O(p_\text{QBX}^2 N_\text{eval,QBX})$ \\
    \bottomrule
  \end{tabular}\\[5pt]
  \begin{minipage}{11.8cm}\small
    Time complexities for evaluating the double layer
    potential~$\vec{\Dp}$, excluding precomputation time.
    Here,
    \begin{itemize}
      \item
        $N_\text{grid}$ is the total number of grid points on the
        part of the surface included in the special quadrature
        method (i.e.\ $\tilde{\Gamma}$ as defined in
        section~\ref{sec:global-local-qbx}),
      \item
        $N_\text{eval,D}$, $N_\text{eval,U$i$}$ and $N_\text{eval,QBX}$
        are the number of evaluation points in the direct
        quadrature region, the $i$th upsampled quadrature region
        and the QBX region, respectively,
      \item
        $N_\text{exp}$ is the number of expansion centres that
        are to be used for the evaluation points in the QBX region.
    \end{itemize}
  \end{minipage}
\end{table}%

\begin{table}[h!]
  \centering
  \caption{Example of actual evaluation times [seconds].}
  \label{tab:example-time}
  \begin{tabular}{p{4.5cm}@{}l}
    \toprule
    Direct quadrature \\
    \quad Evaluate & $T_\text{DE} = \num{5.6e-9} \times N_\text{grid} N_\text{eval,D}$ \\
    \midrule
    Upsampled quadrature \\
    \quad Interpolate density & $T_\text{UI,$i$} = \num{8.0e-10}
    \times \kappa_\text{U$i$}^2 N_\text{grid}^2$ \\[2pt]
    \quad Evaluate & $T_\text{UE,$i$} = \num{5.7e-9} \times
    \kappa_\text{U$i$}^2 N_\text{grid} N_\text{eval,U$i$}$ \\
    \midrule
    QBX \\
    \quad Compute coefficients & $T_\text{QC} = \num{5.3e-8} \times
    (0.071 p_\text{QBX}^2 + 0.56 p_\text{QBX} + 1) N_\text{grid} N_\text{exp}$ \\[2pt]
    \quad Evaluate expansion & $T_\text{QE} = \num{2.2e-5} \times
    (0.0053p_\text{QBX}^2 - 0.0027p_\text{QBX} + 1) N_\text{eval,QBX}$ \\
    \bottomrule
  \end{tabular}\\[5pt]
  \begin{minipage}{0.8\textwidth}\small
    These times are for a modern workstation with a 6-core Intel Core i7-8700 CPU
    (4.6~GHz).
  \end{minipage}
\end{table}%

\begin{table}[h!]
  \centering
  \caption{Evaluation times for $N_\text{grid} = 2250$ and
  a single evaluation point (seconds).}
  \label{tab:example-time-single}
  \begin{tabular}{c}
    {\small\bfseries Direct quadrature} \\
    \toprule
    Time [s] \\
    \midrule
    \num{1.3e-05} \\
    \bottomrule
    \vphantom{\num{1.3e-05}}\hspace*{0.001pt} \\
    \vphantom{\num{1.3e-05}}\hspace*{0.001pt} \\
    \vphantom{\num{1.3e-05}}\hspace*{0.001pt} \\
    \vphantom{\num{1.3e-05}}\hspace*{0.001pt} \\
  \end{tabular}\qquad%
  \begin{tabular}{cc}
    \multicolumn{2}{c}{\small\bfseries Upsampled quadrature} \\
    \toprule
    $\kappa_\text{U$i$}$ & Time [s] \\
    \midrule
    2 & \num{1.6e-02} \\
    3 & \num{3.7e-02} \\
    4 & \num{6.5e-02} \\
    5 & \num{1.0e-01} \\
    6 & \num{1.5e-01} \\
    \bottomrule
  \end{tabular}\qquad%
  \begin{tabular}{cc}
    \multicolumn{2}{c}{\small\bfseries QBX (with $N_\mathrm{exp}=1$)} \\
    \toprule
    $p_\text{QBX}$ & Time [s] \\
    \midrule
    10 & \num{1.7e-03} \\
    20 & \num{4.9e-03} \\
    30 & \num{9.9e-03} \\
    40 & \num{1.7e-02} \\
    50 & \num{2.5e-02} \\
    \bottomrule
  \end{tabular}
\end{table}%

It can also be seen in Table~\ref{tab:example-time-single} that
QBX is often faster than upsampled quadrature. For instance, QBX
with $p_\text{QBX} = 40$ takes about as much time as upsampled
quadrature with $\kappa_\text{U$i$}=2$ and is faster than any
$\kappa_\text{U$i$} \geq 3$. However, note that this conclusion
may not hold when there are more than one evaluation point, since
the evaluation time depends in an intricate way on both the
number of evaluation points in each region and the number of
expansions needed for QBX. In particular, if many expansions are
needed (large $N_\text{exp}$) and there are few evaluation points
per expansion, QBX will tend to be slower than upsampled
quadrature due to the large cost of computing coefficients.

\pagebreak
\section{Summary of the parameter selection strategy}
\label{sec:param-general-summary}

The parameter selection strategy can, in the general case, be
summarized as follows. In all steps, the stresslet
identity~\eqref{eq:stresslet-identity} is used to estimate the
error.

\begin{framedbox}
  \textit{Input:} Discretization of the geometry, error
  tolerance~$\etol$\\[3pt]
  \textit{Output:} Parameters $N_\text{U}$, $(d_\text{U$i$},
  \kappa_\text{U$i$})_{i=1}^{N_\text{U}}$, $d_\text{QBX}$, $r_\text{QBX}$,
  $p_\text{QBX}$, $\kappa_\text{QBX}$ for the special
  quadrature

  \vspace{2\itemsep}
  \noindent
  For each distinct geometrical object:
  \begin{enumerate}
    \item
      Put $\kappa_\text{U$i$} = i+1$. Numerically determine the
      threshold distances $d_\text{U$i$}$ to keep the error below
      $\etol$, as in Figure~\ref{fig:param1-upsamp-line} and
      Table~\ref{tab:param1-upsamp}, up to the first $i$ such
      that $d_\text{U$i$} \leq 2 \gamma h$, where $\gamma = 0.85$
      and $h$ is the grid spacing (defined for particles in
      equation \eqref{eq:grid-spacing-particles} and for walls in
      section~\ref{sec:param-example3}).
    \item
      Put $r_\text{QBX} = h$. If the last (smallest) $d_\text{U$i$}$
      computed in step~1 lies in the interval $[h, 2\gamma h]$,
      put $d_\text{QBX}$ equal to it. Otherwise, put
      $d_\text{QBX} = h$. This also sets $N_\text{U}$, the number
      of upsampled quadrature regions.
    \item
      Choose $p_\text{QBX}$ such that the truncation error is
      below $\etol$, based on a parameter study such as in
      Figure~\ref{fig:param1-qbx-p-kr}~(b). Choose
      $\kappa_\text{QBX}$ such that the coefficient error is no
      larger than the truncation error.
  \end{enumerate}
\end{framedbox}

\noindent
This strategy is designed to keep the error relative to
$\max_\Gamma \absi{\vec{q}}$ below $\etol$ when evaluating the
layer potential, provided that the density~$\vec{q}$ is well-resolved
by the discretization. While there is no guarantee that
the error stays strictly below $\etol$, empirical evidence in
sections~\ref{sec:param1-results}, \ref{sec:param-example2},
\ref{sec:param-example3}, \ref{sec:res1-special-quadrature} and
\ref{sec:res2-convergence} indicates that the error is typically
close to the tolerance, and in any case of the correct order of
magnitude. Note that the procedure, including the parameter
studies, must be repeated every time a new geometrical object,
such as a rod particle with a different aspect ratio, is used. We
will now apply the procedure to two additional examples.

\section{Example II: a rod particle with higher aspect ratio}
\label{sec:param-example2}

\begin{figure}[b!]
  \centering\small
  \begin{minipage}{0.38\textwidth}%
    \centering
    \includegraphics[scale=0.4]{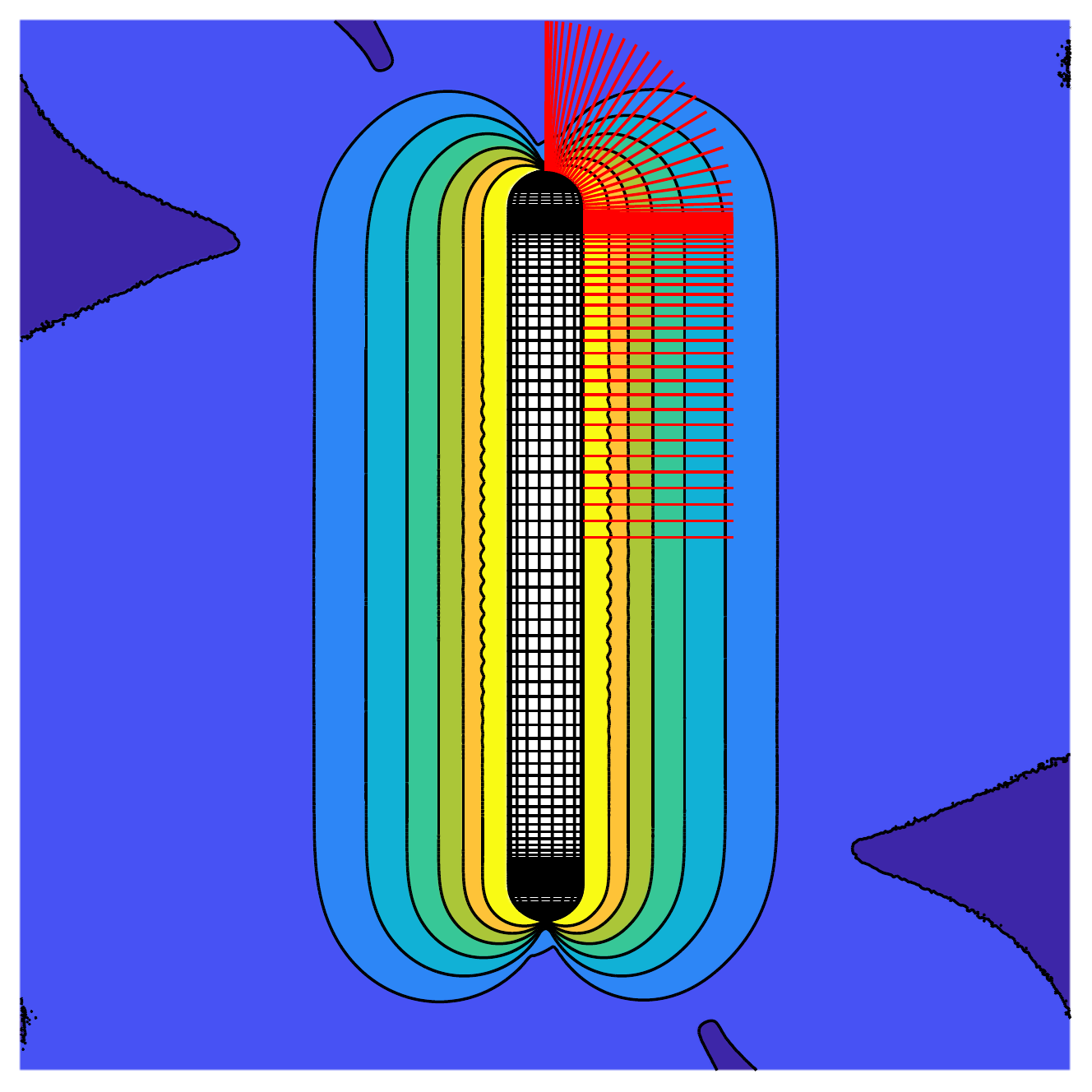}\\%
    (a) Direct quadrature error\\\vphantom{tolerance}
  \end{minipage}%
  \begin{minipage}{0.38\textwidth}%
    \centering
    \includegraphics[scale=0.4]{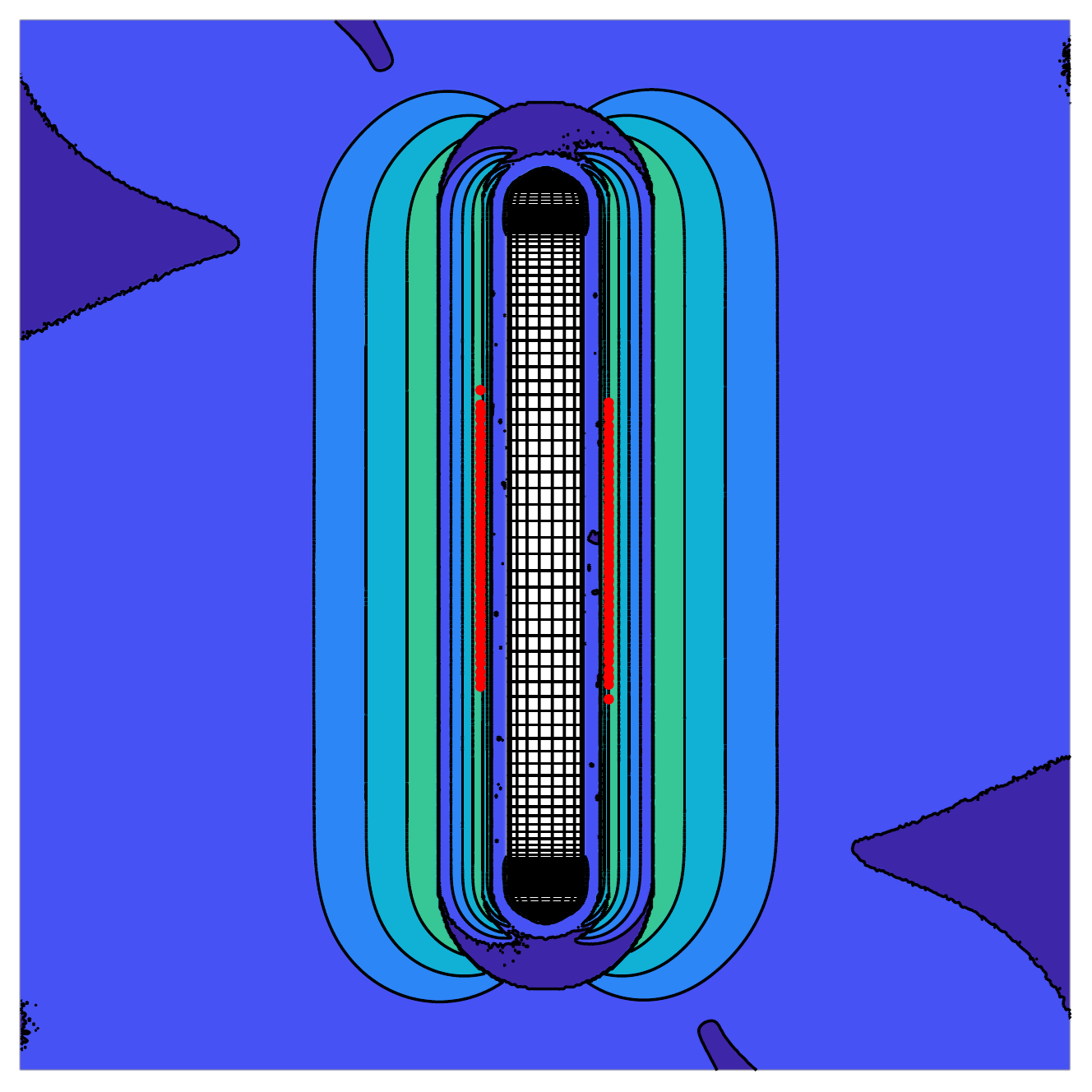}\\%
    (b) Special quadrature error,\\tolerance \num{e-6}
  \end{minipage}%
  \begin{minipage}{0.1\textwidth}%
    \centering
    \includegraphics[scale=1]{fig/colorbar.pdf}\\%
    \vphantom{(c)}\hspace*{0.001pt}\\\vphantom{tolerance}
  \end{minipage}%
  \caption{Error in a slice through the rod particle, when
  evaluating the stresslet identity~\eqref{eq:stresslet-identity}
  in free space using (a) direct quadrature, and (b) combined
  special quadrature with tolerance \num{e-6}. In (b), the
  tolerance is exceeded in 148~points, marked red; the evaluation
  grid consists of $500 \times 500$ points and the largest error
  is \num{1.455e-6}. In (a), the 65~lines used to select parameters
  are shown in red.}
  \label{fig:param2-slice}
\end{figure}

For the second example, we consider a more slender rod particle,
namely the rod of length~$L=10$ and radius~$R=0.5$ (aspect ratio~10)
shown in Figure~\ref{fig:param2-slice}. The grid has parameters
$n_1 = 35$, $n_2 = 60$ and $n_\varphi = 18$, in total 2340 grid points.
The grid spacing as defined by \eqref{eq:grid-spacing-particles}
is $h = \pi/18 \approx 0.1745$. It should be noted that while $h$
is based only on the grid resolution in the azimuthal direction,
the resolution in the polar direction must not be much coarser.
Otherwise, the distance between QBX expansion centres would be too
large in the polar direction, and the balls of convergence would
not cover the QBX region. (This limitation is due to having one
expansion centre per grid point.)

Applying the same constant density $\tilde{\vec{q}} = (1,1,1) /
\sqrt{3}$, the direct quadrature error is shown in
Figure~\ref{fig:param2-slice}~(a). The special quadrature
parameters are selected as described in
section~\ref{sec:param-general-summary}, with the error along the
65~lines shown in red in Figure~\ref{fig:param2-slice}~(a) used
to select the threshold distances. The parameters for error
tolerance~\num{e-6} are
\begin{equation}
  \label{eq:param2-params-tol6}
  \begin{array}{r@{\:}c@{\:}l@{\qquad}r@{\:}c@{\:}l@{\qquad}r@{\:}c@{\:}l}
    r_\text{QBX}/h &=& 1, & N_\text{U} &=& 2, \\
    d_\text{QBX} &=& 0.234, & d_\text{U1} &=& 0.930, & \kappa_\text{U1} &=& 2, \\
    p_\text{QBX} &=& 28, & d_\text{U2} &=& 0.355, & \kappa_\text{U2} &=& 3. \\
    \kappa_\text{QBX} &=& 20, \\
  \end{array}
\end{equation}
The error when using these parameters is shown in
Figure~\ref{fig:param2-slice}~(b). As before, the tolerance is
not strictly enforced, but the error stays within a factor~2 of
the tolerance.

The parameters \eqref{eq:param2-params-tol6} for the slender rod
can be compared with the parameters \eqref{eq:param1-params-tol6}
for the less slender rod with
the same tolerance; the threshold distances are relative to the
diameter of the cylindrical part of the rod in both cases. Note
that the slender rod \eqref{eq:param2-params-tol6} has larger
threshold distances than the other rod
\eqref{eq:param1-params-tol6}. This reflects the fact that the
error of the underlying direct quadrature at a fixed distance
from the rod is higher for the slender rod, since it has
lower overall resolution (grid points per surface area). The
slender rod also requires a higher $p_\text{QBX}$ since
$r_\text{QBX} = h$ is larger for the slender rod.

\section{Example III: a pair of plane walls}
\label{sec:param-example3}

\begin{figure}[h!]
  \centering\small
  \begin{minipage}[b]{0.43\textwidth}%
    \centering
    \includegraphics[scale=0.5,trim=0mm 22mm 0mm 22mm,clip]%
    {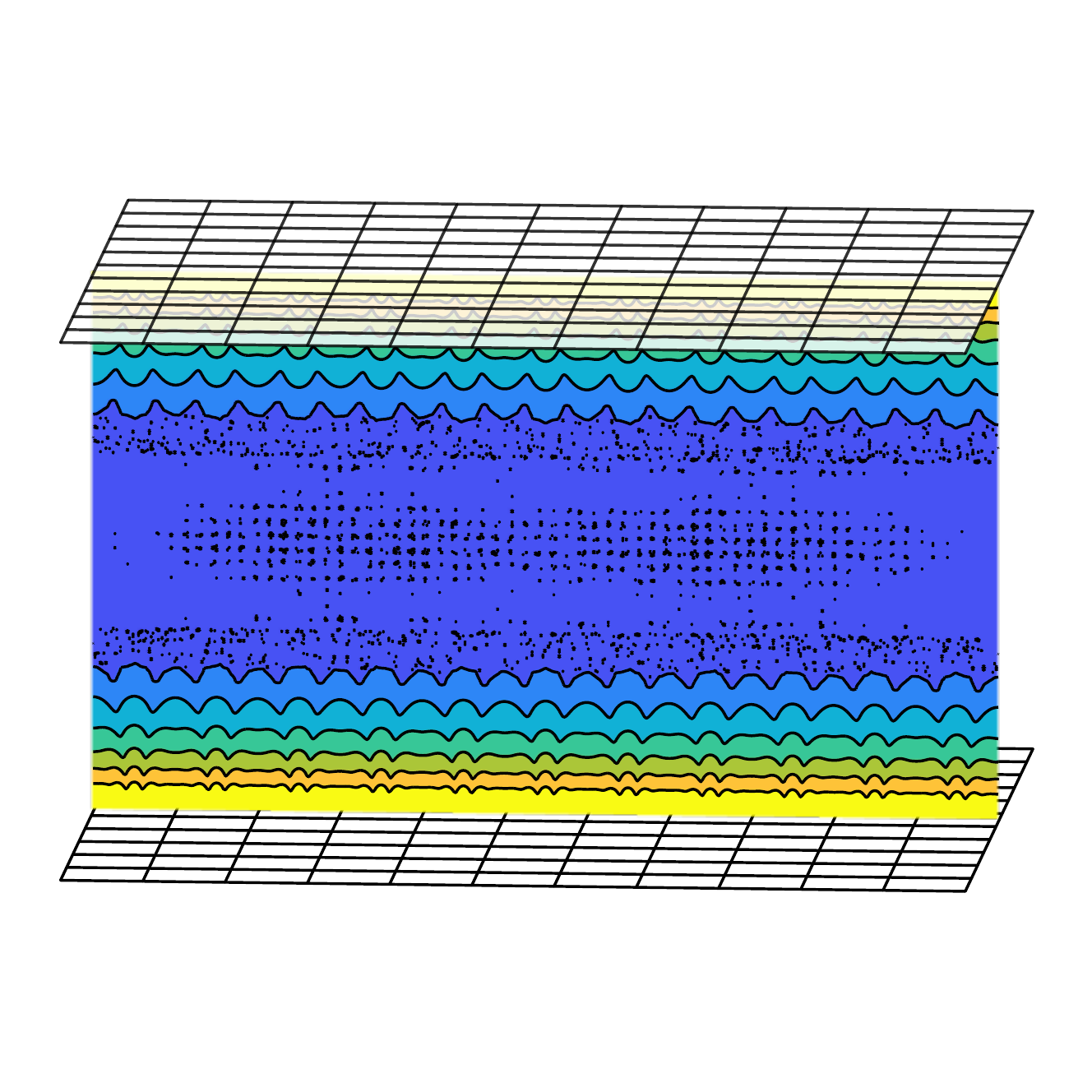}\\[5mm]%
    (a) Direct quadrature error
  \end{minipage}%
  \begin{minipage}[b]{0.12\textwidth}%
    \includegraphics[scale=1]{fig/colorbar.pdf}\\[8mm]%
    \vphantom{(c)}
  \end{minipage}%
  \begin{minipage}[b]{0.45\textwidth}%
    \centering
    \includegraphics[scale=0.75]{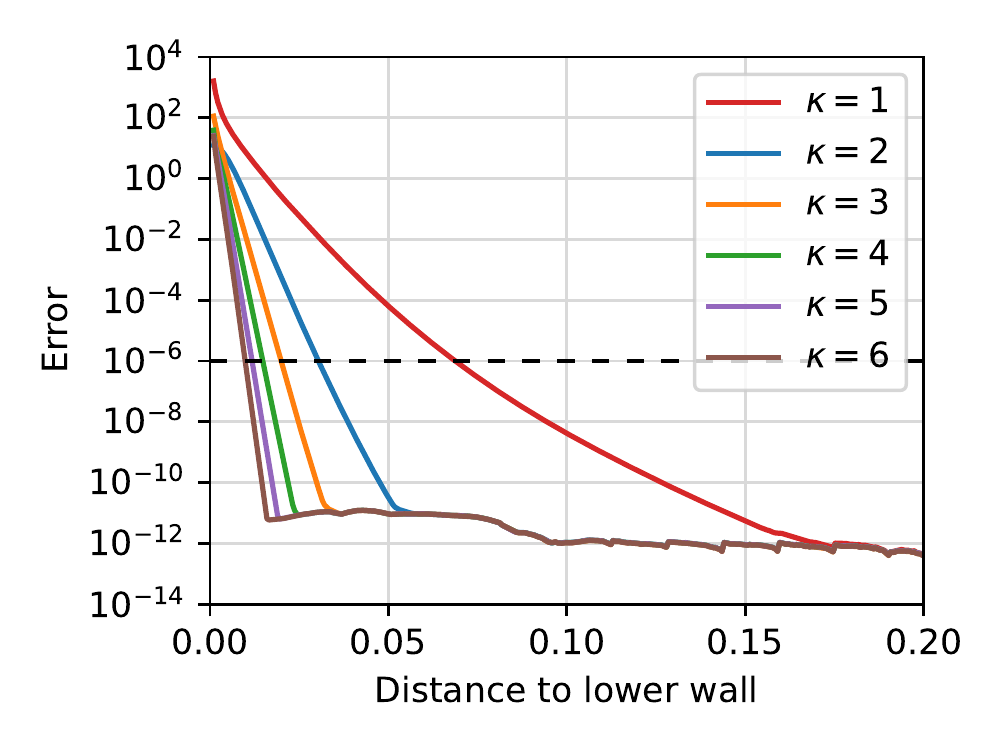}\\[-5pt]%
    \hspace{12mm}(b) Upsampled quadrature error
  \end{minipage}%
  \caption{Errors when evaluating the stresslet
  identity~\eqref{eq:stresslet-identity} in the periodic setting.
  (a)~Error in the centre plane for direct quadrature.
  (b)~Largest error for upsampled quadrature with different
  upsampling factors~$\kappa$, as a function of the distance to
  the lower wall.}
  \label{fig:param3-direct-upsamp}
\end{figure}

In this third example, we select parameters for a plane wall.
Since we always consider walls in a periodic setting, we will do
so here as well, and use the Spectral Ewald (SE) method described
in section~\ref{sec:periodicity}. We will here select the SE
parameters such that the error from SE is completely negligible
compared to the quadrature errors which we strive to control
here.%
\footnote{
  Specifically, the SE parameters used here are $\xi = 15.245$,
  $r_\text{c} = 0.4$, $P = 24$, and the uniform grid used for the
  Fourier-space part has $64 \times 64 \times 64$ grid points (see
  section~\ref{sec:periodicity} and
  appendix~\ref{app:streamlines} for an explanation). These parameters
  should keep the SE error around \num{e-15} according to
  \cite{klinteberg14}.
}
Since the problem is periodic in all three spatial directions, we
must have a pair of walls so that the fluid domain can be
confined between them. The periodic cell is here of size $\vec{B}
= (1,1,1)$ and the two walls are placed at a distance of 0.6 from
each other. The walls are discretized using $11 \times 11$ patches
each, with $8 \times 8$ grid points on each patch (as described
in section~\ref{sec:direct-quad-walls}), in total 7744 grid points
per wall. For a wall, we define
the grid spacing as $h = \max(h_1, h_2)$, where $h_1$ and
$h_2$ are the largest spacings between grid points in each of the
two tensorial directions of the patches. The walls considered here
have grid spacing $h \approx 0.01668$.
The constant density $\tilde{\vec{q}} = (0, 0, 1)$ is applied in
the direction of the normal of the lower wall (pointing into the
fluid domain). The direct quadrature error is shown in
Figure~\ref{fig:param3-direct-upsamp}~(a).

We follow the procedure in section~\ref{sec:param-general-summary}.
The threshold distances of the upsampled quadrature regions are
computed by evaluating the stresslet identity error along normal
lines of the walls, with each line centred at a grid point. The
error is plotted in Figure~\ref{fig:param3-direct-upsamp}~(b),
and the resulting threshold distances for error tolerance
\num{e-6} are $d_\text{U1} = 0.0687$, $d_\text{U2} = 0.0304$ and
$d_\text{QBX} = 0.0198$.

\begin{figure}[h!]
  \centering\small
  \begin{minipage}[b]{0.5\textwidth}%
    \centering
    \includegraphics[scale=0.76]{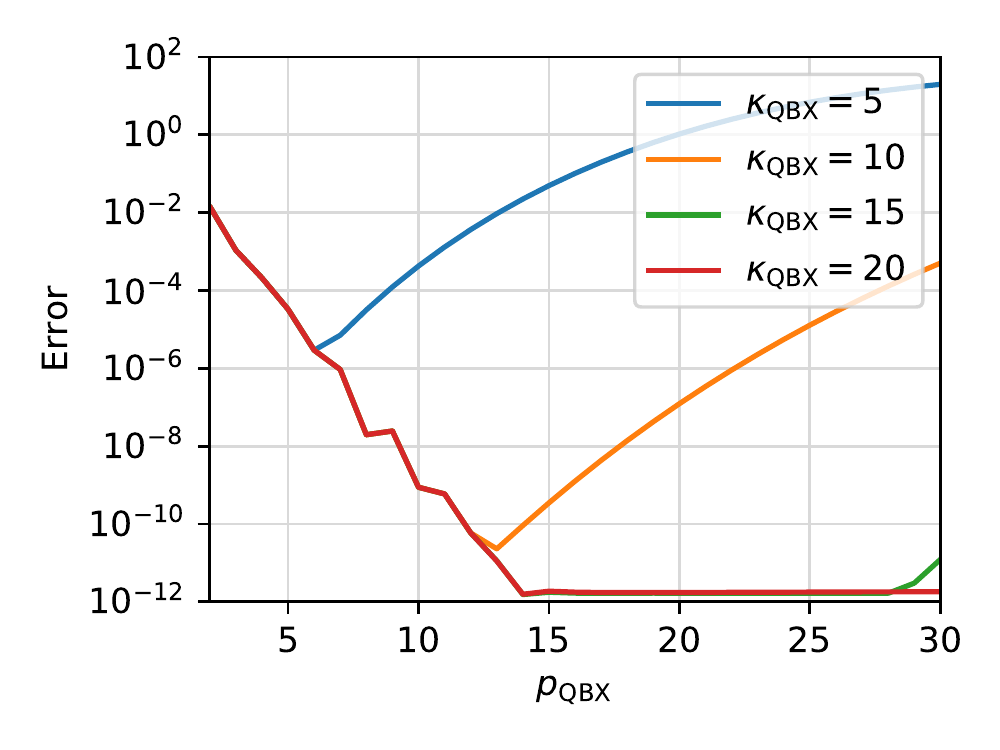}\\[-3mm]%
    \hspace*{6mm}(a) Onsurface QBX error for different $\kappa_\text{QBX}$
  \end{minipage}%
  \begin{minipage}[b]{0.5\textwidth}%
    \centering
    \includegraphics[scale=0.95]{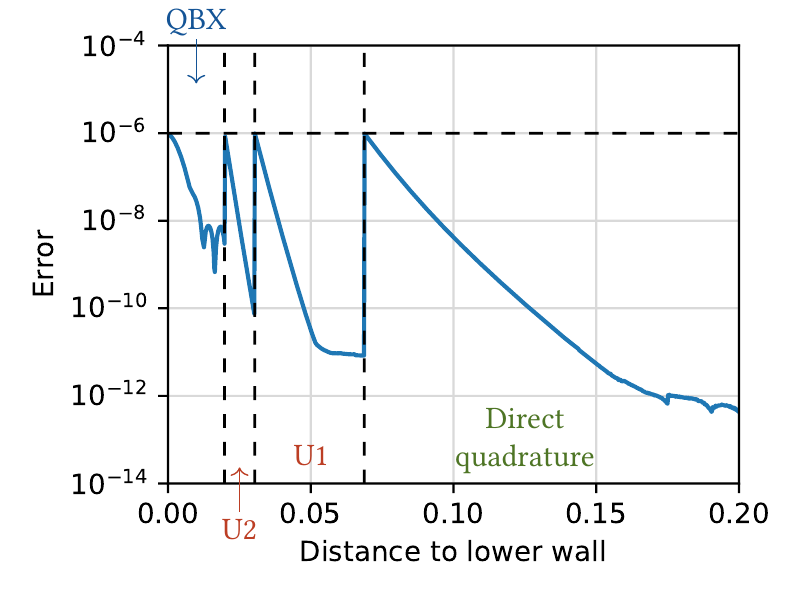}\\[-3mm]%
    \hspace*{6mm}(b) Special quadrature error, tolerance \num{e-6}
  \end{minipage}%
  \caption{(a)~Maximal onsurface error when evaluating the
  stresslet identity using QBX, for different $\kappa_\text{QBX}$
  with $r_\text{QBX}/h = 1$ fixed.
  (b)~Maximal stresslet identity error as a function
  of the distance to the lower wall (for 1000 equispaced
  distances in $[0, 0.2]$), using the combined special quadrature with
  tolerance~\num{e-6}. The largest error is \num{9.812e-7}.}
  \label{fig:param3-special-distance}
\end{figure}

To determine $p_\text{QBX}$ and $\kappa_\text{QBX}$, we do a
parameter study, shown in Figure~\ref{fig:param3-special-distance}~(a).
Note that the plane wall needs a significantly lower
$p_\text{QBX}$ than the rod particles to reach a given error. The
error curves in Figure~\ref{fig:param3-special-distance}~(a)
level out at around \num{e-12} due to other errors not controlled
by the QBX parameters. In order to reach the tolerance \num{e-6}
it is\pagebreak\ sufficient to choose $p_\text{QBX} = 7$ and
$\kappa_\text{QBX} = 10$. The selected parameters are thus
\begin{equation}
  \label{eq:param3-params-tol6}
  \begin{array}{r@{\:}c@{\:}l@{\qquad}r@{\:}c@{\:}l@{\qquad}r@{\:}c@{\:}l}
    r_\text{QBX}/h &=& 1, & N_\text{U} &=& 2, \\
    d_\text{QBX} &=& 0.0198, & d_\text{U1} &=& 0.0687, & \kappa_\text{U1} &=& 2, \\
    p_\text{QBX} &=& 7, & d_\text{U2} &=& 0.0304, & \kappa_\text{U2} &=& 3. \\
    \kappa_\text{QBX} &=& 10, \\
  \end{array}
\end{equation}
The error when using these parameters is shown in
Figure~\ref{fig:param3-special-distance}~(b).

\chapter{Numerical results}
\label{sec:results}

Our numerical method can be summarized as follows:
\begin{enumerate}
  \item
    The geometry is discretized as in
    section~\ref{sec:disc-quadrature}. Parameters for the
    combined special quadrature method are selected as in
    section~\ref{sec:param-general-summary}. Parameter selection
    is done in free space for particles, but the same parameters
    can also be used in the periodic setting. For walls,
    parameter selection is done in the periodic setting.

  \item
    The matrices $\mat{M}_i$ and $\mat{R}_i$ for offsurface QBX
    and onsurface QBX, respectively, are precomputed as in
    section~\ref{sec:qbx-precomp}. Interpolation matrices for the
    upsampled quadrature regions are also precomputed. At this
    point, the special quadrature is ready to be used to evaluate
    the layer potential.

  \item
    The boundary integral equation for either a resistance
    problem or a mobility problem is solved iteratively for
    $\vec{q}$ using GMRES. The Spectral Ewald method is used for
    periodic problems, as described in section~\ref{sec:periodicity}.
    A preconditioner is used in all cases, as described below.
    \begin{itemize}
      \item
        For a resistance problem, velocities are given for all
        particles and the boundary integral equation is given by
        equation~\eqref{eq:resistance-disc}.
      \item
        For a mobility problem, forces and torques are applied to
        all particles and the boundary integral equation is given
        by equation~\eqref{eq:mobility-disc}.
    \end{itemize}

  \item
    The flow field in the fluid domain may be computed in a
    postprocessing step using $\vec{q}$ from step~3. For a
    resistance problem, the forces and torques acting on all
    particles may also be computed here, while for a mobility
    problem, the particle velocities may be computed.
\end{enumerate}
To improve the convergence of GMRES in step~3, we use a
block-diagonal preconditioner similar to the one used in
\cite{klinteberg16b}. The preconditioner is constructed by
computing the explicit inverse of a single-particle system as
well as a system consisting of a single wall patch (if walls are
present in the simulation). These two types of blocks are then
placed along the diagonal and rotated according to the geometry.
This preconditioner has been seen to reduce the number of GMRES
iterations by as much as a factor~17 for some systems with many
particles, such as the ones in section~\ref{sec:res3-packed-rods}.

In this section, we test some aspects of our numerical method,
with focus on the special quadrature. First, in
section~\ref{sec:res1-special-quadrature}, we test the quadrature
on its own (i.e.\ steps~1--2 above) with geometries containing
both particles and walls. This serves as a continuation of the
tests in section~\ref{sec:parameters}, where geometrical objects
were considered only separately. In section~\ref{sec:res2-convergence},
we test the special quadrature in the context of the full numerical
method (steps~1--4 above) and in particular how the quadrature
tolerance influences the accuracy. Finally, in
section~\ref{sec:res3-packed-rods}, we test the computational
complexity of the method on a more complicated problem, and
compute streamlines.

\section{Special quadrature with composite geometries}
\label{sec:res1-special-quadrature}

We consider two geometrical setups, shown in Figures~\ref{fig:ch71-slice1}
and \ref{fig:ch71-B-direct}. Both problems are periodic with a
periodic cell of size $\vec{B} = (1, 1, 1)$, and the Spectral
Ewald parameters are as in section~\ref{sec:param-example3}.
As in section~\ref{sec:parameters}, we use the stresslet
identity~\eqref{eq:stresslet-identity} to estimate the error.
This is the same test used to select the parameters, so it mainly
serves as a consistency check (tests with nonconstant densities
will follow in section~\ref{sec:res2-convergence}).

\subsection{Geometry 1: Two rods between a pair of plane walls}
\label{sec:res1-special-quadrature-geometry-1}

The first geometry consists of two plane walls discretized as in
section~\ref{sec:param-example3}, at a distance~0.6 from each
other. Between these walls are two rod particles of length~$L =
0.5$ and radius $R = L/20$ (aspect ratio~10), discretized as in
section~\ref{sec:param-example2} (but scaled down a factor~20),
oriented such that their axes lie in the centre plane.

\begin{figure}[H]
  \centering\small
  \begin{minipage}{0.45\textwidth}%
    \centering
    \includegraphics[scale=0.45,trim=0mm 22mm 0mm 22mm,clip]%
    {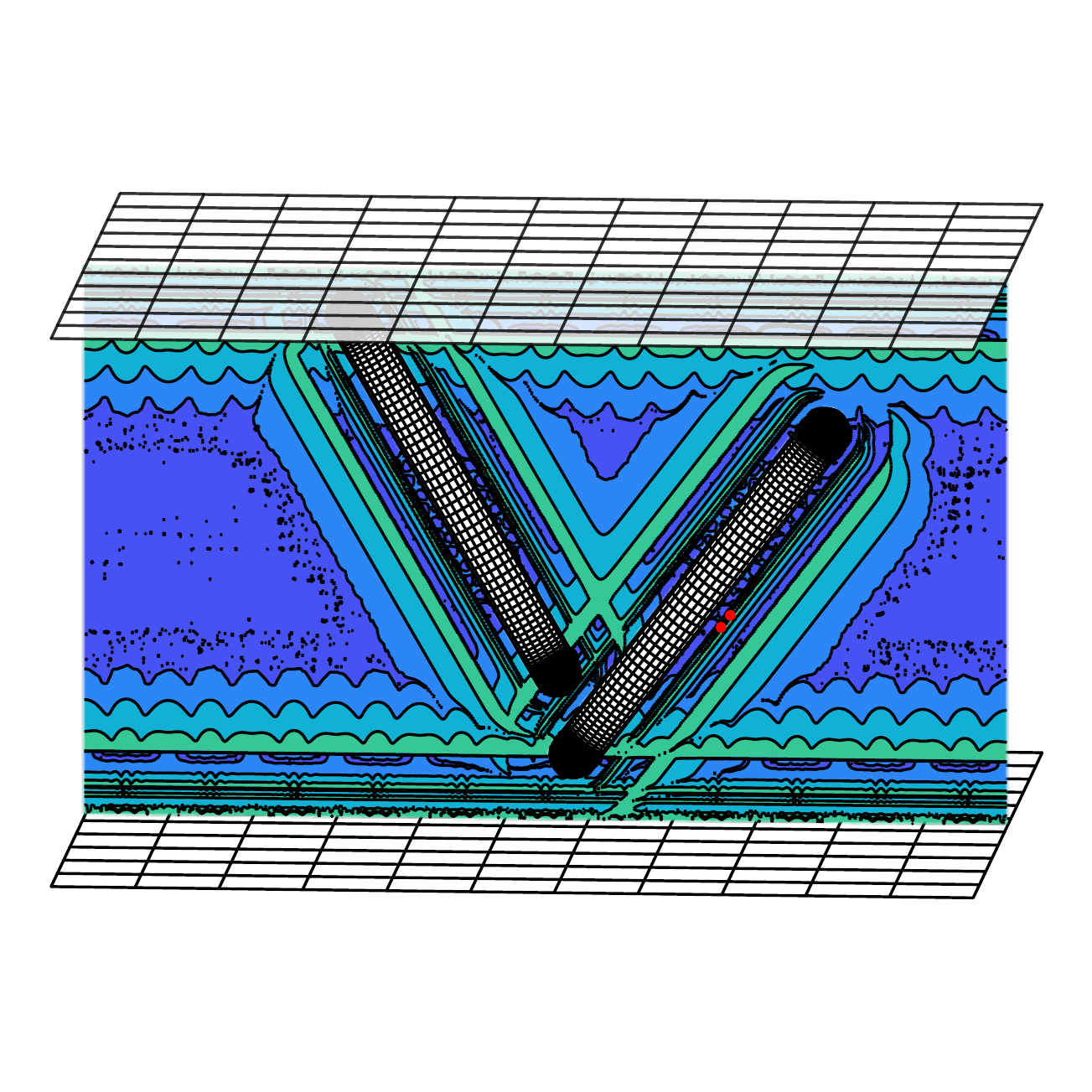}\\%
    (a) Special quadrature error, tolerance \num{e-6}
  \end{minipage}%
  \begin{minipage}{0.45\textwidth}%
    \centering
    \includegraphics[scale=0.45,trim=0mm 22mm 0mm 22mm,clip]%
    {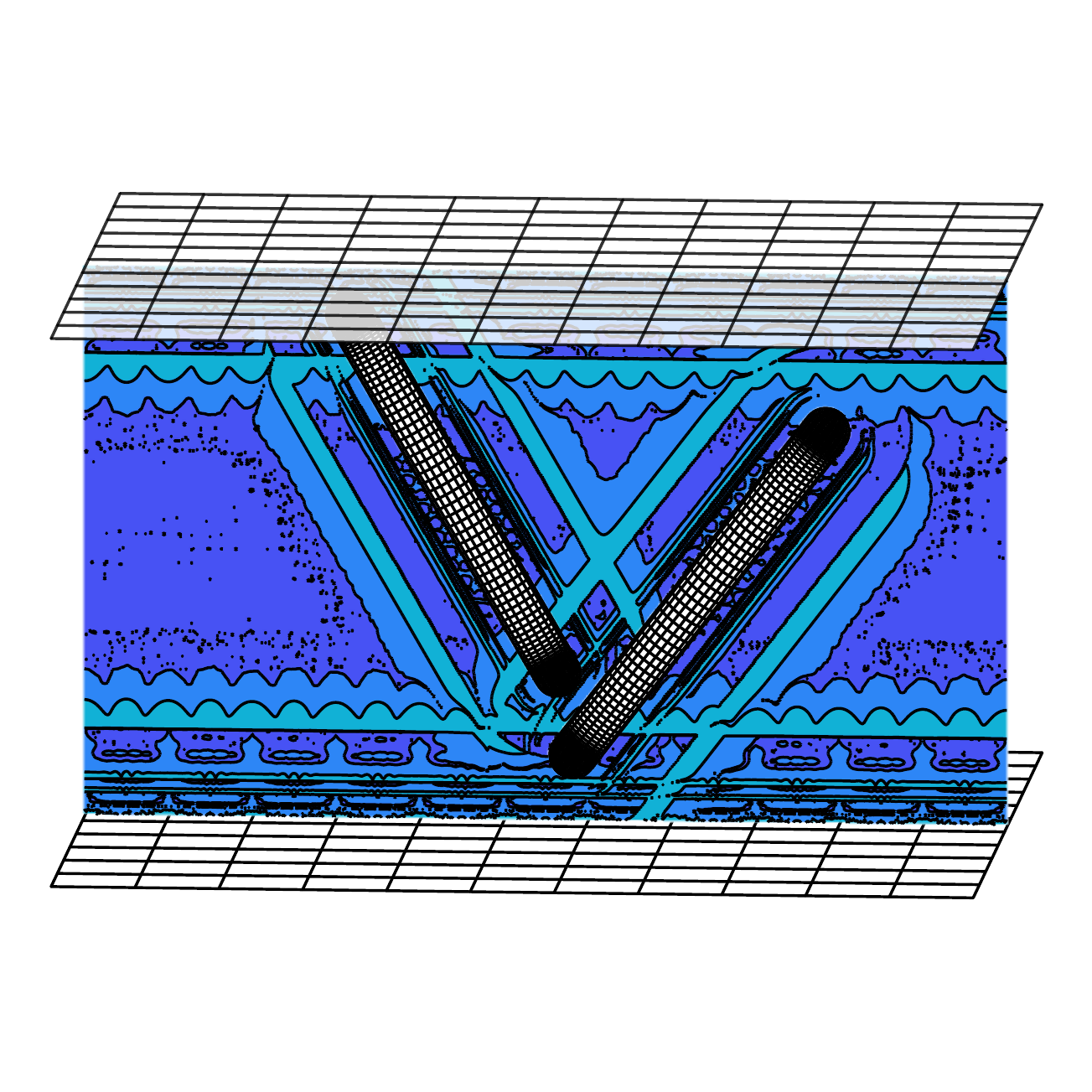}\\%
    (b) Special quadrature error, tolerance \num{e-8}
  \end{minipage}%
  \begin{minipage}{0.1\textwidth}%
    \centering
    \includegraphics[scale=1]{fig/colorbar.pdf}\\
    \vphantom{(c)}
  \end{minipage}%
  \caption{Stresslet identity error in the centre plane, for
  geometry~1, in (a) for tolerance \num{e-6} and in (b) for
  tolerance \num{e-8}. The largest error is \num{1.205e-6}
  in (a) and \num{9.389e-9} in (b). In (a), the error exceeds the
  tolerance in 2 points (the evaluation grid has $500 \times 500$
  points), marked red.}
  \label{fig:ch71-slice1}
\end{figure}

\begin{figure}[H]
  \centering
  \includegraphics[scale=0.75,trim=0mm 5mm 0mm 3.5mm,clip]{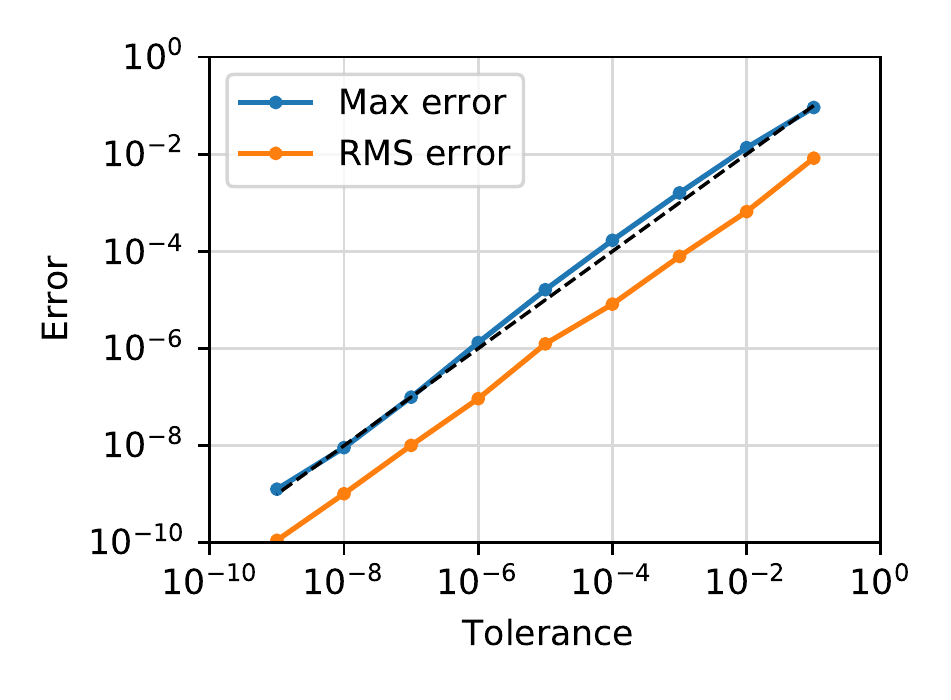}
  \caption{Maximal and root-mean-square (RMS) stresslet identity
  error in the centre plane as a function of the special
  quadrature tolerance, for geometry~1. As observed already in
  section~\ref{sec:parameters}, the tolerance is sometimes
  exceeded slightly at the threshold distances, which causes the
  max error curve to lie above the identity line Error =
  Tolerance (dashed).}
  \label{fig:ch71-SI-conv}
\end{figure}

The stresslet identity error is shown for two different
quadrature tolerances $\etol$ in Figure~\ref{fig:ch71-slice1}.
In (a), $\etol=\num{e-6}$, the quadrature parameters for the
walls are as in section~\ref{sec:param-example3}, i.e.\ given by
\eqref{eq:param3-params-tol6}; for the rods, the parameters are
as in section~\ref{sec:param-example2} but with all distances
scaled by $1/20$ to account for the difference in size. Thus, the
parameters for $\etol = \num{e-6}$ for the rods are $d_\text{QBX}
= 0.0117$, $p_\text{QBX} = 28$, $\kappa_\text{QBX} = 20$ and
$N_\text{U} = 2$, $d_\text{U1} = 0.0465$, $d_\text{U2} = 0.0178$.
In (b), for tolerance $\etol = \num{e-8}$, the parameters
selected according to section~\ref{sec:param-general-summary}
are, for the walls $d_\text{QBX} = 0.0248$, $p_\text{QBX} = 10$,
$\kappa_\text{QBX} = 10$ and $N_\text{U} = 2$, $d_\text{U1} =
0.0944$, $d_\text{U2} = 0.0386$;
and for the rods $d_\text{QBX} = 0.0146$, $p_\text{QBX} = 40$,
$\kappa_\text{QBX} = 25$ and $N_\text{U} = 2$, $d_\text{U1} =
0.0676$, $d_\text{U2} = 0.0237$.
The maximal error in the centre plane is plotted for varying
quadrature tolerance in Figure~\ref{fig:ch71-SI-conv}. This shows
that the error more or less follows the tolerance, as expected.


\FloatBarrier
\subsection{Geometry 2: Two spheroids in a pipe}

The second geometry consists of a pipe of radius~0.3, discretized
using $5 \times 10$ patches with $6 \times 6$ grid points each.
Inside the pipe are two spheroids with semiaxes $a=0.05$ and
$c=0.1$, discretized with parameters $n_\theta = 36$ and
$n_\varphi = 25$ (900~grid points per spheroid).

\begin{figure}[h!]
  \centering\small
  \begin{minipage}{0.41\textwidth}%
    \centering
    \includegraphics[width=\textwidth,trim=10mm 32mm 10mm 32mm,clip]%
    {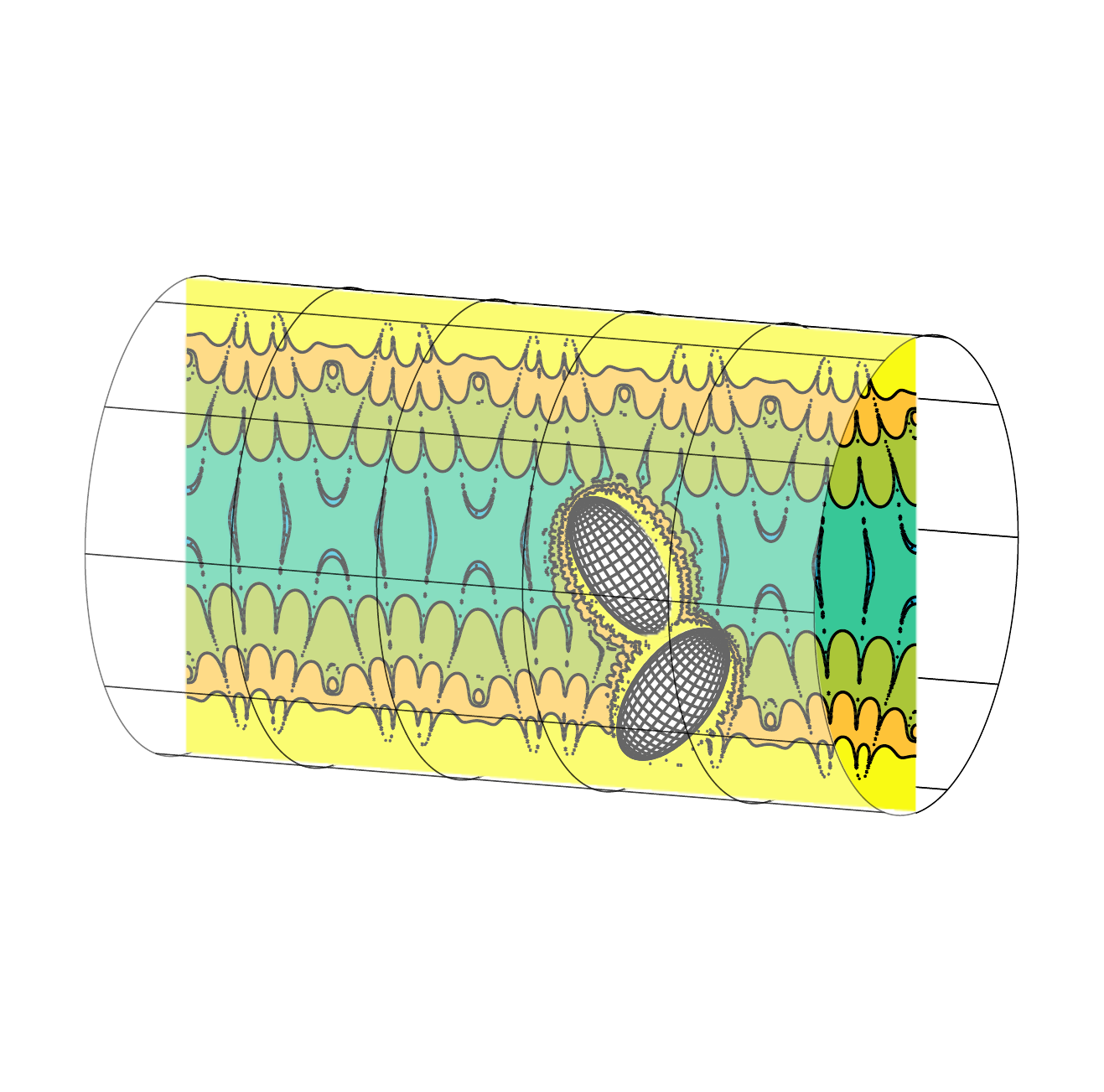}\\%
    (a) Direct quadrature error\\\vphantom{tolerance}
  \end{minipage}%
  \begin{minipage}{0.05\textwidth}%
    \centering
    \phantom{.}
  \end{minipage}%
  \begin{minipage}{0.41\textwidth}%
    \centering
    \includegraphics[width=\textwidth,trim=10mm 32mm 10mm 32mm,clip]%
    {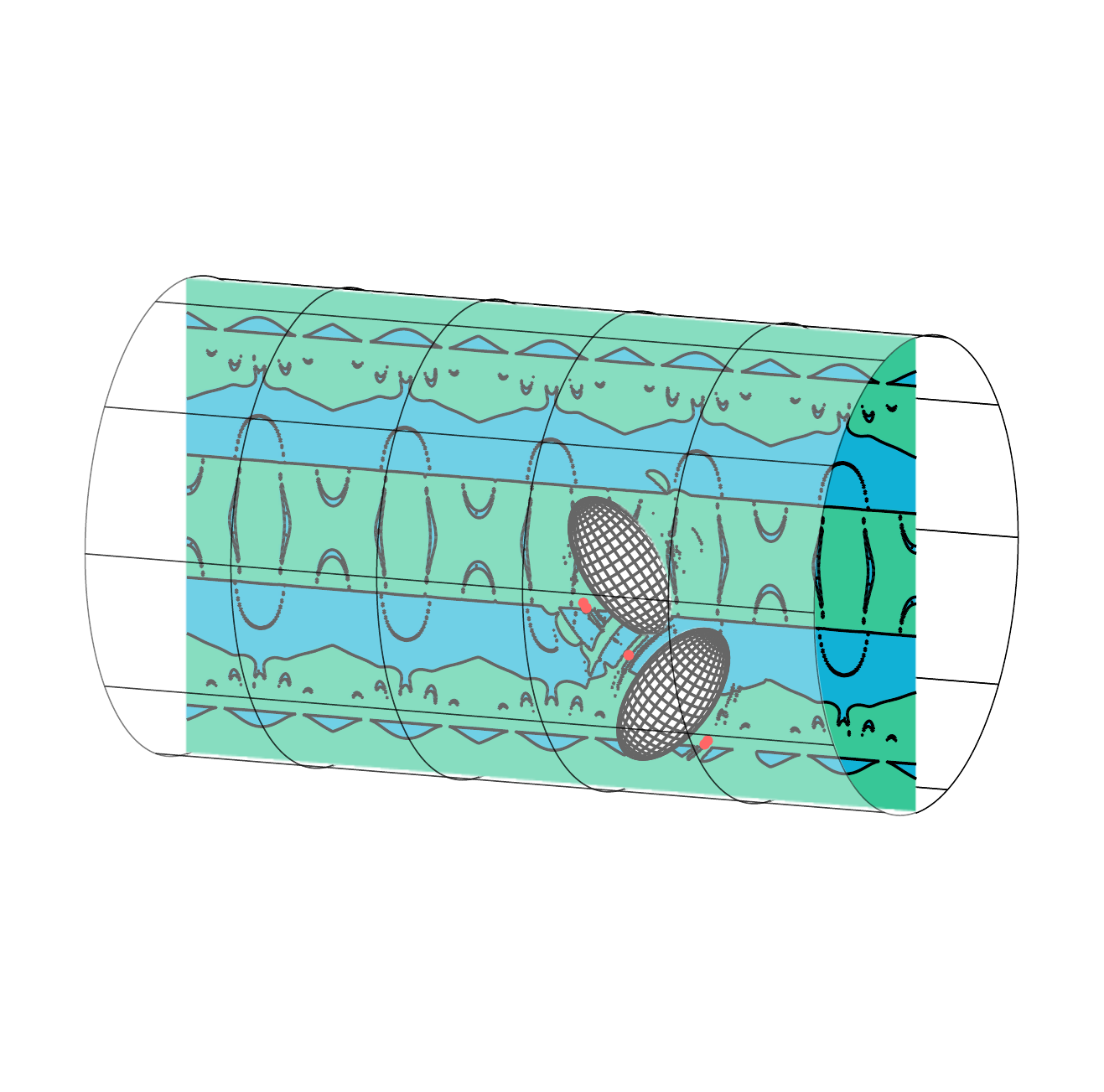}\\%
    (b) Special quadrature error,\\tolerance \num{e-6}
  \end{minipage}%
  \begin{minipage}{0.13\textwidth}%
    \centering
    \vphantom{\textbf{Upq}}\hspace*{0.001pt}\\%
    \includegraphics[scale=1]{fig/colorbar.pdf}\\
    \vphantom{(c)}
  \end{minipage}%
  \caption{Stresslet identity error in the centre plane, for
  geometry~2, in (a) using direct quadrature and in (b) using
  combined special quadrature with tolerance \num{e-6}. In (b),
  the tolerance is exceeded in 5~points (the evaluation grid
  has $500 \times 500$ points), marked red; the largest
  error in the slice is \num{1.323e-6}.}
  \label{fig:ch71-B-direct}
\end{figure}

We select the error tolerance $\etol = \num{e-6}$.
The parameters selected according to
section~\ref{sec:param-general-summary} are, for the pipe
$d_\text{QBX} = 0.0614$, $p_\text{QBX} = 12$, $\kappa_\text{QBX} = 10$
and $N_\text{U} = 2$, $d_\text{U1} = 0.222$, $d_\text{U2} = 0.0888$;
and for the spheroids $d_\text{QBX} = 0.0153$, $p_\text{QBX} = 27$,
$\kappa_\text{QBX} = 15$ and $N_\text{U} = 2$, $d_\text{U1} =
0.0568$, $d_\text{U2} = 0.0235$.
The error when using these parameters is shown in
Figure~\ref{fig:ch71-B-direct}~(b), together with the direct
quadrature error in Figure~\ref{fig:ch71-B-direct}~(a). Note that
we have selected a much lower resolution for the pipe in
comparison to the walls in geometry~1
(section~\ref{sec:res1-special-quadrature-geometry-1}), which is
reflected in the larger threshold distances compared to
\eqref{eq:param3-params-tol6}.


\section{Solving the boundary integral equation}
\label{sec:res2-convergence}

Here, we investigate how the special quadrature tolerance
influences the accuracy of the numerical method, i.e.\ when
solving the boundary integral equation. We will use the mobility
problem as our model problem, and apply the force $\vec{F} =
(0, 0, -1)$ to all particles, with zero torque and no background
flow. In order to get the expected accuracy when solving the
boundary integral equation, the double layer density must
be well-resolved by the geometry discretization. It turns out
that for elongated particles, the density and how easy it is to
resolve depends heavily on the number of completion sources
$N_\text{src}$ (defined in section~\ref{sec:bif}). Therefore, we
begin in section~\ref{sec:res2-convergence-Nsrc} by investigating
how large $N_\text{src}$ must be to ensure that the density is
well-resolved for a given discretization. Then, in
section~\ref{sec:res2-convergence-quad-tol}, we study how the
special quadrature tolerance influences the accuracy of the
method.

\subsection{Selecting the number of completion sources}
\label{sec:res2-convergence-Nsrc}

\begin{figure}[H]
  \vspace*{-9mm}
  \centering\small
  \begin{minipage}[b]{0.5\textwidth}
    \centering
    \includegraphics[scale=1]{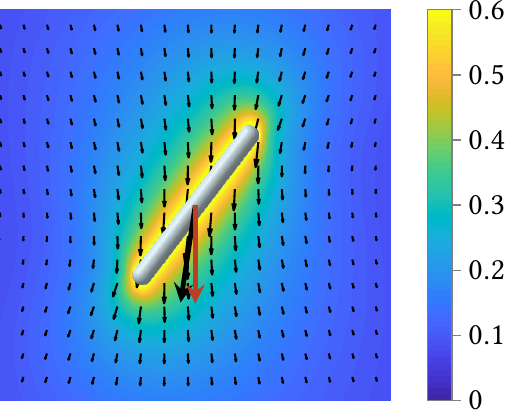}\\[1cm]%
    \hspace{-12mm}(a) Flow field
  \end{minipage}%
  \begin{minipage}[b]{0.5\textwidth}
    \centering
    \includegraphics[scale=0.8]{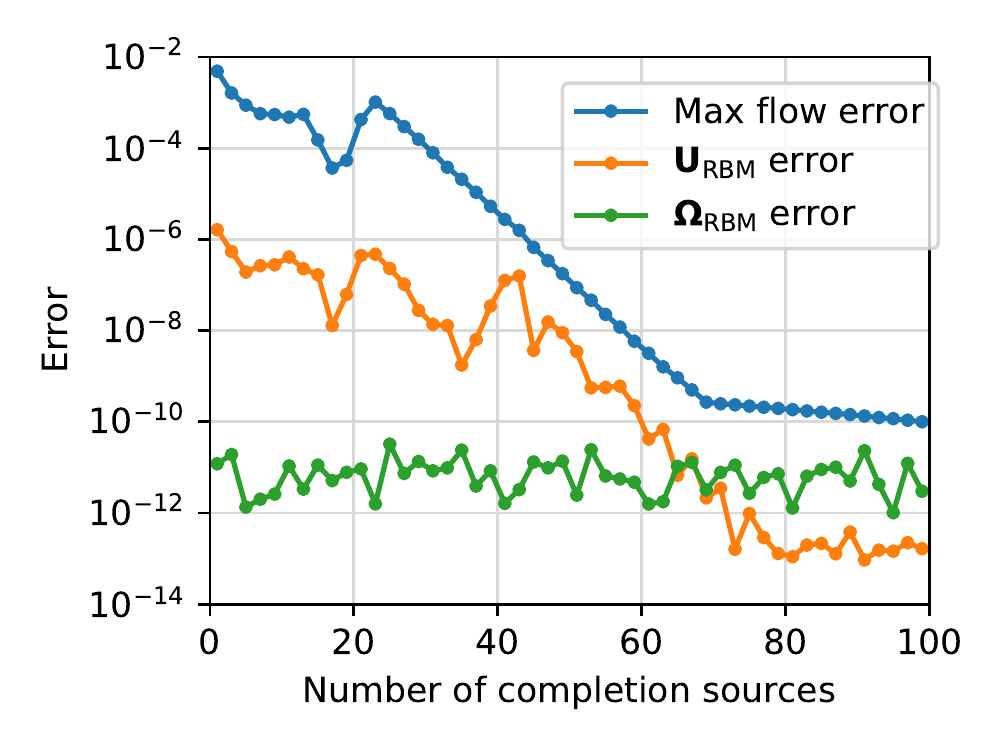}\\%
    \hspace{12mm}(b) Error contribution
  \end{minipage}%
  \caption{(a) Flow field resulting from the mobility problem for a
  single rod in free space (colour indicates velocity magnitude,
  small black arrows indicate velocity direction). The large red
  arrow indicates the applied force, and the large black arrow
  indicates the velocity of the rod (not to scale with the small
  arrows). (b) Contribution to the absolute error from the way the
  completion sources are distributed, as a function of
  $N_\text{src}$ (for a rod particle of aspect ratio~10). The error
  is estimated as the difference to a reference solution with
  $N_\text{src}=135$. Note that the max flow error flattens out
  around \num{e-9}, the special quadrature error tolerance.}
  \label{fig:ch72-1-comb}
\end{figure}

To study the influence of $N_\text{src}$, we consider a single
rod particle with length~$L=0.5$ and radius~$R=L/20$ (aspect
ratio~10) in free space, shown in Figure~\ref{fig:ch72-1-comb}~(a)
together with the flow field resulting from the force $\vec{F} =
(0, 0, -1)$. We now solve this mobility problem for varying
$N_\text{src}$, with the special quadrature error tolerance fixed
to \num{e-9} here.%
\footnote{
  The rod particle is discretized as in
  section~\ref{sec:res1-special-quadrature}, and the special
  quadrature parameters for $\etol=\num{e-9}$ are $r_\text{QBX}/h=1$,
  $d_\text{QBX} = 0.0119$, $p_\text{QBX}=45$, $\kappa_\text{QBX}
  = 25$, and $N_\text{U}=3$, $d_\text{U1} = 0.0810$, $d_\text{U2}
  = 0.0271$, $d_\text{U3} = 0.0164$.
}

\begin{figure}[H]
  \centering
  \includegraphics[scale=1]{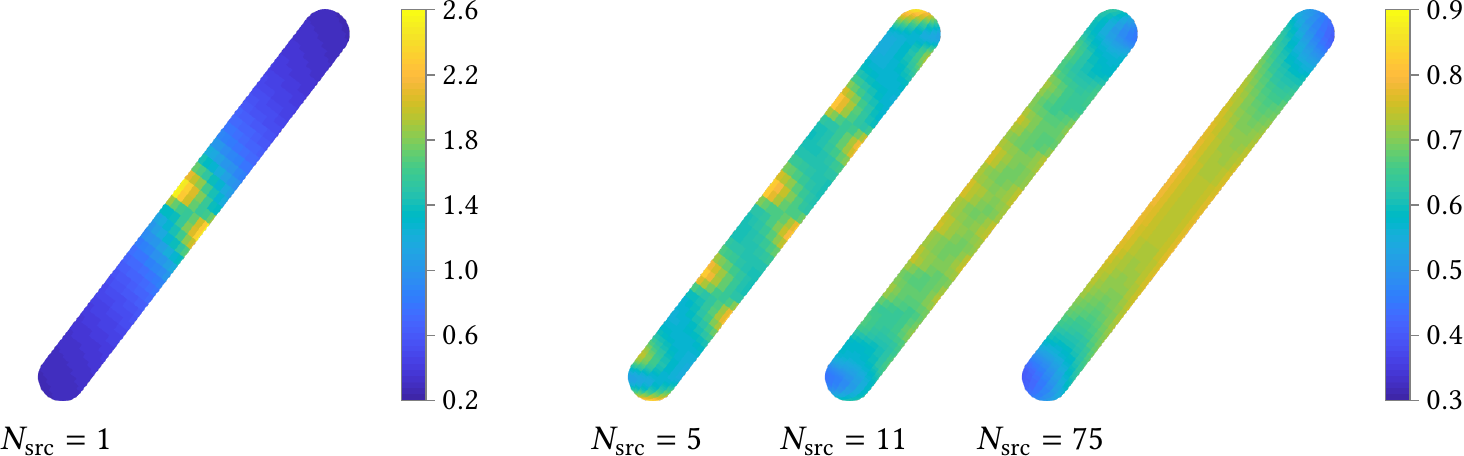}
  \caption{The magnitude of the completion flow
  $\vec{\Vp}^{(\alpha=1)}[\vec{F}, \vec{0}](\vec{x})$ on the
  surface of the rod, for a few different values of
  $N_\text{src}$. Since the background flow is zero, this is
  exactly the right-hand side of the boundary integral
  equation~\eqref{eq:mobility-bie}. Note that the colour scale is
  different for $N_\text{src}=1$ compared to the other values.}
  \label{fig:ch72-2-comp}
\end{figure}

The completion flow $\vec{\Vp}^{(\alpha)}$, which appears in the
right-hand side of the boundary integral
equation~\eqref{eq:mobility-bie}, will change drastically as
$N_\text{src}$ grows from small values, as shown in
Figure~\ref{fig:ch72-2-comp}; the completion flow becomes
increasingly smoother as $N_\text{src}$ increases.
Naturally, this means that the density~$\vec{q}$ itself will
change as $N_\text{src}$ grows. However, the real physical
quantities -- the particle velocity and the flow field -- should
not change since the net force and torque on the particle does
not change. Thus, these physical quantities can be used to gauge
how $N_\text{src}$ affects the accuracy of the solution.
As Figure~\ref{fig:ch72-1-comb}~(b) shows, the effect is quite large,
and most pronounced in the fluid flow velocity (the blue curve).
(For this problem, the magnitude of the fluid flow velocity and
particle velocity $\vec{U}_\text{RBM}$ are both around~0.6, while
the angular velocity $\vec{\Omega}_\text{RBM}$ is zero.)

Thus, it is important to select $N_\text{src}$ high enough for
the error in Figure~\ref{fig:ch72-1-comb}~(b) to satisfy the error
tolerance. The effect of $N_\text{src}$ on the accuracy is
stronger the more elongated the particle is; for particles with
low aspect ratio, $N_\text{src}=1$ may be sufficient. The effect
is very similar for the resistance problem, to the degree that
the max flow error in Figure~\ref{fig:ch72-1-comb}~(b) can be used to
determine $N_\text{src}$ for both the mobility problem and
resistance problem. Note that $N_\text{src}$ does not affect the
size of the linear system, i.e.\ \eqref{eq:resistance-disc} or
\eqref{eq:mobility-disc}.

\pagebreak

\vspace*{-3mm}
\subsection{Effect of the special quadrature on the accuracy}
\label{sec:res2-convergence-quad-tol}

\begin{figure}[t!]
  \centering\small
  \begin{minipage}{0.5\textwidth}%
    \centering
    \includegraphics[scale=1]{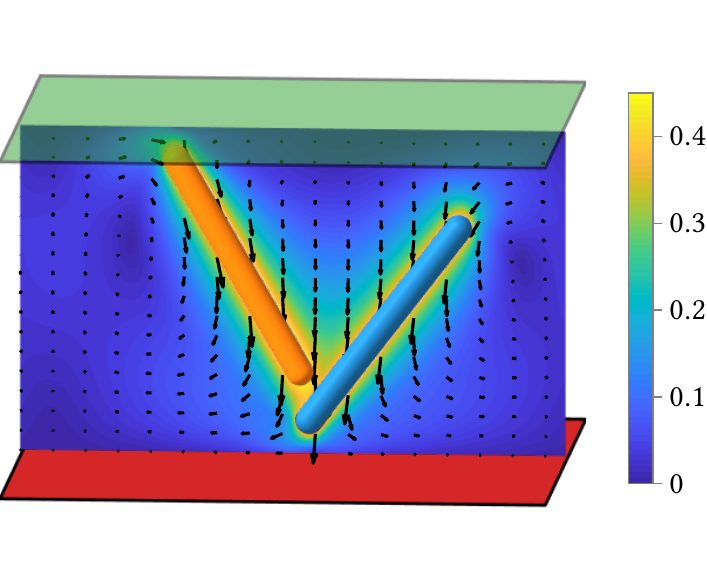}\\%
    \hspace{-12mm}(a) Geometry and flow field
  \end{minipage}%
  \begin{minipage}{0.5\textwidth}%
    \centering
    \includegraphics[scale=0.8]{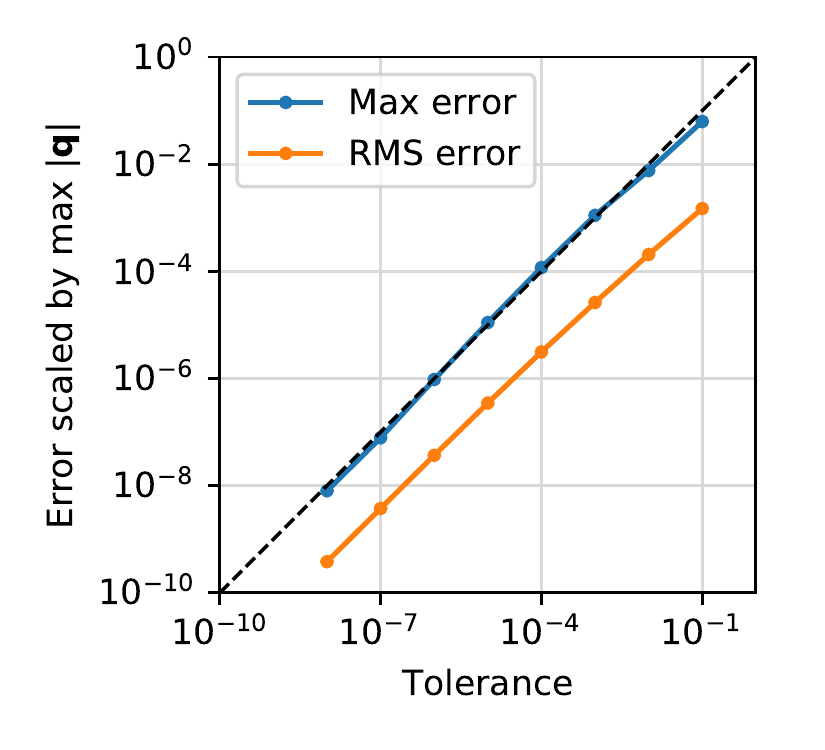}\\%
    \hspace{8mm}(b) Flow field error
  \end{minipage}%
  \caption{(a) Flow field from the periodic mobility problem with
  force $\vec{F} = (0,0,-1)$ applied to both particles (colour
  indicates velocity magnitude, arrows indicate velocity
  direction). (b) Maximal and root-mean-square (RMS) flow field
  error in the centre plane (estimated using a reference solution
  with tolerance \num{e-9}), scaled by the maximal density
  magnitude $\max\,\absi{\vec{q}} \approx 5.4$. The dashed line
  indicates Scaled error = Tolerance.}
  \label{fig:ch72-5-geom}
\end{figure}

\begin{figure}[t!]
  \centering\small
  \begin{minipage}{0.5\textwidth}%
    \centering
    \includegraphics[scale=0.8]{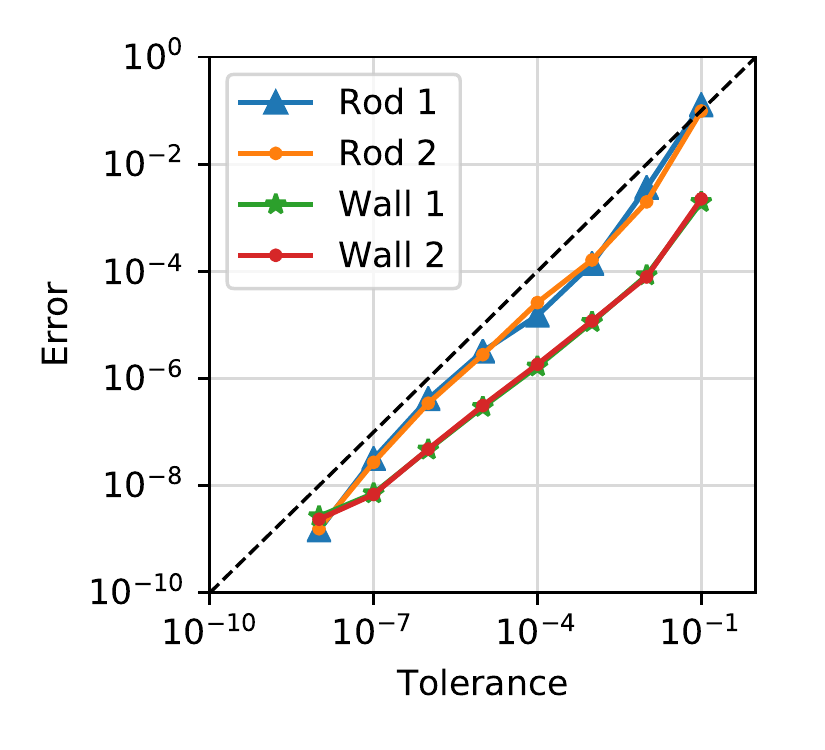}\\%
    \hspace{8mm}(a) Layer density error
  \end{minipage}%
  \begin{minipage}{0.5\textwidth}%
    \centering
    \includegraphics[scale=0.8]{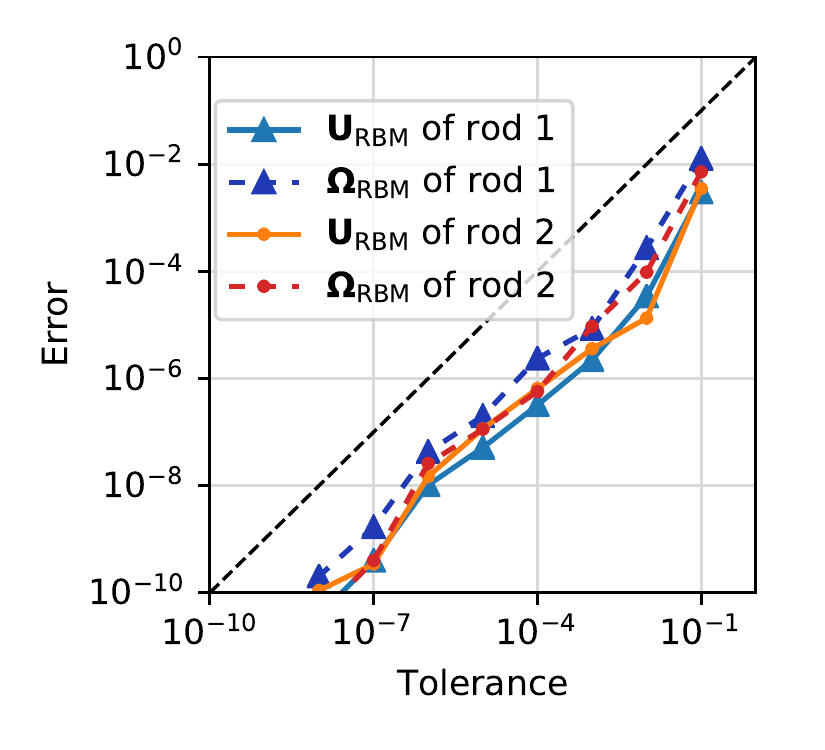}\\%
    \hspace{8mm}(b) Particle velocity error
  \end{minipage}%
  \caption{(a) Maximal absolute error of the layer density on each
  geometrical object (colors as in
  Figure~\ref{fig:ch72-5-geom}~(a)), estimated using a reference
  solution with tolerance \num{e-9}. (b) Absolute error of the
  particle velocities.}
  \label{fig:ch72-4-conv}
\end{figure}

We continue to study the mobility problem, but now add another
rod particle and a pair of plane walls, as shown in
Figure~\ref{fig:ch72-5-geom}. We fix $N_\text{src} = 65$, which
was enough to get the error below \num{e-9} in the previous
problem. The walls are discretized using $22 \times 22$ patches
with $8 \times 8$ grid points each (30\,976 grid points per
wall), and the rod particles are discretized as in
section~\ref{sec:res1-special-quadrature}. We set the special
quadrature tolerance to different values between \num{e-1} and
\num{e-8}, solve the mobility problem, and compute the flow field
and particle velocities. The errors in the flow field, density and
particle velocities are estimated using a reference solution with
special quadrature tolerance \num{e-9}; these are shown in
Figures~\ref{fig:ch72-5-geom}~(b) and
\ref{fig:ch72-4-conv}~(a)--(b). Note that the tolerance sets the
flow field error relative to $\max\,\absi{\vec{q}}$ quite
accurately; the density and particle velocity errors are even
smaller. We would like to point out that the value of the scale
factor $\max\,\absi{\vec{q}}$ is not known a priori, but it is of
course known after having solved the boundary integral equation.

It should be noted that the error cannot be expected to follow
the tolerance unless the density is well-resolved by the geometry
discretization, since otherwise the interpolated density will be
inaccurate. It has been observed that the density becomes hard to
resolve, with either sharp peaks or high-frequency oscillations,
when particles come very close to each other or the walls (where
``very close'' is measured relative to the grid resolution).
Thus, one may be forced to increase the grid resolution in these
cases.

For elongated particles, the set of matrices $\mat{M}_i$,
$i=1,\ldots,n_\theta/2$, which are precomputed for offsurface QBX
tends to become quite large for strict quadrature tolerances. The
reason for this is that with particle-global QBX, the whole set
of matrices consists of
\begin{equation}
  3 n_\varphi n_\theta^2 (p_\text{QBX}+1)(p_\text{QBX}+2)
\end{equation}
complex numbers, i.e.\ it is quadratic in both $p_\text{QBX}$ and
$n_\theta$ (where $n_\theta$ is the number of grid points in the
axial direction,\pagebreak\ which we define as $2n_1 + n_2$ for rod
particles). For elongated particles, $n_\theta$ tends to be
large; for example, for the rods considered in this section
(aspect ratio~10),
$n_\theta=130$ and $n_\varphi = 18$. The set of matrices for
tolerance~\num{e-8} ($p_\text{QBX}=40$) then takes up around
25~gigabytes when stored in double precision, while for
tolerance~\num{e-6} ($p_\text{QBX}=28$) the matrices take up
around 13~gigabytes. To reduce the size of the matrices in this
situation, a local patch-based discretization could be used also
for the particles, in the same way it is already used for the
walls. This would reduce the number of grid points included in
the special quadrature and thus the size of the matrices.

\section{Computational complexity and computation of streamlines}
\label{sec:res3-packed-rods}

\subsection{Computational complexity of the method}
\label{sec:res3-packed-rods-complexity}

The computational cost of our special quadrature method is
quadratic in the number of grid points per particle (or patch),
but linear in the number of particles (patches) if their
discretization is kept fixed. For the Spectral Ewald method, the
computational cost scales like $O(N \log N)$, where $N$ is the
number of unknowns in the system (i.e.\ three times the total
number of grid points), assuming that the number of grid points
within a ball of radius $r_\text{c}$ does not change.
In other words, for fixed grid resolution and particle
concentration, the time required per GMRES iteration when solving
the boundary integral equation is expected to scale like $O(N
\log N)$.

\begin{figure}[h!]
  \centering\small
  \begin{minipage}[b]{0.3\textwidth}%
    \centering
    \includegraphics[scale=0.4,trim=7.1cm 0cm 7.5cm 0cm,clip]%
        {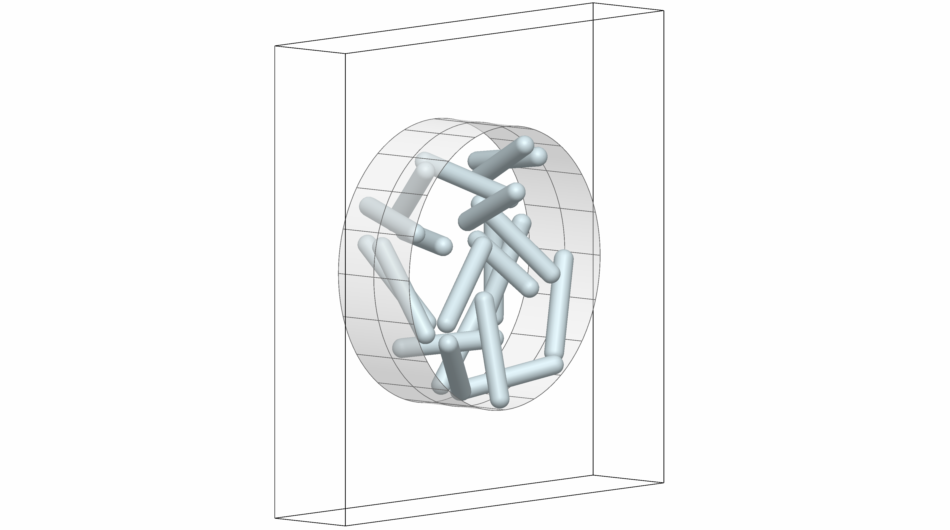}\\%
    (a) 1~segment
  \end{minipage}%
  \begin{minipage}[b]{0.7\textwidth}%
    \centering
    \includegraphics[scale=0.4]{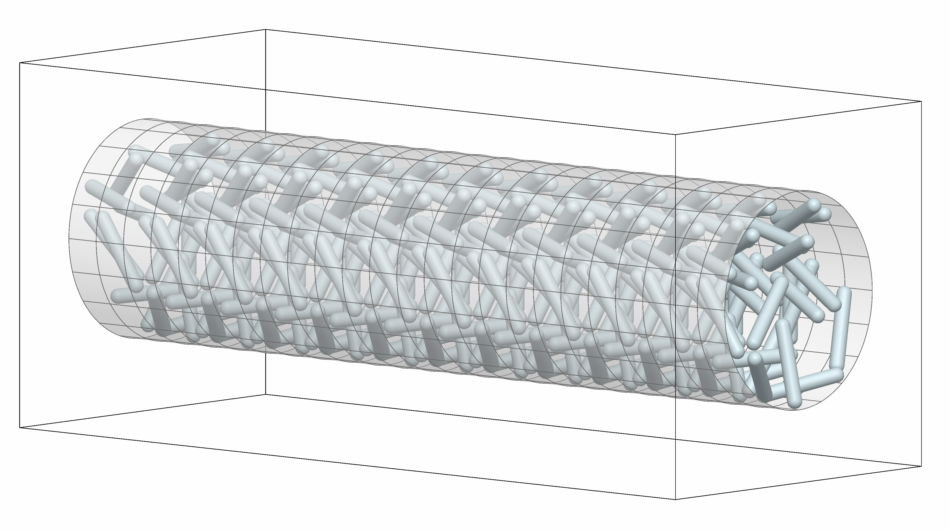}\\%
    (b) 12~segments
  \end{minipage}%
  \caption{The geometry in
  section~\ref{sec:res3-packed-rods-complexity} is made up of
  stacked identical segments, where each segment contains 20~rods
  and $2 \times 20$ pipe patches. In total, there are 27\,360~grid points
  and 82\,080~unknowns in each segment.}
  \label{fig:res3-geometry}
\end{figure}

To test this scaling, we consider a problem with many rod
particles confined in a pipe, shown in
Figure~\ref{fig:res3-geometry}: one segment (a) consists of a
pipe segment of radius~0.3 and length~0.2 confined in a periodic
cell of size $\vec{B} = (0.2, 1, 1)$, with 20~rods of
length~$L=0.25$ and radius~$R=L/12$ (aspect ratio~6) inside the pipe.%
\footnote{
  The discretization, special quadrature parameters and Spectral
  Ewald parameters are fixed as follows. Each pipe patch has $6
  \times 6$ grid points, and each rod is discretized using $n_1 =
  16$, $n_2 = 40$ and $n_\varphi = 18$ (1296~grid points per
  rod). The special quadrature parameters are selected for
  tolerance $\etol = \num{e-4}$, and they are for the pipe
  $d_\text{QBX} = 0.0319$, $p_\text{QBX} = 8$, $\kappa_\text{QBX}
  = 10$ and $N_\text{U} = 1$, $d_\text{U1} = 0.0720$;
  and for the rods $d_\text{QBX} = 0.0110$, $p_\text{QBX} = 19$,
  $\kappa_\text{QBX} = 10$ and $N_\text{U} = 1$, $d_\text{U1} = 0.0254$.
  The Spectral Ewald parameters are $\xi=52.954$, $r_\text{c}=0.0897$,
  $P=16$, and the uniform grid for the Fourier-space part
  has $32n_\text{s} \times 160 \times 160$ grid points, where $n_\text{s}$ is
  the number of segments. Furthermore, we use $N_\text{src} = 45$
  completion sources per rod.
}
This segment is replicated to create a longer pipe, up to 12
times the original length (shown in (b)), with the same grid point
concentration as the original segment.
For $1, 2, 3, \ldots, 12$ segments, we solve a resistance problem
in which all particles are stationary and a quadratic background
flow
\begin{equation}
  \vec{u}_\text{bg}(\vec{x}) = \gp{
    \frac{A^2 - x_2^2 - x_3^2}{A^2},%
    \hspace{.5em}0,\hspace{.5em}0
  }
\end{equation}
is applied, where $A=0.3$ is the radius of the pipe and
$(x_2, x_3) = (0,0)$ is its centre line.

\begin{figure}[h!]
  \centering\small
  \begin{minipage}[b]{0.5\textwidth}%
    \centering
    \includegraphics[scale=0.76]{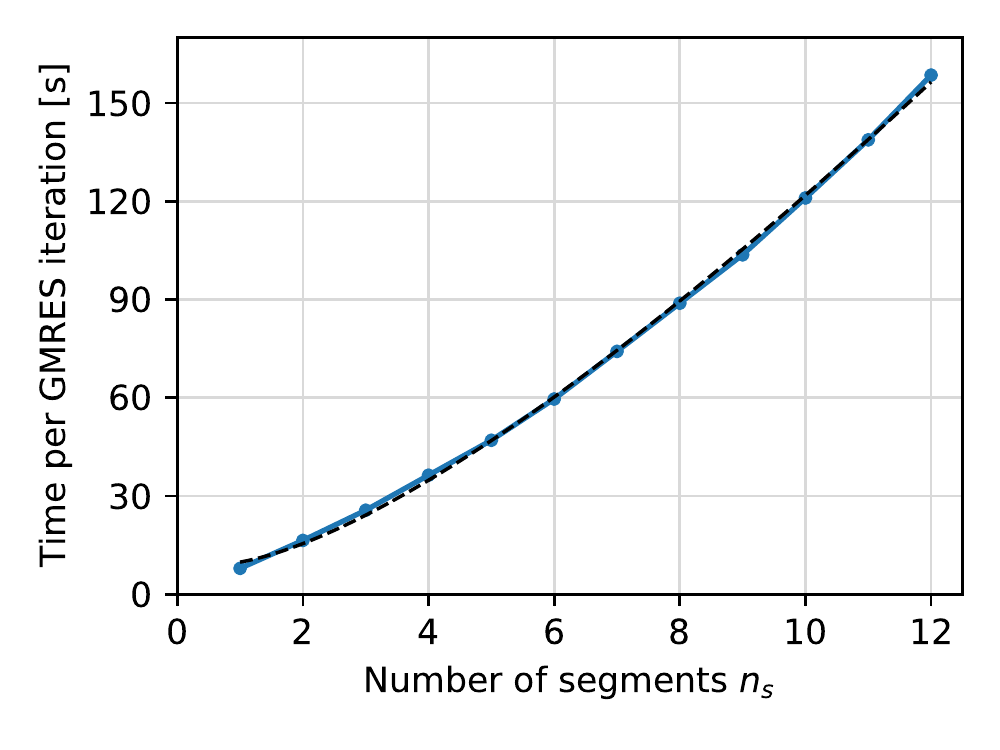}\\[-3mm]%
    \hspace*{6mm}(a) Time per GMRES iteration
  \end{minipage}%
  \begin{minipage}[b]{0.5\textwidth}%
    \centering
    \includegraphics[scale=0.76]{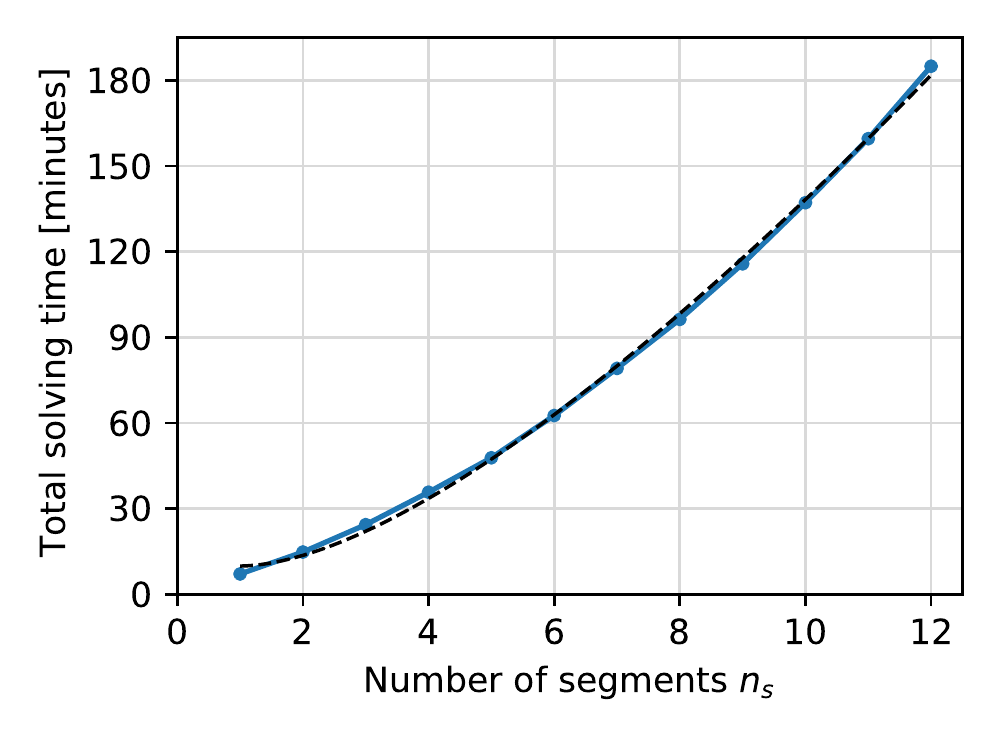}\\[-3mm]%
    \hspace*{6mm}(b) Total solving time
  \end{minipage}%
  \caption{Time required to solve the resistance problem for the
  geometry in Figure~\ref{fig:res3-geometry} with GMRES tolerance
  \num{e-6} in MATLAB, (a) per GMRES iteration and (b) in total.
  With 82\,080 unknowns per segment, the number of unknowns
  ranges from 82\,080 to 984\,960. The dashed curves are
  least-squares fits of $T = A\,n_\text{s} \log n_\text{s} +
  B\,n_\text{s} + C$ to the data, where $n_\text{s}$ is the number of segments.
  In (a), $A=5.8$, $B=-2.4$, $C=12.2$ (seconds), and in (b),
  $A=9.0$, $B=-8.7$, $C=18.6$ (minutes). Thus, the time scales as
  $O(n_\text{s} \log n_\text{s})$.}
  \label{fig:res3-time-scaling}
\end{figure}

As seen in Figure~\ref{fig:res3-time-scaling}~(a), the time per
GMRES iteration follows the expected scaling $O(N \log N)$. Since
the structure of the linear system changes as the number of
segments $n_\text{s}$ grows, the number of GMRES iterations grows
slightly with $n_\text{s}$. However, this growth is slow enough
for the total solving time to also follow the scaling $O(N \log
N)$, as shown in Figure~\ref{fig:res3-time-scaling}~(b).

\subsection{Streamline computation}

\begin{figure}[h!]
  \centering
  \includegraphics[scale=1]{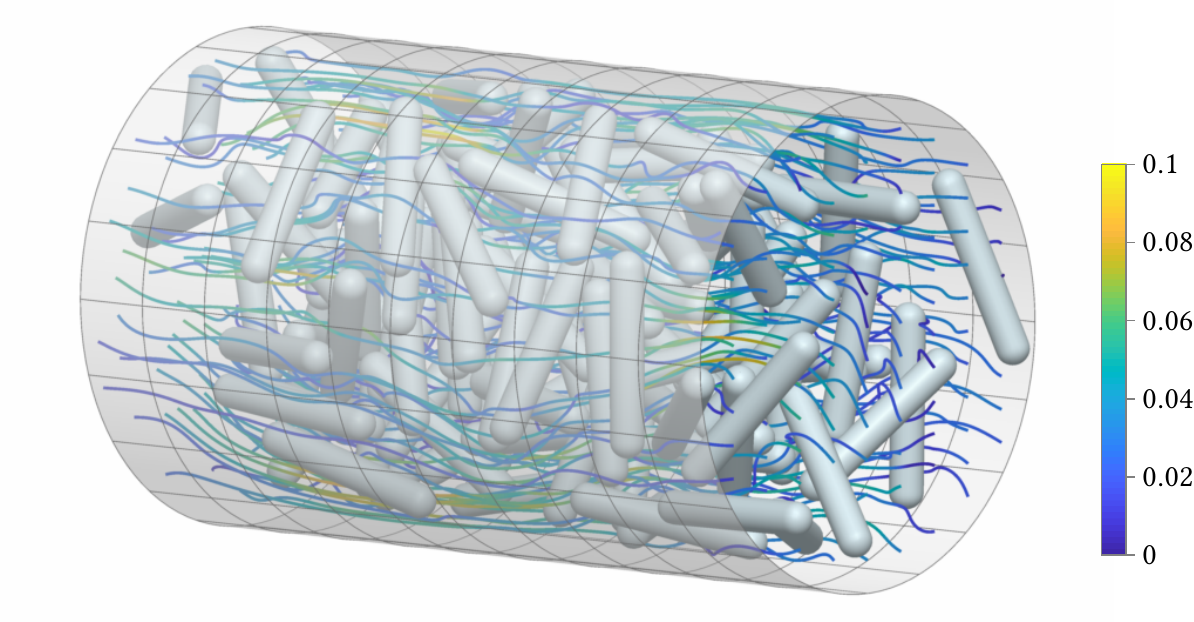}
  \caption{Streamlines for the resistance problem (colour
  indicate velocity magnitude).}
  \label{fig:res3-streamlines}
\end{figure}

\begin{figure}[h!]
  \centering
  \begin{minipage}{0.333\textwidth}%
    \centering
    \includegraphics[width=\textwidth]{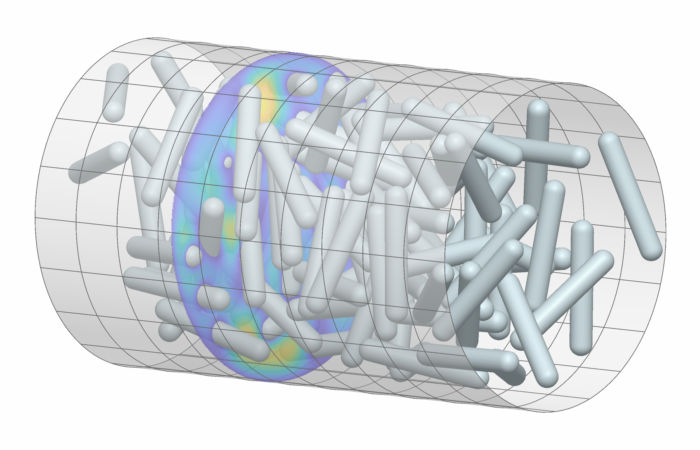}
  \end{minipage}%
  \begin{minipage}{0.333\textwidth}%
    \centering
    \hspace*{-30pt}{\small\bfseries Streamwise velocity}\\[-2pt]%
    \includegraphics[scale=1]{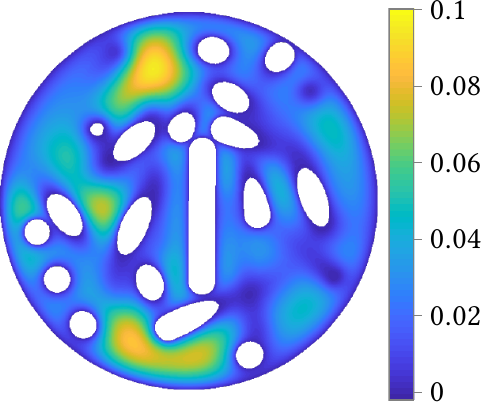}
  \end{minipage}%
  \begin{minipage}{0.333\textwidth}%
    \centering
    \hspace*{-30pt}{\small\bfseries Transverse velocity}\\%
    \includegraphics[scale=1]{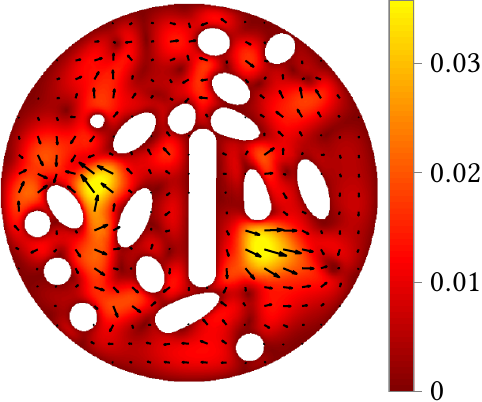}
  \end{minipage}
  \caption{Flow field shown in a slice at streamwise position $x=2/3$.}
  \label{fig:res3-slice1}
\end{figure}

In the postprocessing step (step~4 of the method summary),
streamlines may be computed to visualize the flow field. When
using the Spectral Ewald method, the Fourier-space part on the
uniform grid can be reused to reduce the computation time, as
described in appendix~\ref{app:streamlines}. Here, we compute
streamlines for a periodic resistance problem with 100~rods in a
pipe segment of length~1, otherwise identical to the problem in
section~\ref{sec:res3-packed-rods-complexity} (including
all parameters). Figure~\ref{fig:res3-streamlines} shows 95
streamlines; a typical streamline consists of around 3000~points
and takes around 2~minutes to compute (i.e.\ around 0.04~seconds
per time step). A slice of the same flow field is shown in
Figure~\ref{fig:res3-slice1}.



\FloatBarrier
\chapter{Effects of nonsmooth geometries}
\label{sec:smoothness}

\begin{figure}[H]
  \centering
  \begin{minipage}[b]{0.26\textwidth}%
    \centering
    \begin{tikzpicture}[inner sep=0pt]
      \node [anchor=south] at (0,0)
          {\rotatebox{90}{\includegraphics[trim=16mm 4mm 11mm 2mm,clip,width=4.7cm]%
                          {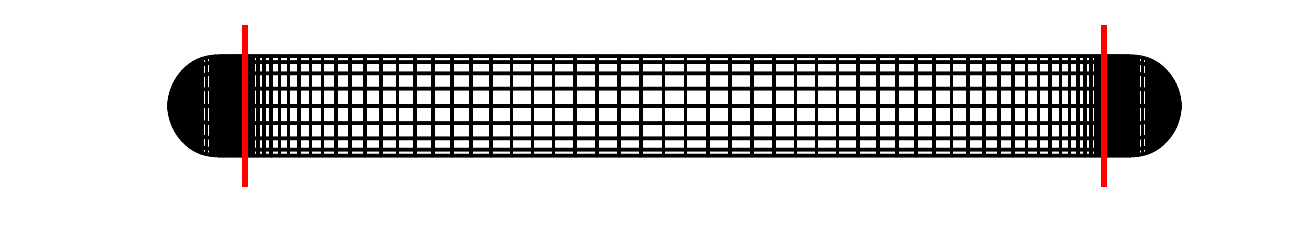}}};
      \node [anchor=south] at (2,0)
          {\rotatebox{90}{\includegraphics[trim=16mm 4mm 11mm 2mm,clip,width=4.7cm]%
                          {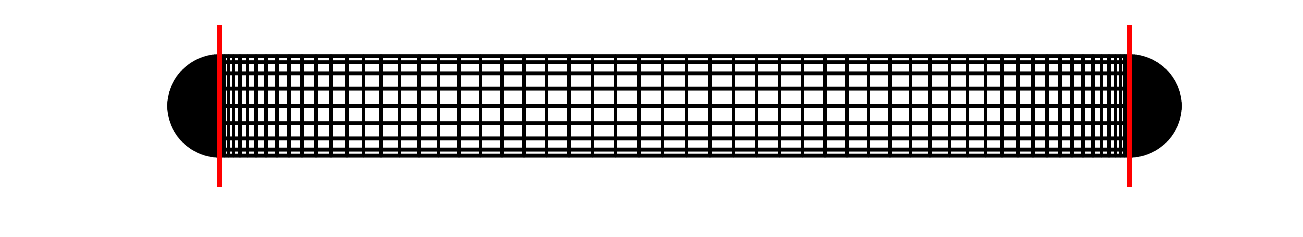}}};
      \node[align=center] at (0,-0.5) {Smooth\\rod};
      \node[align=center] at (2,-0.5) {Nonsmooth\\rod};
    \end{tikzpicture}\\[10pt]%
    \small(a) Geometries
  \end{minipage}%
  \begin{minipage}[b]{0.74\textwidth}%
    \centering\small
    \includegraphics[width=\textwidth]{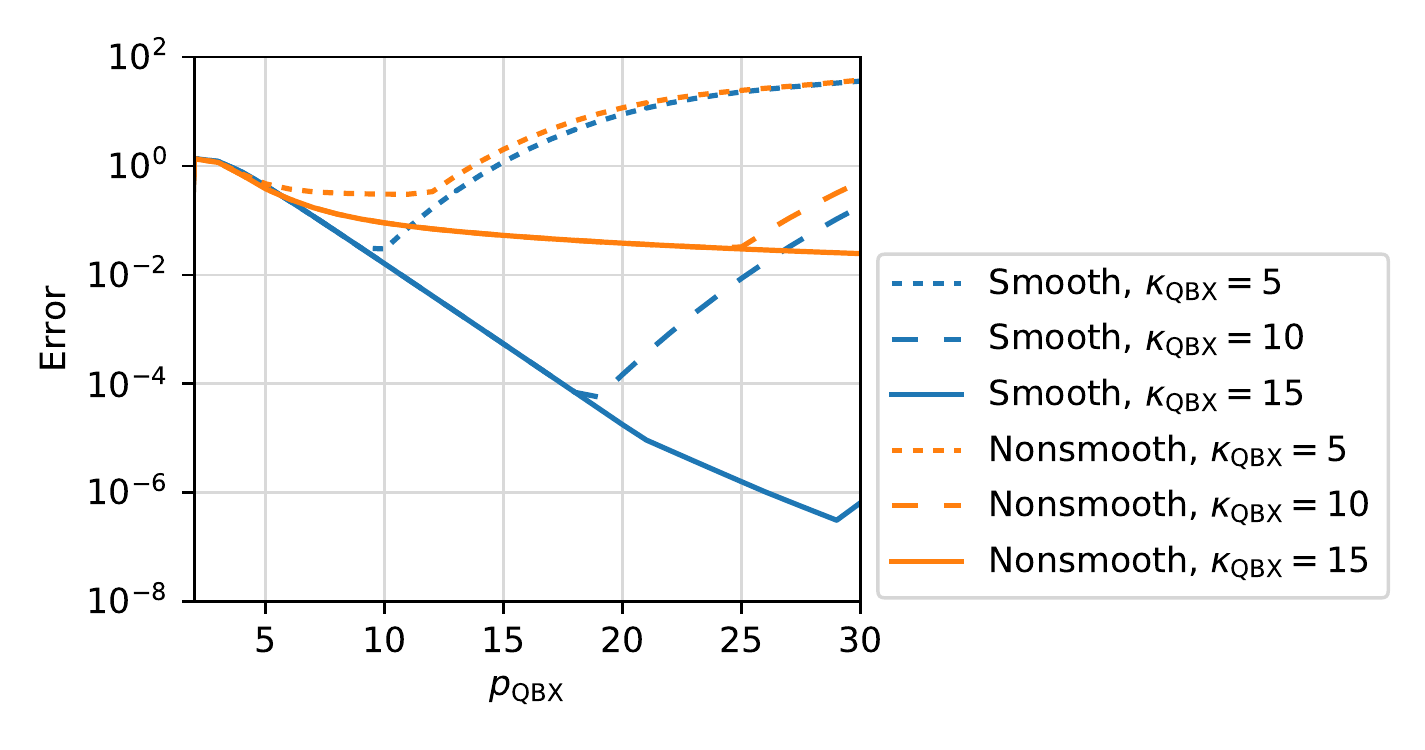}\\%
    \hspace*{-30mm}(b) QBX onsurface error
  \end{minipage}%
  \caption{(a) A smooth and a nonsmooth rod, both of length~$L$
  and radius~$R$. Note that each cap of the smooth rod (marked
  with red lines) has length $1.5 R$, while each cap of the
  nonsmooth rod has length~$R$. (b) Maximal onsurface QBX
  stresslet identity error for the smooth and nonsmooth rod, with
  expansion radius $r_\text{QBX} = h = 2\pi R / n_\varphi$
  fixed.}
  \label{fig:ch8-geom-qbx}
\end{figure}

So far, all geometrical objects considered in this paper have
been smooth. In fact, special care has been taken to ensure that
the rod particles, constructed in
appendix~\ref{app:smooth-rod}, are everywhere smooth. The reason
is that, as noted in \cite{epstein13}, the convergence of the
local expansions used in QBX depends on the smoothness of the
boundary close to the expansion centre. In this section, we
demonstrate this using two different rod particles: one smooth
and one nonsmooth, shown in Figure~\ref{fig:ch8-geom-qbx}~(a). The
rods are both of length~$L$ and radius~$R$, but the smooth rod is
constructed as in appendix~\ref{app:smooth-rod}, while the
nonsmooth rod consists of a cylinder of radius~$R$ and
length~$L-2R$ joined to two half-spherical caps of radius~$R$.
The nonsmooth rod is thus of class $C^1$, since the curvature is
discontinuous where the cylinder meets the spherical caps.

\begin{figure}[H]
  \vspace*{-10mm}
  \begin{tikzpicture}
    \node[inner sep=0pt] at (0,0)
    {\includegraphics[scale=0.45,trim=1.5cm 7.3cm 1.5cm 1cm,clip]{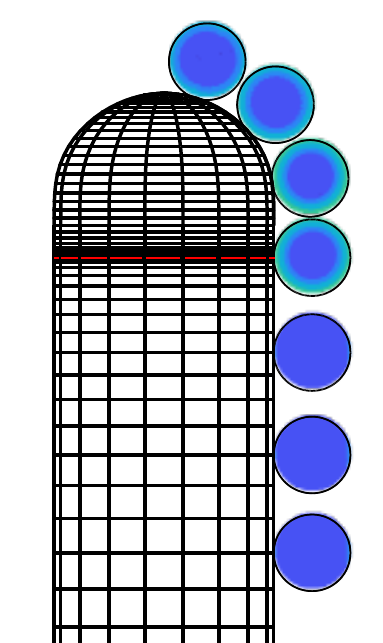}};
    \node[inner sep=0pt] at (3,0)
    {\includegraphics[scale=0.45,trim=1.5cm 7.3cm 1.5cm 1cm,clip]{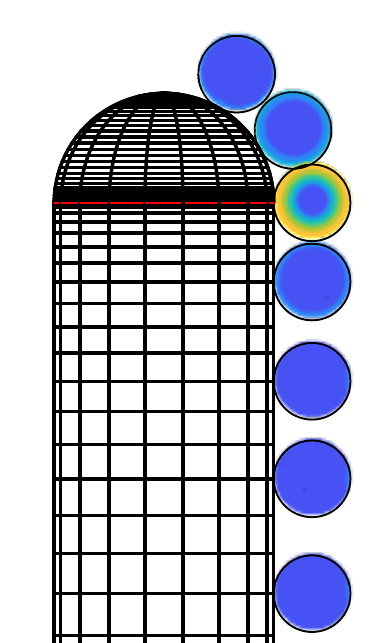}};
    \node[align=center] at (-0.1,-2.0) {Smooth\\rod};
    \node[align=center] at (2.9,-2.0) {Nonsmooth\\rod};
    \node at (5.5,-0.5) {\includegraphics{fig/colorbar.pdf}};
    \node[font=\small] at (1.4,-3.0) {(a)};
  \end{tikzpicture}\hfill%
  \begin{tikzpicture}
    \node[inner sep=0pt] at (0,0) {\includegraphics[scale=0.30]{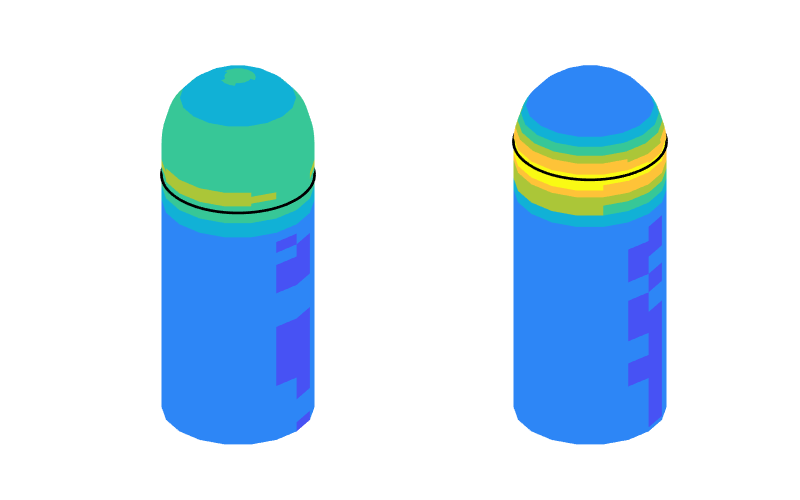}};
    \node[align=center] at (-1.27,-2.0) {Smooth\\rod};
    \node[align=center] at (1.50,-2.0) {Nonsmooth\\rod};
    \node at (3.8,-0.5) {\includegraphics{fig/colorbar.pdf}};
    \node[font=\small] at (0.1,-3.0) {(b)};
  \end{tikzpicture}
  \caption{(a) The QBX stresslet identity error shown in a slice
  through selected balls of convergence, for $r_\text{QBX} = h$,
  $p_\text{QBX} = 25$ and $\kappa_\text{QBX} = 15$. (b) The
  onsurface QBX error on the rods, with the same parameters as in
  (a).}
  \label{fig:ch8-qbx-demo}
\end{figure}

To illustrate the convergence of QBX, consider rods with $L/R =
20$ (aspect ratio~10), discretized using $n_1 = 35$, $n_2 = 60$,
$n_\varphi = 18$ as described in
section~\ref{sec:direct-quad-particles}. In
Figure~\ref{fig:ch8-geom-qbx}~(b), the onsurface QBX stresslet
identity error is plotted as a function of $p_\text{QBX}$, in the
same way as in section~\ref{sec:param-example2} (where this was
done for the smooth rod). Clearly, the convergence with respect
to $p_\text{QBX}$ is much worse for the nonsmooth rod compared to
the smooth one. The reason for this can be seen in
Figure~\ref{fig:ch8-qbx-demo}: the error decays extremely slowly
close to the boundary between the cylinder and the caps, where
the curvature is discontinuous. This is clearly a local effect,
since the convergence is fine a little bit away from the
discontinuity.

It should be noted that it is entirely possible to use QBX on a
nonsmooth geometry, but it requires special measures to be taken.
In \cite{klockner13}, QBX was applied to a geometry with a corner.
In that example, the discretization was dyadically refined around
the corner, to ensure that the layer potential appears locally
smooth on the scale of the discretization. The same approach
could likely be taken also for the nonsmooth rod particle, i.e.\
refining the grid dyadically around the discontinuity. However,
constructing the rod to be smooth in the first place has a clear
advantage in this case, since no grid refinement is needed.


\vspace*{-10mm}
\chapter{Conclusions}

We have presented a numerical method based on a boundary
integral formulation that can be used to simulate rigid
particles in Stokes flow with confining walls. A parameter
selection strategy has also been presented for the combined
special quadrature used in this method. We have demonstrated that
the error of the method is controlled by the special quadrature
tolerance as long as the layer density is well-resolved, and that
the method scales as $O(N \log N)$ in the number of unknowns $N$ for
fixed grid point concentration. This makes it
possible to simulate systems with a large number of particles.
The method can deal with particles and walls of different shapes;
we have here considered spheroids, rod particles, pipes and plane
walls, but it is straightforward to extend the method to any
smooth geometry with sufficient symmetry.

The method could be further improved for example by using local
patch-based quadrature for elongated particles to reduce the size
of the QBX matrices, and allowing the size of the wall patches to
be set adaptively so that the resolution can be focused where
particles are close to the wall.
It could also be useful to allow parameters such as $p_\text{QBX}$
to vary along the particle surface (in response to differences in
the convergence rate of the local expansions, as seen e.g.\ in
Figure~\ref{fig:ch8-qbx-demo}), and to allow the expansion
centres for QBX to be placed independently of the grid points of
the discretization, so that the centres can be placed closer to
the surface in order to decrease the expansion order
$p_\text{QBX}$.
Furthermore, if analytical quadrature error estimates were
available, these could replace the numerical experiments used to
select threshold distances and the QBX upsampling factor.

\FloatBarrier
\chapter*{Acknowledgements}

This work has been supported by the Göran Gustafsson Foundation
for Research in Natural Sciences and by the Swedish Research
Council under grant no. 2015-04998. The authors gratefully
acknowledge this support.

The authors are grateful to Dr.\ Ludvig af Klinteberg for
providing us with an implementation of QBX for spheroidal
particles, parts of which were reused in this work.

\appendices
\renewcommand*{\chapnumfont}{\chaptitlefont Appendix~}

\FloatBarrier
\chapter{The stresslet identity for plane walls and pipes}
\label{app:stresslet-identity}

Here we show that a variant of the stresslet identity
\eqref{eq:stresslet-identity} holds for a pair of parallel
infinite plane walls (in section~\ref{app:stresslet-identity-walls})
and an infinitely long pipe (in
section~\ref{app:stresslet-identity-pipe}).

\section{A pair of parallel plane walls}
\label{app:stresslet-identity-walls}

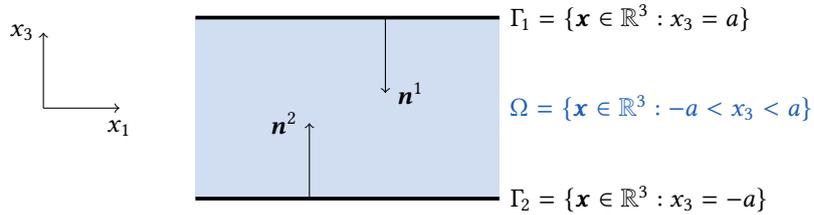
\begin{figure}[h!]
  \centering
  \begin{tikzpicture}
    \coordinate (A) at (-2,1.2);
    \coordinate (B) at (2,1.2);
    \coordinate (C) at (-2,-1.2);
    \coordinate (D) at (2,-1.2);
    \fill [myBlue!20] (C) rectangle (B);
    \draw [line width=1.4pt] (A) -- (B) node [right]
    {$\Gamma_1=\seti{\vec{x} \in \Real^3 : x_3=a}$};
    \draw [line width=1.4pt] (C) -- (D) node [right]
    {$\Gamma_2=\seti{\vec{x} \in \Real^3 : x_3=-a}$};
    \draw [->] (-4,0) -- (-4,1) node [left] {$x_3$};
    \draw [->] (-4,0) -- (-3,0) node [below=1pt] {$x_1$};
    \draw [->] (-0.5,-1.2) -- (-0.5,-0.2) node [left=1pt] {$\vec{n}^2$};
    \draw [->] (0.5,1.2) -- (0.5,0.2) node [right=1pt] {$\vec{n}^1$};
    \node [anchor=west,myBlue] at (2,0) {$\Omega=\seti{\vec{x}
    \in \Real^3 : -a < x_3 < a}$};
  \end{tikzpicture}
  \caption{Two parallel infinite planes $\Gamma_1$ and $\Gamma_2$.}
  \label{fig:appA-walls}
\end{figure}

Let $\Gamma_1$ and $\Gamma_2$ be two parallel infinite planes
oriented as shown in
Figure~\ref{fig:appA-walls}, one wall placed at $x_3=a$ and the
other at $x_3=-a$ for some $a>0$. Let the domain between the two
walls (which we will think of as the fluid domain) be denoted by
$\Omega$.

\noindent
Let $\tilde{\Gamma} = \Gamma_1 \cup \Gamma_2$, and let
$\tilde{\vec{q}} \in \Real^3$ be any constant vector. We shall
show that
\begin{equation}
  \label{eq:appA-stresslet-identity-1}
  \vec{\Dp}[\tilde{\Gamma}, \tilde{\vec{q}}](\vec{x}) =
  \begin{cases}
    -8\pi \tilde{\vec{q}}, & \text{if $\vec{x} \in \Omega$}, \\
    -4\pi \tilde{\vec{q}}, & \text{if $\vec{x} \in \tilde{\Gamma}$}, \\
    \vec{0}, & \text{otherwise}, \\
  \end{cases}
\end{equation}
where the double layer potential $\vec{\Dp}$ is given by
\eqref{eq:double-layer}.

Since $\tilde{\vec{q}}$ is constant and the normals are given by
$\vec{n}^1 = (0,0,-1)$ and $\vec{n}^2 = (0,0,1)$, we can write
\begin{align}
  \Dp_i[\tilde{\Gamma},\tilde{\vec{q}}](\vec{x}) &=
  \tilde{q}_j
  \int_{\tilde{\Gamma}} T_{ijk}(\vec{x}-\vec{y}) n_k(\vec{y}) \, \D S(\vec{y})
  \\
  \label{eq:appA-double-layer-step}
  &=
  -\tilde{q}_j
  \underbrace{\int_{\Gamma_1} T_{ij3}(\vec{x}-\vec{y})
  \, \D S(\vec{y})}_{J^1_{ij}(\vec{x})}
  \:+\:
  \tilde{q}_j
  \underbrace{\int_{\Gamma_2} T_{ij3}(\vec{x}-\vec{y})
  \, \D S(\vec{y})}_{J^2_{ij}(\vec{x})}.
\end{align}
The two integrals which we have called $J^1_{ij}$ and $J^2_{ij}$
can both be expressed in terms of the integral
\begin{equation}
  \label{eq:appA-J0}
  J^0_{ij}(\vec{x}) = \int_{\Gamma_0}
  T_{ij3}(\vec{x} - \vec{y}) \, \D S(\vec{y}),
\end{equation}
where $\Gamma_0 = \seti{\vec{x} \in \Real^3 : x_3 = 0}$.
The integrals are related through $J^0_{ij}(\vec{x}) =
J^1_{ij}(\vec{x} + a\vec{e}_3) = J^2_{ij}(\vec{x} - a\vec{e}_3)$,
with $\vec{e}_3 = (0,0,1)$. In fact, since $\Gamma_0$ is
infinite, the integral $J^0_{ij}(\vec{x})$ as given by
\eqref{eq:appA-J0} depends only on $x_3$, i.e.
\begin{equation}
  \label{eq:appA-J0-simple}
  J^0_{ij}(\vec{x}) = J^0_{ij}(x_3) = \int_{\Gamma_0}
  T_{ij3}(x_3 \vec{e}_3 - \vec{y}) \, \D S(\vec{y}).
\end{equation}
Inserting the expression for the stresslet $\mat{T}$ from
\eqref{eq:double-layer} into \eqref{eq:appA-J0-simple}, we find
\begin{equation}
  J^0_{ij}(x_3) = -6x_3 \int_{-\infty}^\infty \int_{-\infty}^\infty
  \frac{(x_3 \delta_{i3} - y_i)(x_3 \delta_{j3} - y_j)}
  {(y_1^2 + y_2^2 + x_3^2)^{5/2}} \D y_1 \D y_2,
\end{equation}
where $y_3 = 0$.
This double integral can be computed analytically, and the result
is
\begin{equation}
  \label{eq:appA-J0-result}
  J^0_{ij}(x_3) = -4\pi \sgn(x_3) \delta_{ij},
\end{equation}
where $\sgn(\cdot)$ denotes the sign function.
Using the relations $J^1_{ij}(\vec{x}) =
J^0_{ij}(x_3-a)$ and $J^2_{ij}(\vec{x}) = J^0_{ij}(x_3+a)$ and
inserting \eqref{eq:appA-J0-result} into
\eqref{eq:appA-double-layer-step}, we get
\begin{equation}
  \Dp_i[\tilde{\Gamma},\tilde{\vec{q}}](\vec{x}) =
  4\pi \tilde{q}_i \sgn(x_3-a) - 4\pi \tilde{q}_i \sgn(x_3+a).
\end{equation}
From this the result \eqref{eq:appA-stresslet-identity-1} follows.

\section{A pipe}
\label{app:stresslet-identity-pipe}

\begin{figure}[h!]
  \centering
  \begin{tikzpicture}
    \fill [myBlue!50] (-2.5,1.3) -- (2.5,1.2)
    arc [x radius=0.5, y radius=1.2, start angle=90, end angle=-90]
    -- (-2.5,-1.1)
    arc [x radius=0.5, y radius=1.2, start angle=-90, end angle=90];
    \draw [line width=1pt] (2.5,0) ellipse [x radius=0.5, y radius=1.2];
    \draw [line width=1pt] (-2.5,0.1) ellipse [x radius=0.5, y radius=1.2];
    \draw [line width=1pt] (2.5,1.2) -- (-2.5,1.3);
    \draw [line width=1pt] (2.5,-1.2) -- (-2.5,-1.1);
    \draw [dashed] (2.5,0) -- (-2.5,0.1);
    \draw [->] (2.5,0) -- (3.5,-0.02) node [right] {$x_1$};
    \draw [->] (0.5,1.24) -- (0.4812,0.3) node [right=1pt] {$\vec{n}$};
    \draw [->] (-0.5,-1.14) -- (-0.4812,-0.2) node [left=1pt] {$\vec{n}$};
    \fill [myBlue!20, draw=black, line width=1pt, line join=bevel,
           fill opacity=0.4]
    (-2.5,1.3) -- (2.5,1.2)
    arc [x radius=0.5, y radius=1.2, start angle=90, end angle=270]
    -- (-2.5,-1.1)
    arc [x radius=0.5, y radius=1.2, start angle=270, end angle=90];
    \node [anchor=west] at (2.9,1.2)
    {$\tilde{\Gamma}=\seti{\vec{x} \in \Real^3 : x_2^2 + x_3^2 = a^2}$};
    \node [anchor=west,myBlue] at (3.1,-0.6)
    {$\Omega=\seti{\vec{x} \in \Real^3 : x_2^2 + x_3^2 < a^2}$};
  \end{tikzpicture}
  \caption{An infinite cylindrical pipe $\tilde{\Gamma}$.}
  \label{fig:appA-pipe}
\end{figure}
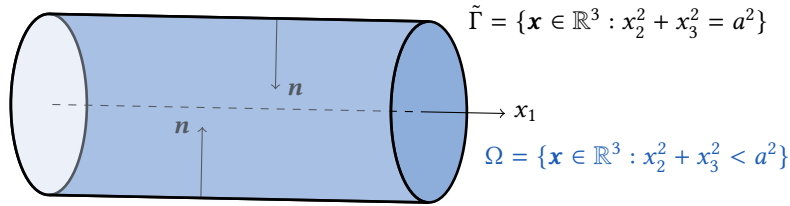

Let now $\tilde{\Gamma}$ be an infinitely long pipe given by the
equation $x_2^2 + x_3^2 = a^2$ for some $a>0$, as shown in
Figure~\ref{fig:appA-pipe}. Let the domain inside the pipe be
denoted by $\Omega$. We shall show that for any constant vector
$\tilde{\vec{q}} \in \Real^3$, the identity
\eqref{eq:appA-stresslet-identity-1} holds.

Let us introduce cylindrical coordinates and write $\vec{x} = x_1
\vec{e}_1 + r \vec{e}_\varphi$ and $\vec{y} = y_1 \vec{e}_1 + a
\vec{e}_\theta$ for the evaluation point and integration
variable, respectively. The unit vectors are given by
\begin{equation}
  \vec{e}_1 = (1,0,0), \qquad
  \vec{e}_\varphi = (0, \cos \varphi, \sin \varphi) \qquad\text{and}\qquad
  \vec{e}_\theta = (0, \cos \theta, \sin \theta),
\end{equation}
and $r \geq 0$.
Using the fact that the normal vector is given by
$\vec{n}(\vec{y}) = -\vec{e}_\theta$, we can write the double
layer potential from \eqref{eq:double-layer} as
\begin{equation}
  \label{eq:appA-pipe-first-step}
  \Dp_i[\tilde{\Gamma},\tilde{\vec{q}}](\vec{x}) =
  -\tilde{q}_j \int_0^{2\pi} \int_{-\infty}^\infty
  T_{ijk}(x_1\vec{e}_1 + r\vec{e}_\varphi
  -y_1 \vec{e}_1 - a \vec{e}_\theta) (\vec{e}_\theta)_k
  \, a \, \D y_1 \D \theta,
\end{equation}
where $(\vec{e}_\theta)_k$ denotes the $k$th component of $\vec{e}_\theta$.
Using the variable substitution $y_1 - x_1 = u$, we can eliminate
$x_1$. Writing out the stresslet $\mat{T}$ from \eqref{eq:double-layer},
and using a few trigonometric identities, the integral in
\eqref{eq:appA-pipe-first-step} can be written as
\begin{equation}
  \label{eq:appA-pipe-double-step}
  \Dp_i[\tilde{\Gamma},\tilde{\vec{q}}](\vec{x}) =
  6 \tilde{q}_j a \underbrace{\int_0^{2\pi} \int_{-\infty}^\infty
  \frac{
    (r\vec{e}_\varphi-u \vec{e}_1 - a \vec{e}_\theta)_i
    (r\vec{e}_\varphi-u \vec{e}_1 - a \vec{e}_\theta)_j
    (r \cos(\varphi - \theta) - a)
  }{\gps{\big}{u^2 + r^2 + a^2 - 2ra \cos (\varphi -
  \theta)}^{5/2}}
  \, \D u \, \D \theta}_{I_{ij}(r,\varphi)}.
\end{equation}
At this point it is not immediately apparent that the integral
which we have called $I_{ij}(r,\varphi)$ is independent of
$\varphi$, but that does indeed turn out to be the case. We
expect the offdiagonal elements of $I_{ij}$ to be zero, which can
be verified by first integrating in $u$ and then in $\theta$.
It thus remains to compute the diagonal elements of $I_{ij}$.

To compute $I_{11}$, first integrate in $u$ using the formula
\begin{equation}
  \int_{-\infty}^\infty \frac{u^2}{(u^2 + C)^{5/2}} \, \D u
  = \frac{2}{3C}, \qquad C > 0.
\end{equation}
The outer integral becomes
\begin{equation}
  I_{11}(r,\varphi) = \frac{2}{3}
  \int_0^{2\pi}
  \frac{r \cos(\varphi - \theta) - a}{r^2 + a^2 - 2ra \cos (\varphi - \theta)}
  \, \D\theta.
\end{equation}
The variable $\varphi$ can now be eliminated using the
substitution $\theta-\varphi = \nu$ (and the limits shifted back
to $[0, 2\pi]$ due to periodicity). The value of the integral can
then be calculated to be
\begin{equation}
  I_{11}(r, \varphi) = \frac{2\pi}{3a} \gps{\big}{\sgn(r-a) - 1},
  \qquad r \geq 0,
\end{equation}
where $\sgn(\cdot)$ is the sign function.

To compute $I_{22}$, first integrate in $u$ using the formula
\begin{equation}
  \label{eq:appA-known-integral-1}
  \int_{-\infty}^\infty \frac{1}{(u^2 + C)^{5/2}} \, \D u =
  \frac{4}{3C^2}, \qquad C > 0,
\end{equation}
to get
\begin{equation}
  I_{22}(r,\varphi) = \frac{4}{3}
  \int_0^{2\pi}
  \frac{(r\cos\varphi-a\cos\theta)^2
  \gp{r \cos(\varphi - \theta) - a}}
  {(r^2 + a^2 - 2ra \cos (\varphi - \theta))^2}
  \, \D\theta.
\end{equation}
Using the substitution $\theta - \varphi = \nu$ and shifting the
limits back to $[0, 2\pi]$ yields the integral
\begin{equation}
  \label{eq:appA-crucial-step-I22}
  I_{22}(r,\varphi) = \frac{4}{3}
  \int_0^{2\pi}
  \frac{(r\cos\varphi - a \cos \nu \cos \varphi + a \sin \nu \sin \varphi)^2
  \gp{r \cos \nu - a}}
  {(r^2 + a^2 - 2ra \cos \nu)^2}
  \, \D\nu,
\end{equation}
which we compute by expanding the square in the numerator, thus
splitting the integral into six terms, after which each term can
be integrated separately. The result is
\begin{equation}
  \label{eq:appA-I22-result}
  I_{22}(r, \varphi) = \frac{2\pi}{3a} \gps{\big}{\sgn(r-a) - 1},
  \qquad r \geq 0.
\end{equation}
Note that the dependence on $\varphi$ disappears when summing the
six terms to get the above result.

Finally, to compute $I_{33}$, we again start by integrating in $u$
using \eqref{eq:appA-known-integral-1}, after which we use the
substitution $\theta - \varphi = \nu$ to get
\begin{equation}
  \label{eq:appA-crucial-step-I33}
  I_{33}(r, \varphi) = \frac{4}{3} \int_0^{2\pi}
  \frac{(r\sin\varphi - a\cos \nu \sin \varphi - a\sin \nu \cos \varphi)^2 (r \cos \nu - a)}
  {(r^2 + a^2 - 2ra \cos \nu)^2}
  \, \D\nu.
\end{equation}
Comparing \eqref{eq:appA-crucial-step-I22} and
\eqref{eq:appA-crucial-step-I33}, note that $I_{33}(r, \varphi +
\pi/2) = I_{22}(r, \varphi)$. But as we saw in
\eqref{eq:appA-I22-result}, $I_{22}$ does not depend on
$\varphi$, so $I_{33} = I_{22}$.

To summarize, we have shown that
\begin{equation}
  I_{ij}(r, \varphi) = \frac{2\pi}{3a} \! \gps{\big}{\!\sgn(r-a) - 1}
  \delta_{ij}.
\end{equation}
Inserting this into \eqref{eq:appA-pipe-double-step}, we find
that
\begin{equation}
  \Dp_i[\tilde{\Gamma},\tilde{\vec{q}}](\vec{x}) = 4\pi
  \tilde{q}_i \! \gps{\big}{\!\sgn(r-a) - 1},
\end{equation}
from which the result \eqref{eq:appA-stresslet-identity-1}
follows for the pipe.

\FloatBarrier
\chapter{Construction of smooth rod particles}
\label{app:smooth-rod}

\begin{figure}[t!]
  \centering
  \includegraphics{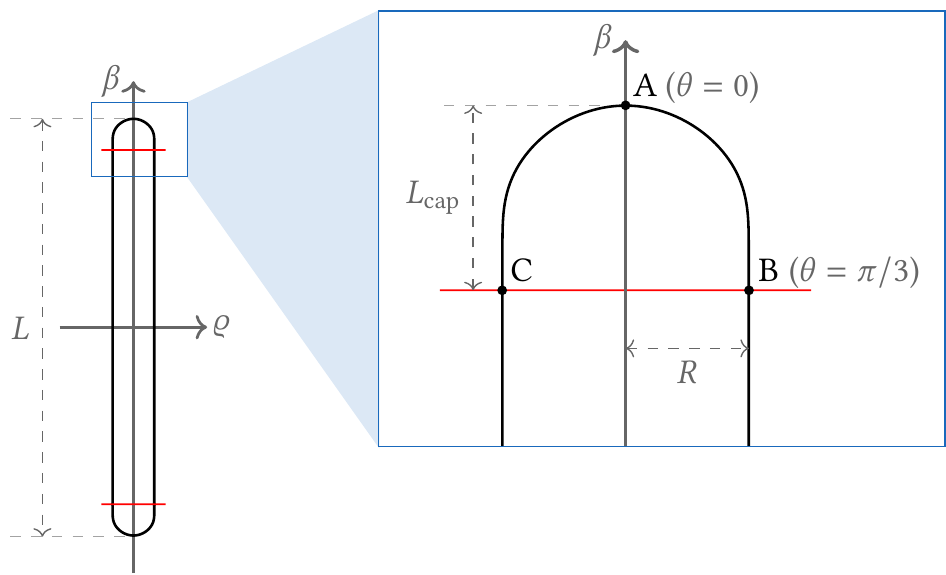}
  \caption{The shape of the smooth rod, here with $L=10$ and
  $R=0.5$.}
  \label{fig:smooth-rod-construction}
\end{figure}

In this section, we describe how the rod particles are
constructed to ensure that they are smooth everywhere. Recall
from section~\ref{sec:direct-quad-particles} the parametrization
\begin{equation}
  \label{eq:rod-parametrization-2}
  \begin{cases}
    x_1 = \varrho(\theta; L, R) \cos\varphi, \\
    x_2 = \varrho(\theta; L, R) \sin\varphi, \\
    x_3 = \beta(\theta; L, R),
  \end{cases}
\end{equation}
of the rod, where $\varphi \in [0,2\pi)$ and $\theta \in [0,
\pi]$ are parameters, $L$ is the length of the rod and $R$ the
radius. The goal here is to derive the shape functions
$\varrho(\cdot\,;L,R) : [0,\pi] \to
[0,R]$ and $\beta(\cdot\,;L,R) : [0,\pi] \to [-\tfrac{1}{2}L,
\tfrac{1}{2}L]$ so that the rod has the smooth shape shown in
Figure~\ref{fig:smooth-rod-construction}.
The rod consists of three smoothly joined parts: a top cap,
corresponding to $\theta \in I_1 = [0, \pi/3]$; a middle
cylinder, corresponding to $\theta \in I_2 = [\pi/3, 2\pi/3]$;
and a bottom cap, corresponding to $\theta \in I_3 = [2\pi/3,
\pi]$. Let the length of each cap be $L_\text{cap}$, as shown in
Figure~\ref{fig:smooth-rod-construction}. The ratio
$L_\text{cap}/R$ determines the aspect ratio of the cap. Here, we
fix this ratio by setting $L_\text{cap}$ to
\begin{equation}
  L_\text{cap} = 1.5 R,
\end{equation}
which gives the cap a shape similar to a half-sphere. The length
of the middle cylinder is then
\begin{equation}
  L_\text{mid} = L - 3R.
\end{equation}
However, note that the derivation below is valid for any value of
$L_\text{cap} \in (0, L/2)$, with $L_\text{mid} = L-2L_\text{cap}$.

Let us for fixed $L$ and $R$ define $\vec{g}(\theta) =
(g_1(\theta), g_2(\theta)) = (\varrho(\theta; L,R), \beta(\theta;
L,R))$. For the middle cylinder, the parametrization is
\begin{equation}
  g_1(\theta) = R, \qquad
  g_2(\theta) = \gp{1 - \frac{3}{\pi}\theta}\!L_\text{mid}
  + \frac{L_\text{mid}}{2}, \qquad
  \theta \in I_2 = [\pi/3, 2\pi/3].
\end{equation}
Note that $g_2$ is simply an affine function of $\theta$. At the
endpoints of the interval $I_2$ we have
\begin{equation}
  \label{eq:smooth-boundary-conditions-1}
  \begin{array}{r@{\:}c@{\:}l}
    \vec{g}(\pi/3) &=& (R, L_\text{mid}/2), \\[4pt]
    \vec{g}'(\pi/3) &=& (0, -(3/\pi)L_\text{mid}), \\[4pt]
    \vec{g}^{(n)}(\pi/3) &=& (0,0), \quad n \geq 2,
  \end{array}
  \qquad\qquad\text{and}\qquad\qquad
  \begin{array}{r@{\:}c@{\:}l}
    \vec{g}(2\pi/3) &=& (R, -L_\text{mid}/2), \\[4pt]
    \vec{g}'(2\pi/3) &=& (0, -(3/\pi)L_\text{mid}), \\[4pt]
    \vec{g}^{(n)}(2\pi/3) &=& (0,0), \quad n \geq 2.
  \end{array}
\end{equation}
Our goal is now to extend the parametrization $\vec{g}(\theta)$
to $I_1$ and $I_3$ in a way such that the unit tangent vector
$\vec{g}'(\theta)/\absi{\vec{g}'(\theta)}$ and its higher
derivatives are continuous everywhere. As an intermediate step we
introduce an auxiliary function $\widehat{\vec{g}}(t) =
(\widehat{g}_1(t), \widehat{g}_2(t))$ with a different parameter
$t \in [-1,1]$. The function $\widehat{\vec{g}}$ should trace the
curve from C to B via A in
Figure~\ref{fig:smooth-rod-construction}, with C corresponding to
$t=-1$, A corresponding to $t=0$ and B corresponding to $t=1$. We
will later relate $t \in [0,1]$ to $\theta \in [0,\pi/3]$ to get
the final parametrization. At this point, note that to match
\eqref{eq:smooth-boundary-conditions-1} we must require
\begin{equation}
  \label{eq:smooth-boundary-conditions-2}
  \begin{array}{r@{\:}c@{\:}l}
    \widehat{\vec{g}}(1) &=& (R, L_\text{mid}/2), \\[4pt]
    \widehat{\vec{g}}'(1) &=& (0, -b), \\[4pt]
    \widehat{\vec{g}}^{(n)}(1) &=& (0,0), \quad n \geq 2,
  \end{array}
  \qquad\qquad\text{and}\qquad\qquad
  \begin{array}{r@{\:}c@{\:}l}
    \widehat{\vec{g}}(-1) &=& (-R, L_\text{mid}/2), \\[4pt]
    \widehat{\vec{g}}'(-1) &=& (0, b), \\[4pt]
    \widehat{\vec{g}}^{(n)}(-1) &=& (0,0), \quad n \geq 2,
  \end{array}
\end{equation}
where $b$ is some positive constant. In order to construct
$\widehat{\vec{g}}(t)$ we will use a bump function $\psi : \Real
\to \Real$, which must satisfy the following:
\begin{itemize}
  \item
    $\psi$ must be infinitely differentiable on $\Real$,
  \item
    $\psi$ must have compact support in $[-1,1]$, i.e.\ $\psi(t) =
    0$ if $t > 1$ or $t < -1$,
  \item
    $\psi(t)$ must be positive for $t \in (-1,1)$,
  \item
    $\psi$ must be even, i.e.\ $\psi(t) = \psi(-t)$ for all $t
    \in \Real$.
\end{itemize}
We also introduce its primitive function
\begin{equation}
  \Psi(t) = \int_0^t \psi(\tau) \, \D \tau, \qquad t \in \Real,
\end{equation}
which is an odd function since $\psi$ is even. We choose a
specific bump function, namely%
\footnote{This function was found at \url{https://math.stackexchange.com/a/101484}.}
\begin{equation}
  \psi(t) = \begin{cases}
    \dfrac{(t^2+1) \exp[4t/(t^2-1)]}{\gbs{\big}{(t^2-1)
    (1+\exp[4t/(t^2-1)])}^2}, & \text{if $t \in (-1,1)$}, \\[12pt]
    0, & \text{otherwise}.
  \end{cases}
\end{equation}
This function has the primitive function
\begin{equation}
  \Psi(t) = \begin{cases}
    -\dfrac{1}{8} \tanh\fp{-\dfrac{2t}{1-t^2}}, & \text{if $t \in (-1,1)$}, \\[12pt]
    -\dfrac{1}{8}, & \text{if $t \leq -1$}, \\[12pt]
    \dfrac{1}{8}, & \text{if $t \geq 1$}.
  \end{cases}
\end{equation}
We now construct $\widehat{\vec{g}}(t)$ as
\begin{equation}
  \label{eq:auxiliary-parametrization}
  \widehat{g}_1(t) = R \frac{\Psi(t)}{\Psi(1)}, \qquad
  \widehat{g}_2(t) = \frac{L_\text{mid}}{2} - b \int_{-1}^t
  \frac{\Psi(\tau)}{\Psi(1)} \, \D \tau, \qquad
  t \in [-1,1],
\end{equation}
which satisfies \eqref{eq:smooth-boundary-conditions-2}. We can
determine $b$ by noting that we must have $\widehat{g}_2(0) =
L/2$ (at point A in Figure~\ref{fig:smooth-rod-construction}),
which yields
\begin{equation}
  \label{eq:determine-b}
  b = L_\text{cap} \frac{\Psi(1)}{\int_0^1 \Psi(\tau) \, \D \tau}.
\end{equation}
The integrals of $\Psi$ in \eqref{eq:auxiliary-parametrization}
and \eqref{eq:determine-b} are computed numerically using
MATLAB's \textit{integral} function.

Finally, we go from the parameter $t$ to the parameter
$\theta$. We would like the discretization points to be
distributed as Gauss--Legendre points in the arclength, and so we
must choose $\theta$ so that it is proportional to the arclength
on the caps. Consider the arclength
\begin{equation}
  s(t) = \int_0^t \absi{\widehat{\vec{g}}'(\tau)} \, \D \tau,
  \qquad t \in [0,1].
\end{equation}
Let us then define
\begin{equation}
  \theta = G(t) = \frac{\pi}{3} \frac{s(t)}{s(1)}, \qquad t \in [0, 1],
\end{equation}
and note that this defines $\theta \in I_1 = [0,\pi/3]$ as an
invertible function of $t \in [0,1]$. We can now define
$\vec{g}(\theta) = \vec{g}(G(t)) = \widehat{\vec{g}}(t)$ for $t
\in [0,1]$, and thus
\begin{equation}
  \vec{g}(\theta) = \widehat{\vec{g}}(G^{-1}(\theta)), \qquad
  \theta \in I_1 = [0,\pi/3].
\end{equation}
The bottom cap should be the reflection of the top cap in the
plane corresponding to $\beta=0$, so
\begin{equation}
  \vec{g}(\theta) = (g_1(\pi-\theta), -g_2(\pi-\theta)),
  \qquad \theta \in I_3 = [2\pi/3, \pi].
\end{equation}
Now that we have defined $\vec{g}(\theta)$ for all $\theta \in
[0,\pi]$, its two components $g_1$ and $g_2$ correspond to the
shape factors $\varrho(\theta;L,R)$ and $\beta(\theta;L,R)$,
respectively, which are to be used in
\eqref{eq:rod-parametrization-2}.

\FloatBarrier
\chapter{Derivation of the safety factor $\gamma$}
\label{app:safety-factor}

\begin{figure}[h!]
  \centering\small
  \begin{minipage}[b]{0.6\textwidth}
    \centering
    \includegraphics{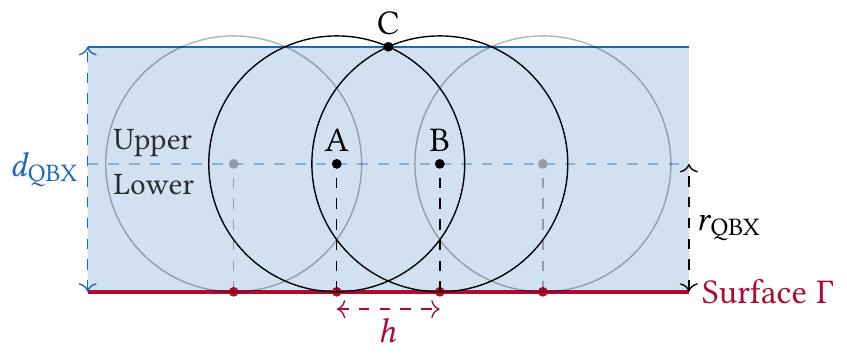}\\
    (a)\hspace*{7mm}
  \end{minipage}%
  \begin{minipage}[b]{0.4\textwidth}
    \centering
    \includegraphics{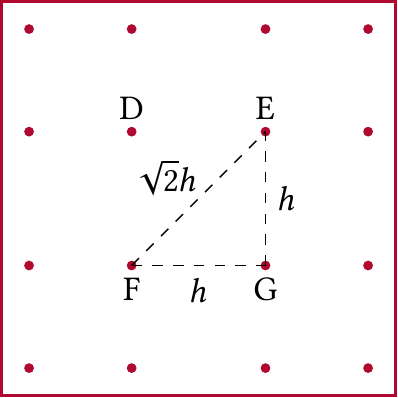}\\
    (b)
  \end{minipage}%
  \caption{(a) Balls of convergence for a flat surface~$\Gamma$
  (c.f.\ Figure~\ref{fig:param1-qbx-ball-sketch}), seen from the
  side.\quad (b) Grid points of $\Gamma$, seen from above. Here,
  $h$ is the largest spacing between grid points in each
  tensorial direction.}
  \label{fig:appC-balls}
\end{figure}

Recall from section~\ref{sec:parameters} that one may want to
select $d_\text{QBX}$ larger than $r_\text{QBX}$ since QBX may be
faster than the upsampled quadrature due to the precomputation
scheme. Let us call the set of points of the QBX region with
distance to $\Gamma$ greater than $r_\text{QBX}$ the \emph{upper}
QBX region, and the set of points with distance to $\Gamma$ smaller than
$r_\text{QBX}$ the \emph{lower} QBX region, as shown in
Figure~\ref{fig:appC-balls}~(a). As noted in
section~\ref{sec:param1-qbx}, putting $d_\text{QBX} = 2
r_\text{QBX}$ would lead to some areas of the upper QBX region
not falling within any ball of convergence. To avoid this, we
introduce a safety factor $\gamma$ and require that
\begin{equation}
  \label{eq:appC-starting-point}
  d_\text{QBX} \leq 2\gamma r_\text{QBX}.
\end{equation}
The goal here is to derive the value of the safety
factor~$\gamma$. We assume for simplicity that $\Gamma$ is a flat
surface.

The key is to choose $d_\text{QBX}$ below the intersection of
neighbouring balls of convergence, marked by the point C in
Figure~\ref{fig:appC-balls}~(a). Since the grid on $\Gamma$ is
two-dimensional, the largest distance between neighbouring grid
points is not $h$ but $\sqrt{2}h$, where $h$ is as shown in
Figure~\ref{fig:appC-balls}~(b). The four balls of convergence of
the expansion centres above the grid points D--G in this figure
intersect at distance
\begin{equation}
  \label{eq:appC-restriction-1}
  d_{\star}(h) = r_\text{QBX} + \sqrt{r_\text{QBX}^2 -
  \bigg(\frac{\sqrt{2}h}{2}\bigg)^2}
\end{equation}
from $\Gamma$. Thus, choosing $d_\text{QBX} \leq d_{\star}(h)$ is
sufficient to ensure that all points in the upper QBX region
fall within a ball of convergence.
This restriction on $d_\text{QBX}$ can be simplified by minimizing
$d_{\star}(h)$ with respect to $h$, subject to the constraint $0 < h
\leq r_\text{QBX}$. The result is
\begin{equation}
  \label{eq:appC-restriction-2}
  d_{\star\star} = \min_{0 < h \leq r_\text{QBX}} d_{\star}(h) =
  d_{\star}(r_\text{QBX}) = \gp{1 + \frac{1}{\sqrt{2}}} r_\text{QBX}.
\end{equation}
It is thus sufficient to require that $d_\text{QBX} \leq
d_{\star\star}$. Comparing \eqref{eq:appC-restriction-2} and
\eqref{eq:appC-starting-point}, we see that the safety factor
should be
\begin{equation}
  \gamma = \frac{1}{2} \gp{1 + \frac{1}{\sqrt{2}}} \approx 0.85.
\end{equation}

This derivation holds when $\Gamma$ is a flat surface, in which
case the requirement \eqref{eq:appC-starting-point} with $\gamma
= 0.85$ guarantees that all points in the upper QBX region fall
within a ball of convergence, as long as $h \leq r_\text{QBX}$.
If $\Gamma$ is curved, this guarantee holds on the concave
side of $\Gamma$, but not necessarily on the convex side, where
$d_\text{QBX}$ may have to be even smaller for the guarantee to
hold. Nonetheless, we use \eqref{eq:appC-starting-point} with
$\gamma=0.85$ also for convex surfaces such as rods and
spheroids, and it seems to work well in practice. Of course, the
parameter selection strategy
(section~\ref{sec:param-general-summary}, step~2) will in most
cases choose $d_\text{QBX}$ less than the upper bound $2\gamma
r_\text{QBX}$.

\FloatBarrier
\chapter{Efficient computation of streamlines in periodic flow}
\label{app:streamlines}

To compute streamlines in a periodic problem such as in
section~\ref{sec:res3-packed-rods}, we must first solve the
periodic boundary integral equation as described in
section~\ref{sec:periodicity} to get the density $\vec{q}$ on
$\Gamma$. We can then compute the flow field
\begin{equation}
  \label{eq:appD-flowfield}
  \vec{u}(\vec{x}_\text{e}) = \vec{u}_\text{bg}(\vec{x}_\text{e}) +
  \vec{\Dp}^\text{3P}[\Gamma, \vec{q}](\vec{x}_\text{e}) +
  \sum_{\alpha=1}^M \vec{\Vp}^{(\alpha),\text{3P}}[\vec{F}^{(\alpha)},
  \vec{\tau}^{(\alpha)}](\vec{x}_\text{e})
\end{equation}
at any evaluation point $\vec{x}_\text{e}$ in the fluid domain. To compute a
streamline we pick any point $\vec{x}_0$ in the fluid domain and
then solve the differential equation
\begin{equation}
  \label{eq:streamline-ode}
  \od{\vec{x}_\text{e}}{t} = \vec{u}(\vec{x}_\text{e}(t)), \qquad
  \vec{x}_\text{e}(0) = \vec{x}_0.
\end{equation}
Of course, \eqref{eq:streamline-ode} is discretized using some
timestepping method, which must evaluate
\eqref{eq:appD-flowfield} at every timestep. Recall that the
periodic double layer potential $\vec{\Dp}^\text{3P}$ is split
into two parts
\begin{equation}
  \vec{\Dp}^\text{3P}[\Gamma,\vec{q}](\vec{x}_\text{e}) =
  \vec{\Dp}^\text{3P,R}[\Gamma,\vec{q}](\vec{x}_\text{e};\xi)
  + \vec{\Dp}^\text{3P,F}[\Gamma,\vec{q}](\vec{x}_\text{e};\xi)
\end{equation}
and similarly for $\vec{\Vp}^{(\alpha),\text{3P}}$. The first part
$\vec{\Dp}^\text{3P,R}$ decays fast and is treated according to
section~\ref{sec:periodicity}. The second part
$\vec{\Dp}^\text{3P,F}$ decays slowly in real space, but since it
is smooth its Fourier coefficients decay fast. In the Spectral
Ewald method, $\vec{\Dp}^\text{3P,F}$ as given by
\eqref{eq:fourier-sum} is first discretized using the direct
quadrature rule \eqref{eq:quadrature-rule} to give
\begin{equation}
  \label{eq:discretized-fourier-sum}
  \Dp_i^\text{3P,F,$h$}[\Gamma, \vec{q}](\vec{x}_\text{e}; \xi) =
  \sum_{\vec{k} \in \Integer^3} \sum_{s=1}^N
  T_{ijl}^\text{F}(\vec{x}_\text{e}+\vec{k}\cdot\vec{B}-\vec{x}_s; \xi)
  q_j(\vec{x}_s) n_l(\vec{x}_s) w_s.
\end{equation}
This is a periodic sum of point sources with strengths
$Z_{jl}(\vec{x}_s) = q_j(\vec{x}_s) n_l(\vec{x}_s) w_s$. The
Spectral Ewald method \cite{klinteberg14,lindbo10,lindbo11}
computes the periodic sum \eqref{eq:discretized-fourier-sum} in
five steps:
\begin{enumerate}
  \item
    \textbf{Spreading point sources to a grid:} A
    three-dimensional uniform grid is constructed over the
    primary cell. A window function $W(\vec{r})$ is convolved
    with the point sources in the primary cell to give
    \begin{equation}
      \label{eq:appD-convolution1}
      H_{jl}(\vec{x}) = \sum_{s=1}^N Z_{jl}(\vec{x}_s)
      W([\vec{x}-\vec{x}_s]_*).
    \end{equation}
    Here, $[\cdot]_*$ denotes that the shortest periodic distance
    should be used, i.e.
    \begin{equation}
      [\vec{r}]_* = \vec{r} + \vec{B} \cdot \argmin_{\vec{k} \in
      \Integer^3} \absi{\vec{r}+\vec{B}\cdot\vec{k}},
    \end{equation}
    where $\vec{B} = (B_1, B_2, B_3)$ is the size of the periodic cell.
    In this work the window function is a truncated Gaussian,
    given by $W(\vec{r}) = w(r_1)w(r_2)w(r_3)$,
    where
    \begin{equation}
      w(r) = \begin{cases}
        \E^{-A (r/r_\text{trunc})^2}, & \text{if $\absi{r} \leq r_\text{trunc} =
        h_\text{g} P/2$}, \\
        0, & \text{otherwise}.
      \end{cases}
    \end{equation}
    Here, $h_\text{g}$ is the grid spacing of the uniform grid,
    $P$ is the number of grid points within the support of $w$,
    and $A = 0.9^2 \pi P/2$. The parameter $P$ is chosen as
    discussed in \cite{klinteberg14}.
    It is also possible to use other window functions than the
    Gaussian, as discussed for example in \cite{shamshirgar17}.

    The function $H_{jl}(\vec{x})$ as given by
    \eqref{eq:appD-convolution1} is evaluated on the uniform grid.

  \item
    \textbf{FFT:} The three-dimensional Fourier transform
    $\widehat{H}_{jl}(\vec{k})$ is computed using the FFT. This
    is possible since $H_{jl}(\vec{x})$ is defined on a uniform grid.

  \item
    \textbf{Scaling:} The result is multiplied by the Fourier
    transform of $\mat{T}^\text{F}$, and divided by the Fourier
    transform of the window function~$W$ to undo the convolution in
    step~1. Since we will convolve again in step~5, this division
    is done twice. Thus, we here compute
    \begin{equation}
      \widehat{\widetilde{H}}_i(\vec{k}) =
      \widehat{T}^\text{F}_{ijl}(\vec{k};\xi)
      \frac{1}{[\widehat{W}(\vec{k})]^2}
      \widehat{H}_{jl}(\vec{k}),
    \end{equation}
    where
    \begin{equation}
      \widehat{T}^\text{F}_{ijl}(\vec{k};\xi) = \sqrt{-1}
      \frac{\pi}{\absi{\vec{k}}^2} \fbs{\bigg}{(\delta_{ij}k_{l} +
      \delta_{jl}k_{i} + \delta_{li}k_{j})
      - 2\frac{k_{i}k_{j}k_{l}}{\absi{\vec{k}}^2}}
      \gp{8 + 2 \frac{\absi{\vec{k}}^2}{\xi^2} + \frac{\absi{\vec{k}}^4}{\xi^4}}
      \E^{-\absi{\vec{k}}^2 / (4\xi^2)},
    \end{equation}
    as given in \cite{klinteberg14}.

  \item
    \textbf{IFFT:} An inverse FFT is applied to
    $\smash{\widehat{\widetilde{H}}_i(\vec{k})}$ to compute
    $\widetilde{H}_i(\vec{x})$ on the uniform grid.

  \item
    \textbf{Gathering:} In order to compute the final result at
    the evaluation point $\vec{x}_\text{e}$ (which need not be on
    the uniform grid), another convolution with the window
    function is performed, i.e.\
    \begin{equation}
      \label{eq:gather-convolution}
      \Dp_i^\text{3P,F,$h$}[\Gamma, \vec{q}](\vec{x}_\text{e}; \xi) =
      \int_B \widetilde{H}_i(\vec{x})
      W([\vec{x}_\text{e}-\vec{x}]_*) \, \D \vec{x},
    \end{equation}
    where $B$ denotes the primary cell. The integral in
    \eqref{eq:gather-convolution} is evaluated using the
    trapezoidal rule on the uniform grid, which is spectrally
    accurate since the integrand is periodic.
\end{enumerate}
Since the density $\vec{q}$ does not change during the
computation of the streamlines, and the evaluation point
$\vec{x}_\text{e}$ enters only in step~5 above, it is possible to
do step 1--4 once before starting to compute the
streamlines, and save $\widetilde{H}_i(\vec{x})$ on the
uniform grid from step~4. When the Fourier-space part
$\vec{\Dp}^\text{3P,F}[\Gamma,\vec{q}](\vec{x}_\text{e};\xi)$
is to be evaluated at $\vec{x}_\text{e}(t)$ at every timestep of
solving \eqref{eq:streamline-ode}, it is then enough to do only
step~5. This speeds up the computation of the streamlines since
evaluating \eqref{eq:gather-convolution} is fast for a single
evaluation point. The
real-space part $\vec{\Dp}^\text{3P,R}[\Gamma,\vec{q}](\vec{x}_\text{e};\xi)$
must be computed from scratch at every timestep, but this is fast
since it is a local sum due to its rapid decay.

The periodic completion flow $\vec{\Vp}^{(\alpha),\text{3P}}$ which appears
in \eqref{eq:appD-flowfield} is treated in a very similar way;
for details, we refer to
\cite{lindbo10,klinteberg14,klinteberg16a}. Note that steps 1--5
of the Spectral Ewald method are also what is used when solving
the periodic boundary integral equation as described in
section~\ref{sec:periodicity}, but in that situation all the
evaluation points (i.e.\ the grid points of $\Gamma$) are known
in advance so they can all be fed into step~5 at the same time.


\FloatBarrier
\chapter*{References}

\renewcommand*{\bibsection}{}

\end{document}

%% file: custom-colors.tex
\definecolor{myYellow}{rgb:HTML}{EDDF20}
\definecolor{myBlue}{rgb:HTML}{1B5CB8}
\definecolor{myRed}{rgb:HTML}{C01212}
\definecolor{myGreen}{rgb:HTML}{70C61B}
\definecolor{myOrange}{rgb:HTML}{EC901C}
\definecolor{myCyan}{rgb:HTML}{30D6E0}
\definecolor{myPurple}{rgb:HTML}{832F8E}

\definecolor{KTHblue}{RGB}{25,105,188}
\definecolor{KTHlblue}{RGB}{22,159,219}
\definecolor{KTHyellow}{RGB}{251,186,0}
\definecolor{KTHred}{RGB}{176,9,48}
\definecolor{KTHlred}{RGB}{231,51,57}
\definecolor{KTHgreen}{RGB}{98,146,46}
\definecolor{KTHlgreen}{RGB}{175,202,11}
\definecolor{KTHpink}{RGB}{219,81,151}

\colorlet{fblue}{KTHblue}
\colorlet{fred}{KTHlred!80!KTHyellow}
\colorlet{fyellow}{KTHyellow}
\colorlet{fpurple}{KTHblue!50!KTHlred}
\colorlet{fgreen}{KTHgreen}
\colorlet{fcyan}{KTHlblue}
\colorlet{fmagenta}{KTHpink}

%% file: main.bbl
\begin{thebibliography}{10}\small\setlength{\itemsep}{0pt plus 0.3ex}
  \bibitem{atkinson97}
    \bibauthor{K.~E.~Atkinson}, \bibtitle{The Numerical Solution
    of Integral Equations of the Second Kind},
    Cambridge University Press, Cambridge, 1997,
    \bibdoi{10.1017/CBO9780511626340}.
  \bibitem{bagge17}
    \bibauthor{J.~Bagge, A.-K.~Tornberg},
    \bibarticle{Accurate quadrature methods with application to
    Stokes flow with particles in confined geometries}, in
    D.~J.~Chappell (ed.), \bibtitle{Proceedings of the Eleventh UK Conference on
    Boundary Integral Methods (UKBIM 11)}, Nottingham, UK, July 2017, pp.~15--24,
    ISBN~9780993111297. Available:
    \textcolor{blue}{\url{http://irep.ntu.ac.uk/id/eprint/31463}}
  \bibitem{barnett14}
    \bibauthor{A.~H.~Barnett},
    \bibarticle{Evaluation of layer potentials close to the
    boundary for Laplace and Helmholtz problems on analytic
    planar domains},
    \bibtitle{SIAM J.\ Sci.\ Comput.} \bibvolume{36}
    \bibissue{2}, A427--A451 (2014),
    \bibdoi{10.1137/120900253}.
  \bibitem{barnett15}
    \bibauthor{A.~Barnett, B.~Wu, S.~Veerapaneni},
    \bibarticle{Spectrally accurate quadratures for evaluation of
    layer potentials close to the boundary for the 2D Stokes and
    Laplace equations},
    \bibtitle{SIAM J.\ Sci.\ Comput.} \bibvolume{37}
    \bibissue{4}, B519--B542 (2015),
    \bibdoi{10.1137/140990826}.
  \bibitem{beale01}
    \bibauthor{J.~T.~Beale, M.-C.~Lai},
    \bibarticle{A method for computing nearly singular integrals},
    \bibtitle{SIAM J.\ Numer.\ Anal.} \bibvolume{38}
    \bibissue{6}, 1902--1925 (2001),
    \bibdoi{10.1137/S0036142999362845}.
  \bibitem{beale16}
    \bibauthor{J.~T.~Beale, W.~Ying, J.~R.~Wilson},
    \bibarticle{A Simple Method for Computing Singular or Nearly
    Singular Integrals on Closed Surfaces},
    \bibtitle{Commun.\ Comput.\ Phys.}
    \bibvolume{20} \bibissue{3}, 733--753 (2016),
    \bibdoi{10.4208/cicp.030815.240216a}.
  \bibitem{berrut04}
    \bibauthor{J.-P.~Berrut, L.~N.~Trefethen},
    \bibarticle{Barycentric Lagrange Interpolation},
    \bibtitle{SIAM Rev.} \bibvolume{46} \bibissue{3}, 501--517 (2004),
    \bibdoi{10.1137/S0036144502417715}.
  \bibitem{blake71}
    \bibauthor{J.~R.~Blake},
    \bibarticle{A note on the image system for a stokeslet in a
    no-slip boundary},
    \bibtitle{Proc.\ Camb.\ Phil.\ Soc.} \bibvolume{70} \bibissue{2}, 303--310
    (1971), \bibdoi{10.1017/S0305004100049902}.
  \bibitem{bruno01}
    \bibauthor{O.~P.~Bruno, L.~A.~Kunyansky},
    \bibarticle{A Fast, High-Order Algorithm for the Solution of
    Surface Scattering Problems: Basic Implementation, Tests, and
    Applications},
    \bibtitle{J.\ Comp.\ Phys.}
    \bibvolume{169} \bibissue{1}, 80--110 (2001),
    \bibdoi{10.1006/jcph.2001.6714}.
  \bibitem{carvalho18b}
    \bibauthor{C.~Carvalho, S.~Khatri, A.~D.~Kim},
    \bibarticle{Close evaluation of layer potentials in three dimensions},
    \href{https://arxiv.org/abs/1807.02474}{\bibtitle{\textcolor{blue}{arXiv:1807.02474 [math.NA]}}}
    (2018).
  \bibitem{carvalho18c}
    \bibauthor{C.~Carvalho, S.~Khatri, A.~D.~Kim},
    \bibarticle{Asymptotic approximations for the close
    evaluation of double-layer potentials},
    \href{https://arxiv.org/abs/1810.02483}{\bibtitle{\textcolor{blue}{arXiv:1810.02483 [math.NA]}}}
    (2018).
  \bibitem{corona18}
    \bibauthor{E.~Corona, S.~Veerapaneni},
    \bibarticle{Boundary integral equation analysis for
    suspension of spheres in Stokes flow},
    \bibtitle{J.\ Comp.\ Phys.}
    \bibvolume{362}, 327--345 (2018),
    \bibdoi{10.1016/j.jcp.2018.02.017}.
  \bibitem{deserno98}
    \bibauthor{M.~Deserno, C.~Holm},
    \bibarticle{How to mesh up Ewald sums. I. A theoretical and
    numerical comparison of various particle mesh routines},
    \bibtitle{J.~Chem.~Phys.} \bibvolume{109} \bibissue{18},
    7678--7693 (1998), \bibdoi{10.1063/1.477414}.
  \bibitem{epstein13}
    \bibauthor{C.~L.~Epstein, L.~Greengard, A.~Klöckner},
    \bibarticle{On the convergence of local expansions of layer
    potentials}, \bibtitle{SIAM J.\ Numer.\ Anal.}
    \bibvolume{51} \bibissue{5}, 2660--2679 (2013),
    \bibdoi{10.1137/120902859}.
  \bibitem{gimbutas15}
    \bibauthor{Z.~Gimbutas, L.~Greengard, S.~Veerapaneni},
    \bibarticle{Simple and efficient representations for the
    fundamental solutions of Stokes flow in a half-space},
    \bibtitle{J.\ Fluid Mech.}
    \bibvolume{776}, R1 (2015),
    \bibdoi{10.1017/jfm.2015.302}.
  \bibitem{gontijo15}
    \bibauthor{R.~G.~Gontijo, F.~R.~Cunha},
    \bibarticle{Dynamic numerical simulations of magnetically
    interacting suspensions in creeping flow},
    \bibtitle{Powder Technol.} \bibvolume{279}, 146--165 (2015),
    \bibdoi{10.1016/j.powtec.2015.03.033}.
  \bibitem{greengard87}
    \bibauthor{L.~Greengard, V.~Rokhlin},
    \bibarticle{A fast algorithm for particle simulations},
    \bibtitle{J.\ Comp.\ Phys.}
    \bibvolume{73} \bibissue{2}, 325--348 (1987),
    \bibdoi{10.1016/0021-9991(87)90140-9}.
  \bibitem{greengard97}
    \bibauthor{L.~Greengard, V.~Rokhlin},
    \bibarticle{A new version of the Fast Multipole Method for
    the Laplace equation in three dimensions},
    \bibtitle{Acta Numer.}
    \bibvolume{6}, 229--269 (1997),
    \bibdoi{10.1017/S0962492900002725}.
  \bibitem{guasto10}
    \bibauthor{J.~S.~Guasto, A.~S.~Ross, J.~P.~Gollub},
    \bibarticle{Hydrodynamic irreversibility in particle
    suspensions with nonuniform strain},
    \bibtitle{Phys.\ Rev.\ E} \bibvolume{81} \bibissue{6}, 061401 (2010),
    \bibdoi{10.1103/PhysRevE.81.061401}.
  \bibitem{guasto12}
    \bibauthor{J.~S.~Guasto, R.~Rusconi, R.~Stocker},
    \bibarticle{Fluid Mechanics of Planktonic Microorganisms},
    \bibtitle{Annu.\ Rev.\ Fluid Mech.}
    \bibvolume{44}, 373--400 (2012),
    \bibdoi{10.1146/annurev-fluid-120710-101156}.
  \bibitem{helsing08}
    \bibauthor{J.~Helsing, R.~Ojala},
    \bibarticle{On the evaluation of layer potentials close to
    their sources},
    \bibtitle{J.\ Comp.\ Phys.}
    \bibvolume{227} \bibissue{5}, 2899--2921 (2008),
    \bibdoi{10.1016/j.jcp.2007.11.024}.
  \bibitem{kimkarrila91}
    \bibauthor{S.~Kim, S.~J.~Karrila},
    \bibtitle{Microhydrodynamics : Principles and Selected
    Applications}, Butterworth--Heinemann, Boston, 1991,
    \bibdoi{10.1016/C2013-0-04644-0}.
  \bibitem{klinteberg14}
    \bibauthor{L.~af~Klinteberg, A.-K.~Tornberg},
    \bibarticle{Fast Ewald summation for Stokesian particle suspensions},
    \bibtitle{Int.\ J.\ Numer.\ Methods Fluids} \bibvolume{76}
    \bibissue{10}, 669--698 (2014),
    \bibdoi{10.1002/fld.3953}.
  \bibitem{klinteberg16a}
    \bibauthor{L.~af~Klinteberg},
    \bibarticle{Ewald summation for the rotlet singularity of
    Stokes flow},
    \href{https://arxiv.org/abs/1603.07467}{\bibtitle{\textcolor{blue}{arXiv:1603.07467 [physics.flu-dyn]}}}
    (2016).
  \bibitem{klinteberg16b}
    \bibauthor{L.~af~Klinteberg, A.-K.~Tornberg},
    \bibarticle{A fast integral equation method for solid
    particles in viscous flow using quadrature by expansion},
    \bibtitle{J.\ Comp.\ Phys.} \bibvolume{326}, 420--445 (2016),
    \bibdoi{10.1016/j.jcp.2016.09.006}.
  \bibitem{klinteberg17}
    \bibauthor{L.~af~Klinteberg, A.-K.~Tornberg},
    \bibarticle{Error estimation for quadrature by expansion in
    layer potential evaluation},
    \bibtitle{Adv.\ Comput.\ Math.} \bibvolume{43}, 195--234 (2017),
    \bibdoi{10.1007/s10444-016-9484-x}.
  \bibitem{klinteberg17b}
    \bibauthor{L.~af~Klinteberg, D.~Saffar~Shamshirgar, A.-K.~Tornberg},
    \bibarticle{Fast Ewald summation for free-space Stokes
    potentials},
    \bibtitle{Res.\ Math.\ Sci.} \bibvolume{4} \bibissue{1}, (2017),
    \bibdoi{10.1186/s40687-016-0092-7}.
  \bibitem{klinteberg18}
    \bibauthor{L.~af~Klinteberg, A.-K.~Tornberg},
    \bibarticle{Adaptive quadrature by expansion for layer
    potential evaluation in two dimensions},
    \bibtitle{SIAM J.\ Sci.\ Comput.} \bibvolume{40}
    \bibissue{3}, A1225--A1249 (2018),
    \bibdoi{10.1137/17M1121615}.
  \bibitem{klockner13}
    \bibauthor{A.~Klöckner, A.~Barnett, L.~Greengard, M.~O'Neil},
    \bibarticle{Quadrature by expansion: A new method for the
    evaluation of layer potentials}, \bibtitle{J.\ Comp.\ Phys.}
    \bibvolume{252}, 332--349 (2013),
    \bibdoi{10.1016/j.jcp.2013.06.027}.
  \bibitem{kress14}
    \bibauthor{R.~Kress}, \bibtitle{Linear Integral Equations},
    Springer, New York, 3rd ed., 2014,
    \bibdoi{10.1007/978-1-4614-9593-2}.
  \bibitem{lindbo10}
    \bibauthor{D.~Lindbo, A.-K.~Tornberg},
    \bibarticle{Spectrally accurate fast summation for periodic
    Stokes potentials}, \bibtitle{J.\ Comp.\ Phys.}
    \bibvolume{229} \bibissue{23}, 8994--9010 (2010),
    \bibdoi{10.1016/j.jcp.2010.08.026}.
  \bibitem{lindbo11}
    \bibauthor{D.~Lindbo, A.-K.~Tornberg},
    \bibarticle{Spectral accuracy in fast Ewald-based methods for
    particle simulations}, \bibtitle{J.\ Comp.\ Phys.}
    \bibvolume{230} \bibissue{24}, 8744--8761 (2011),
    \bibdoi{10.1016/j.jcp.2011.08.022}.
  \bibitem{ewald-package}
    \bibauthor{D.~Lindbo, L.~af~Klinteberg, D.~Saffar~Shamshirgar},
    \bibtitle{The spectral Ewald unified package},
    \textcolor{blue}{\url{http://github.com/ludvigak/SE_unified}},
    2018.
  \bibitem{lu19}
    \bibauthor{L.~Lu, M.~J.~Morse, A.~Rahimian, G.~Stadler,
    D.~Zorin},
    \bibarticle{Scalable Simulation of Realistic Volume Fraction
    Red Blood Cell Flows through Vascular Networks}, in
    \bibtitle{Proceedings of the International Conference for
    High Performance Computing, Networking, Storage and Analysis
    (SC '19)}, Denver, CO, USA, November 2019, article~6,
    \bibdoi{10.1145/3295500.3356203}.
  \bibitem{malhotra18}
    \bibauthor{D.~Malhotra, A.~Rahimian, D.~Zorin, G.~Biros},
    \bibarticle{A parallel algorithm for long-timescale
    simulation of concentrated vesicle suspensions in three
    dimensions}, 2018. Preprint available:
    \textcolor{blue}{\url{https://cims.nyu.edu/~malhotra/files/pubs/ves3d.pdf}}
  \bibitem{mittal18}
    \bibauthor{N.~Mittal, F.~Ansari, K.~Gowda.~V, C.~Brouzet,
    P.~Chen, P.~T.~Larsson, S.~V.~Roth, F.~Lundell, L.~Wågberg,
    N.~A.~Kotov, L.~D.~Söderberg},
    \bibarticle{Multiscale Control of Nanocellulose Assembly:
    Transferring Remarkable Nanoscale Fibril Mechanics to
    Macroscale Fibers},
    \bibtitle{ACS Nano} \bibvolume{12} \bibissue{7}, 6378--6388 (2018),
    \bibdoi{10.1021/acsnano.8b01084}.
  \bibitem{ojala15}
    \bibauthor{R.~Ojala, A.-K.~Tornberg},
    \bibarticle{An accurate integral equation method for
    simulating multi-phase Stokes flow},
    \bibtitle{J.\ Comp.\ Phys.}
    \bibvolume{298}, 145--160 (2015),
    \bibdoi{10.1016/j.jcp.2015.06.002}.
  \bibitem{dlmf}
    \bibauthor{F.~W.~J.~Olver, A.~B.~Olde~Daalhuis, D.~W.~Lozier,
    B.~I.~Schneider, R.~F.~Boisvert, C.~W.~Clark, B.~R.~Miller,
    B.~V.~Saunders, H.~S.~Cohl, M.~A.~McClain},
    \bibtitle{NIST Digital Library of Mathematical Functions},
    \textcolor{blue}{\url{http://dlmf.nist.gov/}}, Release 1.0.25 of 2019-12-15.
  \bibitem{palsson19}
    \bibauthor{S.~Pålsson, M.~Siegel, A.-K.~Tornberg},
    \bibarticle{Simulation and validation of surfactant-laden drops
    in two-dimensional Stokes flow}, \bibtitle{J.\ Comp.\ Phys.}
    \bibvolume{386}, 218--247 (2019),
    \bibdoi{10.1016/j.jcp.2018.12.044}.
  \bibitem{perez19}
    \bibauthor{C.~Pérez-Arancibia, L.~M.~Faria, C.~Turc},
    \bibarticle{Harmonic density interpolation methods for high-order
    evaluation of Laplace layer potentials in 2D and 3D},
    \bibtitle{J.\ Comp.\ Phys.}
    \bibvolume{376}, 411--434 (2019),
    \bibdoi{10.1016/j.jcp.2018.10.002}.
  \bibitem{power-miranda}
    \bibauthor{H.~Power, G.~Miranda}, \bibarticle{Second kind
    integral equation formulation of Stokes' flows past a
    particle of arbitrary shape}, \bibtitle{SIAM J.\ Appl.\
    Math.} \bibvolume{47} \bibissue{4}, 689--698 (1987),
    \bibdoi{10.1137/0147047}.
  \bibitem{pozrikidis92}
    \bibauthor{C.~Pozrikidis}, \bibtitle{Boundary integral and
    singularity methods for linearized viscous flow}, Cambridge
    University Press, Cambridge, 1992,
    \bibdoi{10.1017/CBO9780511624124}.
  \bibitem{rachh15}
    \bibauthor{M.~Rachh}, \bibtitle{Integral equation methods for
    problems in electrostatics, elastostatics and viscous flow},
    PhD~Thesis, New York University, 2015,
    ISBN~978-1-321-95485-2. Available:
    \textcolor{blue}{\url{https://search.proquest.com/docview/1710781501}}
  \bibitem{rachh17}
    \bibauthor{M.~Rachh, A.~Klöckner, M.~O'Neil},
    \bibarticle{Fast algorithms for Quadrature by Expansion I:
    Globally valid expansions}, \bibtitle{J.\ Comp.\ Phys.}
    \bibvolume{345}, 706--731 (2017),
    \bibdoi{10.1016/j.jcp.2017.04.062}.
  \bibitem{rahimian18}
    \bibauthor{A.~Rahimian, A.~Barnett, D.~Zorin},
    \bibarticle{Ubiquitous evaluation of layer potentials using
    Quadrature by Kernel-Independent Expansion}, \bibtitle{BIT Numer.\ Math.}
    \bibvolume{58}, 423--456 (2018),
    \bibdoi{10.1007/s10543-017-0689-2}.
  \bibitem{reddig13}
    \bibauthor{S.~Reddig, H.~Stark},
    \bibarticle{Nonlinear dynamics of spherical particles in
    Poiseuille flow under creeping-flow condition},
    \bibtitle{J.~Chem.~Phys.} \bibvolume{138} \bibissue{23},
    234902 (2013), \bibdoi{10.1063/1.4809989}.
  \bibitem{saad86}
    \bibauthor{Y.~Saad, M.~H.~Schultz},
    \bibarticle{GMRES: A generalized minimal residual algorithm
    for solving nonsymmetric linear systems}, \bibtitle{SIAM J.\
    Sci.\ Stat.\ Comput.} \bibvolume{7} \bibissue{3}, 856--869 (1986),
    \bibdoi{10.1137/0907058}.
  \bibitem{shamshirgar17}
    \bibauthor{D.~Saffar~Shamshirgar, A.-K.~Tornberg},
    \bibarticle{Fast Ewald summation for electrostatic potentials
    with arbitrary periodicity},
    \href{https://arxiv.org/abs/1712.04732}{\bibtitle{\textcolor{blue}{arXiv:1712.04732 [math.NA]}}}
    (2017).
  \bibitem{siegel18}
    \bibauthor{M.~Siegel, A.-K.~Tornberg},
    \bibarticle{A local target specific quadrature by expansion
    method for evaluation of layer potentials in 3D},
    \bibtitle{J.\ Comp.\ Phys.}
    \bibvolume{364}, 365--392 (2018),
    \bibdoi{10.1016/j.jcp.2018.03.006}.
  \bibitem{sorgentone18}
    \bibauthor{C.~Sorgentone, A.-K.~Tornberg},
    \bibarticle{A highly accurate boundary integral equation
    method for surfactant-laden drops in 3D},
    \bibtitle{J.\ Comp.\ Phys.}
    \bibvolume{360}, 167--191 (2018),
    \bibdoi{10.1016/j.jcp.2018.01.033}.
  \bibitem{squires05}
    \bibauthor{T.~M.~Squires, S.~R.~Quake},
    \bibarticle{Microfluidics: Fluid physics at the nanoliter
    scale}, \bibtitle{Rev.\ Mod.\ Phys.}
    \bibvolume{77} \bibissue{3}, 977--1026 (2005),
    \bibdoi{10.1103/RevModPhys.77.977}.
  \bibitem{srinivasan18}
    \bibauthor{S.~Srinivasan, A.-K.~Tornberg},
    \bibarticle{Fast Ewald summation for Green's functions of
    Stokes flow in a half-space},
    \bibtitle{Res.\ Math.\ Sci.} \bibvolume{5} \bibissue{3},
    article~35 (2018), \bibdoi{10.1007/s40687-018-0153-1}.
  \bibitem{tlupova13}
    \bibauthor{S.~Tlupova, J.~T.~Beale},
    \bibarticle{Nearly Singular Integrals in 3D Stokes Flow},
    \bibtitle{Commun.\ Comput.\ Phys.}
    \bibvolume{14} \bibissue{5}, 1207--1227 (2013),
    \bibdoi{10.4208/cicp.020812.080213a}.
  \bibitem{tlupova19}
    \bibauthor{S.~Tlupova, J.~T.~Beale},
    \bibarticle{Regularized single and double layer integrals in 3D Stokes flow},
    \bibtitle{J.\ Comp.\ Phys.}
    \bibvolume{386}, 568--584 (2019),
    \bibdoi{10.1016/j.jcp.2019.02.031}.
  \bibitem{tornberg08}
    \bibauthor{A.-K.~Tornberg, L.~Greengard},
    \bibarticle{A fast multipole method for the three-dimensional
    Stokes equations}, \bibtitle{J.\ Comp.\ Phys.}
    \bibvolume{227} \bibissue{3}, 1613--1619 (2008),
    \bibdoi{10.1016/j.jcp.2007.06.029}.
  \bibitem{trefethen14}
    \bibauthor{L.~N.~Trefethen, J.~A.~C.~Weideman},
    \bibarticle{The Exponentially Convergent Trapezoidal Rule},
    \bibtitle{SIAM Rev.} \bibvolume{56} \bibissue{3}, 385--458
    (2014),
    \bibdoi{10.1137/130932132}.
  \bibitem{ying06}
    \bibauthor{L.~Ying, G.~Biros, D.~Zorin},
    \bibarticle{A high-order 3D boundary integral equation solver
    for elliptic PDEs in smooth domains},
    \bibtitle{J.\ Comp.\ Phys.}
    \bibvolume{219}, 247--275 (2006),
    \bibdoi{10.1016/j.jcp.2006.03.021}.
  \bibitem{wala18}
    \bibauthor{M.~Wala, A.~Klöckner},
    \bibarticle{A fast algorithm with error bounds for Quadrature
    by Expansion},
    \bibtitle{J.\ Comp.\ Phys.}
    \bibvolume{374}, 135--162 (2018),
    \bibdoi{10.1016/j.jcp.2018.05.006}.
  \bibitem{wala19}
    \bibauthor{M.~Wala, A.~Klöckner},
    \bibarticle{A fast algorithm for Quadrature by Expansion in
    three dimensions},
    \bibtitle{J.\ Comp.\ Phys.}
    \bibvolume{388}, 655--689 (2019),
    \bibdoi{10.1016/j.jcp.2019.03.024}.
  \bibitem{wala20}
    \bibauthor{M.~Wala, A.~Klöckner},
    \bibarticle{Optimization of fast algorithms for global
    Quadrature by Expansion using target-specific expansions},
    \bibtitle{J.\ Comp.\ Phys.}
    \bibvolume{403}, 108976 (2020),
    \bibdoi{10.1016/j.jcp.2019.108976}.
  \bibitem{zhao10}
    \bibauthor{H.~Zhao, A.~H.~G.~Isfahani, L.~N.~Olson,
    J.~B.~Freund},
    \bibarticle{A spectral boundary integral method for flowing blood cells},
    \bibtitle{J.\ Comp.\ Phys.}
    \bibvolume{229}, 3726--3744 (2010),
    \bibdoi{10.1016/j.jcp.2010.01.024}.
\end{thebibliography}
